\documentclass{siamltex}
\usepackage{graphicx}
\usepackage{amsmath,amssymb}
\usepackage{bm}
\usepackage[latin1]{inputenc}
\usepackage{algorithm,algorithmic}
\usepackage{color}
\usepackage{subfigure}
\usepackage{url}  

\DeclareMathOperator{\range}{Range}   %

\newcommand{\cA}{\mathcal{A}}
\newcommand{\cB}{\mathcal{B}}
\newcommand{\cC}{\mathcal{C}}

\newcommand{\cF}{\mathcal{F}}

\newcommand{\cR}{\mathcal{R}}

\newcommand{\cZ}{\mathcal{Z}}

\newcommand{\RR}{\mathbb{R}}

\renewcommand{\a}{\alpha}
\renewcommand{\b}{\beta}
\newcommand{\g}{\gamma}

\newcommand{\<}{\langle}
\renewcommand{\>}{\rangle}

\newcommand{\x}{\times}

\newcommand{\tp}{^{\sf T}}

\newcommand{\tml}[3][]{\bm{\left(} #2 \bm{\right)}_{ #1} \bm{\cdot} #3}
\newcommand{\tmr}[3][]{#2 \bm{\cdot} \bm{\left(} #3 \bm{\right)}_{ #1}}

\usepackage{listings}
\definecolor{dkgreen}{rgb}{0,0.6,0}
\definecolor{gray}{rgb}{0.5,0.5,0.5}
\definecolor{mauve}{rgb}{0.58,0,0.82}
\title{Analyzing Large and Sparse Tensor Data using Spectral Low-Rank 
  Approximation. }
\author{Lars Eldén and Maryam Dehghan, \today}
\author{Lars Eld\'en\thanks{Department of Mathematics, Link\"{o}ping
    University,   Link\"{o}ping, Sweden
    ({lars.elden@liu.se})} \and Maryam Dehghan\thanks{Department of
    Teleinformatics Engineering, Federal university of Cear\'a,
    Fortaleza, Brazil ({maryamdehghan@ufc.br,
      maryamdehghan880@yahoo.com})} 
   }

\begin{document}

\maketitle

   \begin{center}
   \small  \today
   \end{center}

   \bigskip

\begin{abstract}
  Information is extracted from large and sparse data sets organized
  as 3-mode tensors. Two methods are described, based on best
  rank-(2,2,2) and rank-(2,2,1)   approximation of the tensor. The
first   method can be considered as a generalization of spectral graph
  partitioning to tensors, and it gives a reordering of the tensor that
  clusters the information. The second method gives an expansion of
  the tensor in  sparse rank-(2,2,1) terms, where the terms
  correspond to graphs. The low-rank approximations are computed using
  an efficient Krylov-Schur type algorithm that avoids filling in the sparse
  data. The methods are applied to topic search in news text, a tensor
  representing conference author-terms-years, and network traffic
  logs.  
\end{abstract}

 \begin{keywords}
   tensor, multilinear rank, best rank-(p,q,r) approximation,
    sparse, graph,  block  Krylov-Schur
   algorithm, (1,2)-symmetric tensor, text analysis, traffic logs,
   spectral partitioning
 \end{keywords}       

 \begin{AMS}
 05C50,   15A69, 65F15.
 \end{AMS}

\section{Introduction}
\label{sec:intro}

Finding clusters in sparse data is a standard task in numerous areas
of information sciences. Here we are concerned with data that are
organized according to three categories, represented by real-valued
tensors with three modes, i.e.  objects
$\cA \in \RR^{l \times m \times n}$. In particular  we will be dealing
with large, sparse tensors, which can 
be thought of as a collection of adjacency matrices for undirected
graphs. Spectral partitioning is an important class of methods for
finding interesting clustering structure in a graph, and extracting
information from the data set represented by the graph. In this paper
we present a method for analyzing tensors that can be considered as a
generalization of spectral graph partitioning.

The proposed method is based on the computation of the best
rank-(2,2,2) approximation of the tensor. This is
the analogy  with the computation of two eigenvalues and corresponding
eigenvectors in spectral graph partitioning. By a reordering of certain
vectors and applying the same reordering of the tensor, a clustering
of the information in the tensor is performed. A variant of the
method is also described, where an expansion of the tensor in
rank-(2,2,1) terms is computed. 

\begin{figure}[htbp!]    
\centering
\includegraphics[width=.35\textwidth]{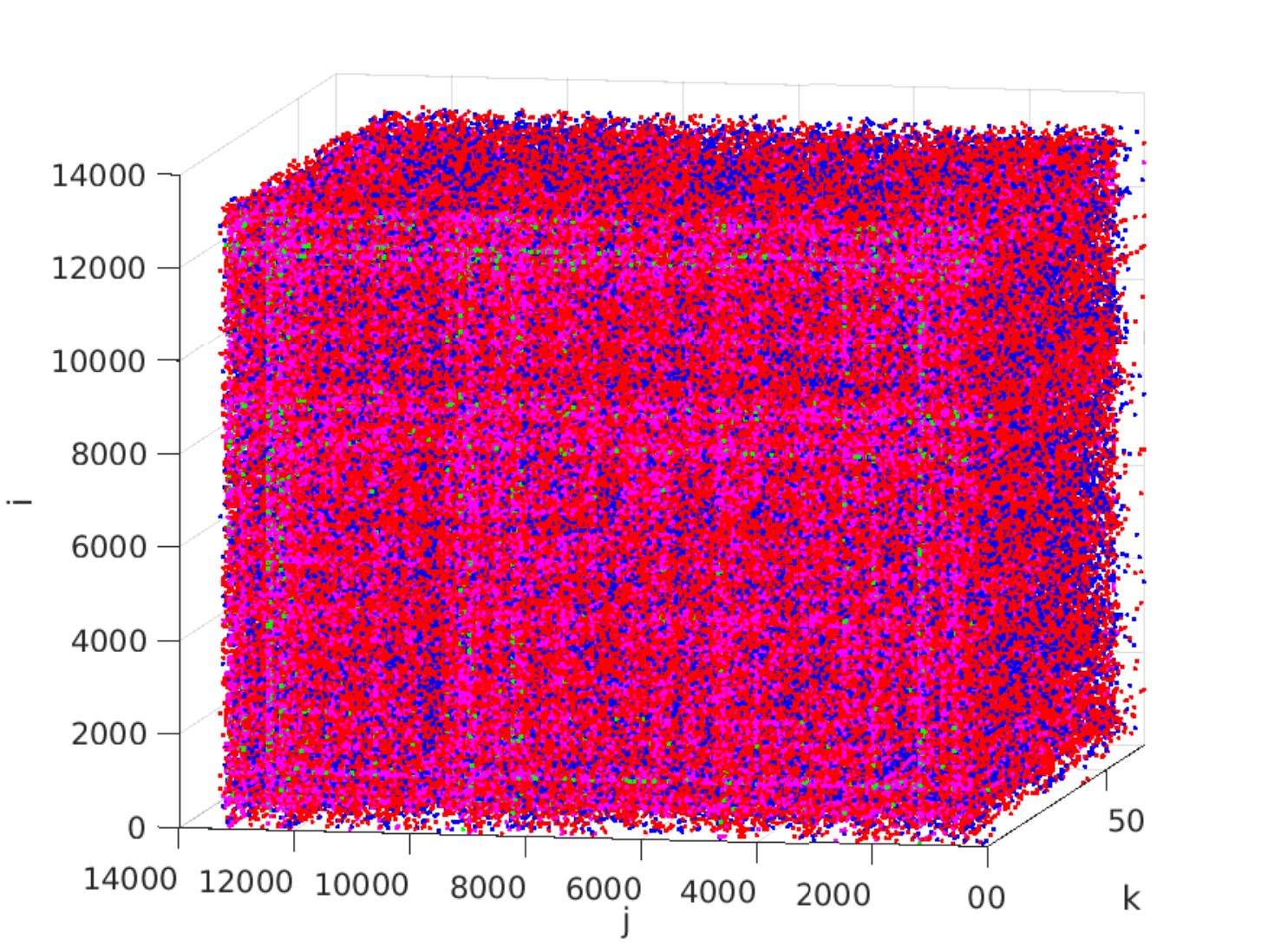}
\includegraphics[width=.35\textwidth]{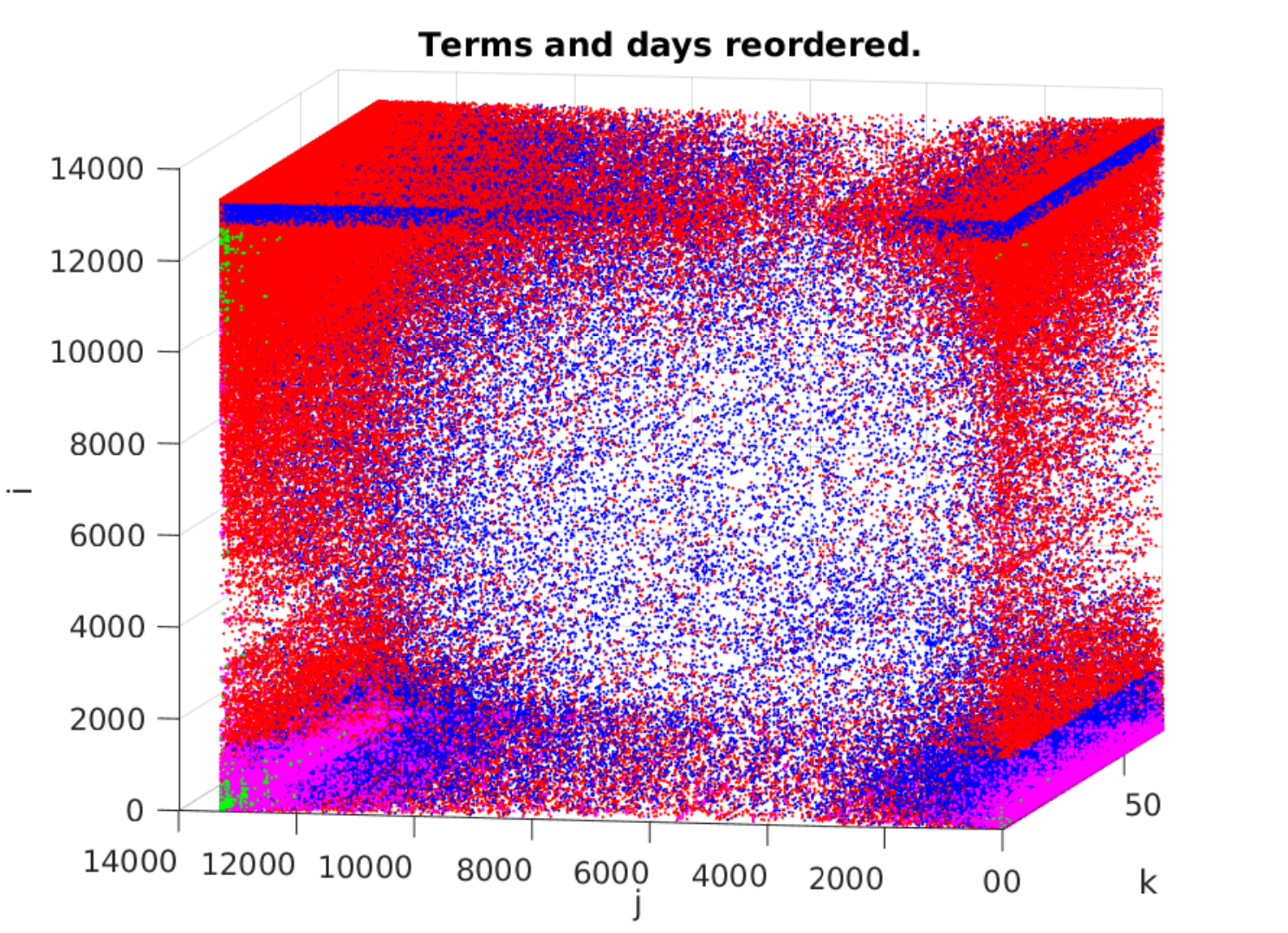}
\caption{The   original Reuters tensor (left) and reordered by our method
  (right).  The two main topics are placed  at the upper left and bottom
  right corners of the right plot. Each dot represents the
  cooccurrence of two terms in news texts from one day.  }
\label{fig:reuters1-2} 
\end{figure}

This paper is an investigation of the usefulness of the method in
three application areas with large and sparse tensors. The first area
is topic search in a collection of news texts from Reuters over 66
days. From this corpus of more than 13000 terms and their cooccurrence
data organized in 66 graphs, we find the two main topics and their
subtopics, and also analyze the topic variation in time. The
clustering of the data tensor by the method is illustrated in Figure
\ref{fig:reuters1-2}. 

The second application is a tensor from a series of
conferences, where the tensor modes represent   authors, terms, and
years. Our methods find the topics and their variation over the
years. It also exhibits the development in time of participating
authors.

The third data set is a tensor of network traffic logs over 371 time
intervals, and the objective is to find the dominating senders
(spammers?) and receivers out of almost 9000. We compute an
approximation of the tensor in terms of an expansion of three
rank-(2,1,1) terms, through which the dominating communication
patterns are identified and visualized. In our fourth example we
compute an expansion of seven rank-(2,2,1) terms of the above news
text tensor. Each term in the expansion is a graph that represents a
dominating topic and its prevalence during the time period.

This paper is the third in a series of three. In the first
\cite{eldehg20b} we develop an algorithm for computing the best
low-rank approximation of large and sparse tensors. The second paper
\cite{eldehg20a} shows that the method used in the present paper can
be considered as a generalization of spectral graph partitioning.  The
main contribution of the present paper is the demonstration that the
best low-rank approximation can be used to analyze large and sparse
tensor data coming from real applications.

More general purpose tensor clustering algorithms are described in
\cite{jsb09,ldjd10,cwhyl13,wzzy19} and applied to synthetic examples
and small to medium size examples with real data. Larger problems are
solved in \cite{cxz20,xzxgg20}.

This paper is organized as follows. We  introduce  some pertinent tensor
concepts and give a very brief outline of our method in Section
\ref{sec:concepts}. Then in Section \ref{sec:topic} we give more
detailed description of the method of computing a best rank-(2,2,2)
approximation as applied to the problem finding
topics in news texts. The conference paper example is given in Section 
\ref{sec:NeurIPS}. The ideas behind the expansion of a tensor in
rank-(2,2,1) terms are described in Section \ref{sec:expansion}, and
the application to the network traffic logs and topics in news texts
in Section \ref{sec:Reuters221}.  The expansion algorithm is
summarized in Section \ref{sec:sum-comp}. Finally we give some
concluding remarks in Section \ref{sec:conclusions}. 

\section{Tensor Concepts}
\label{sec:concepts}

\subsection{Notation and Preliminaries}
Tensors will be denoted by calligraphic letters, e.g
$\mathcal{A},\mathcal{B}$, matrices by capital roman letters and
vectors by lower case roman letters.  In order not to burden the
presentation with too much detail, we sometimes will not explicitly
mention the dimensions of matrices and tensors, and assume that they
are such that the operations are well-defined. The whole presentation
will be in terms of tensors of order three, or equivalently 3-tensors.
The generalization to order-$N$ tensors is obvious.

We will use the term tensor for a 
3-dimensional array of real numbers, $\cA \in \RR^{l \times m \times
  n}$, where the vector space is equipped with some algebraic
structures to be defined. The different ``dimensions'' of the tensor
are referred to as \emph{modes}.  We will use both standard subscripts
and ``MATLAB-like'' notation: a particular tensor element will be
denoted in two equivalent ways, $\cA(i,j,k) = a_{ijk}$.

A subtensor obtained by fixing one of the indices is called a
\emph{slice}, e.g., $\cA(i,:,:)$.
In the case when the index is fixed in the first mode, we call the slice 
a $1-$slice, and correspondingly for the other modes.  
A  slice  can be  considered as an order-3 tensor
(3-tensor) with a singleton mode, and also as a matrix.  
A tensor $\cA$ is called \emph{(1,2)-symmetric} if all 3-slices
$\cA(:,:,k)$  are symmetric.
A \emph{fiber} is a subtensor, where all indices but one are
fixed. For instance, $\cA(i,:,k)$ denotes a mode-2 fiber. 

For a given third order tensor, there are three associated subspaces,
one for each mode. These subspaces are given by  
\begin{align*}
& \range \{ \cA(:,j,k) \; | \; j = 1:m, \; k = 1:n  \},\\
& \range \{ \cA(i,:,k) \; | \; i = 1:l, \; k = 1:n  \},\\
& \range \{ \cA(i,j,:) \; | \; i = 1:l, \; j = 1:m  \}.
\end{align*}
The \emph{multilinear rank} \cite{hit:27,sili08} of the tensor is said
to be equal to $(p,q,r)$ 
if the dimension of these subspaces are $p,$ $q,$ and $r$, respectively. 



\subsection{Tensor-Matrix Multiplication}
\label{sec:ten-matmult}

We define \emph{multilinear multiplication of a tensor by a matrix} as
follows.  For concreteness we first present multiplication by one matrix
along the first mode and later for all three modes simultaneously. The
mode-$1$ product of 
a tensor $\cA \in \RR^{l \times m \times n}$ by a matrix $U \in \RR^{p
  \times l}$ is defined\footnote{The notation \eqref{eq:contra}-\eqref{eq:mat-tensor} was
suggested by de Silva and Lim \cite{sili08}. An alternative notation was
  earlier given in \cite{lmv:00a}.  Our $\tml[d]{X}{\cA}$ is the same
  as $\cA \times_d X$ in the latter system.} 
\begin{equation}\label{eq:contra}
 \RR^{p \times m \times n} \ni \mathcal{B} =  \tml[1]{U}{ \cA} ,\qquad
b_{ijk} = \sum_{\a=1}^{l} u_{i \a} a_{ \a jk}  .
\end{equation}
This means that all mode-$1$ fibers in the $3$-tensor $\cA$ are
multiplied by the matrix $U$. Similarly, mode-$2$ multiplication by a
matrix $V \in \RR^{q \x m}$ means that all mode-$2$ fibers are
multiplied by the matrix $V$.  Mode-$3$ multiplication is analogous.
With a third matrix $W\in \RR^{r \x n}$, the tensor-matrix
multiplication
in all modes is given by
\begin{equation}
  \label{eq:mat-tensor}
\RR^{p \x q \x r}	\ni \cB = \tml{U,V,W}{\cA}, \qquad b_{ijk} = \sum_{\a,\b,\g=1}^{l,m,n} u_{i \a} v_{j \b} w_{k \g}a_{ \a \b \g},
\end{equation}
where the mode of each multiplication is understood from the order in
which the matrices are given.

It is convenient to introduce a separate notation for multiplication
by a transposed matrix $\bar{U} \in \RR^{l \times p}$:
\begin{equation}
  \label{eq:cov}
\RR^{p \times m \times n} \ni \mathcal{C}= \tml[1]{\bar{U}\tp }{\cA}
=\tmr[1]{\cA}{\bar{U}}, \qquad
c_{ijk} = \sum_{\a=1}^{l} a_{\a jk} \bar{u}_{\a i}.
\end{equation}
Let $u \in \RR^l $ be a  vector and $\cA \in \RR^{l \times m
  \times n}$ a tensor. Then
\begin{equation}
  \label{eq:ident-1}
 \RR^{1 \times m   \times n} \ni \cB := \tml[1]{u\tp}{\cA} =
\tmr[1]{\cA}{u} \equiv B \in \RR^{m   \times  n}.
\end{equation}
Thus we identify a tensor with a singleton dimension with  a
matrix.

\subsection{Inner Product and Norm}
\label{sec:norm}

Given two tensors $\cA$ and $\cB$ of the same dimensions, we define
the \emph{inner product},
\begin{equation}\label{eq:inner-prod}
\< \cA , \cB \> = \sum_{\a,\b,\g} a_{\a \b \g} b_{\a \b \g},
\end{equation}
and the analogous for 2-tensors. 
The  corresponding \emph{tensor norm} is
\begin{equation}
  \label{eq:norm}
  \| \cA \| = \< \cA, \cA \>^{1/2}.
\end{equation}
This \emph{Frobenius norm} will be used throughout the paper.
As in the matrix case, the norm is invariant under orthogonal
transformations, i.e.
\[
  \| \cA \| = \left\| \tml{U,V,W}{\cA} \right\| = \| \tmr{\cA}{P,Q,S} \|,
\]
for orthogonal matrices $U$, $V$, $W$, $P$, $Q$, and $S$. This is obvious 
from the fact that multilinear multiplication by orthogonal matrices
does not change the Euclidean length of the corresponding fibers of the tensor.

\subsection{Best Rank-$(r_1,r_2,r_3)$  Approximation}
\label{sec:best-appr}
The problem of approximating  tensor $\cA \in \RR^{l \x m \x n}$ by
another tensor $\cB$ of lower multi-linear rank, 
\begin{equation}
  \label{eq:best-appr-min}
\min_{\rank(\cB)=(r_1,r_2,r_3)} \| \cA - \cB \|^2,  
\end{equation}
is treated in \cite{lmv:00b,zhgo:00,elsa09,idav08}. It is shown in
\cite{lmv:00b} that \eqref{eq:best-appr-min} is equivalent to 
\begin{equation}
  \label{eq:best-appr-max}
  \max_{X,Y,Z} \| \tmr{\cA}{X,Y,Z} \|^2, \quad \mbox{subject to} \quad
  X\tp X = I_{r_1}, \quad
  Y\tp Y =I_{r_2}, \quad Z\tp Z = I_{r_3},
\end{equation}
where $X \in \RR^{l \x r_1}$, $Y \in \RR^{m \x r_2}$, and 
$Z \in \RR^{n \x r_3}$.  The problem \eqref{eq:best-appr-max} can be
considered as a \emph{Rayleigh quotient maximization problem}, in
analogy with the matrix case \cite{amsd:02}.
To
simplify the terminology somewhat we will refer to a solution
$(U,V,W)$ of the approximation problem as \emph{the best
  rank-$(r_1,r_2,r_3)$ approximation of the tensor}; it is
straightforward to  prove \cite{lmv:00b} 
that  the corresponding core tensor is $\cF =
\tmr{\cA}{U,V,W}$, so, strictly speaking, $\cB=\tml{U,V,W}{\cF}$ is
the best approximation.
In the case when the tensor is (1,2)-symmetric, the solution has the
form $(U,U,W)$, cf. \cite{eldehg20a}. 

The constrained maximization problem \eqref{eq:best-appr-max} can be
thought of as an unconstrained maximization problem on the
\emph{Grassmann manifold}, cf. \cite{eas98,ams:07,elsa09}. 

\subsection{Computing the Best Low-Rank Approximation}
\label{sec:compute}
The simplest algorithm for computing the best low-rank approximation
of a tensor is the HOOI method, which  an alternating orthogonal
iteration algorithm 
\cite{lmv:00b}. 
For small and medium-size tensors Grassmann variants of standard
optimization algorithms often converge faster
\cite{idav08,elsa09,sali10,iahd11}. For large and sparse tensors a
Block Krylov-Schur-like (BKS) method has been developed \cite{eldehg20b}. It is
analogous to an algorithm \cite{lesoya:98}, which is a standard method
for computing a low-rank approximation of matrices.  The method is
memory efficient as it uses the tensor only in tensor-matrix
multiplications, where the matrix has a small number of
columns. Convergence is often quite fast for tensors with data from
applications. It is shown in \cite{eldehg20b} that the BKS method is more
robust and faster than HOOI for large and sparse tensors, and we use
it throughout this paper.

The method developed in \cite{eldehg20b} is intended for
(1,2)-symmetric tensors. The method can be modified for non-symmetric
tensors $\cA$ by applying (implicitly) the symmetric BKS method to the
(1,2)-symmetric tensor
\begin{equation}
  \label{eq:0AA0}
  \begin{pmatrix}
    0 & \cA \\
    \cA' & 0
  \end{pmatrix},
\end{equation}
where the 3-slices of $\cA'$ are the transposes of the corresponding
slices of $\cA$.

\subsection{Partitioning Graphs by Low-Rank
  Approximation}
\label{sec:low-rank}

\begin{figure}[htbp!]    
\centering
\includegraphics[width=.4\textwidth]{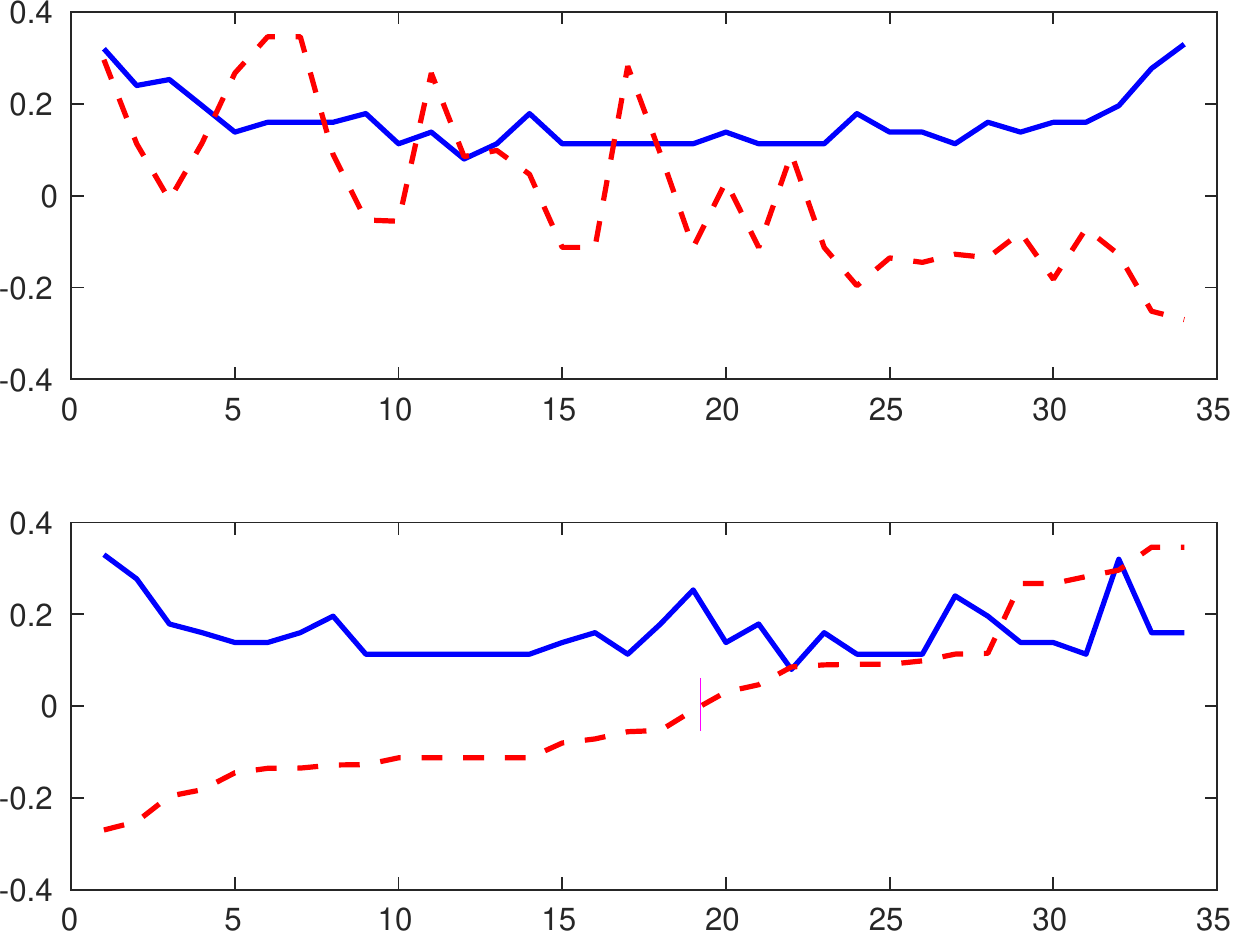}\quad
\includegraphics[width=0.4\textwidth,height=5cm]{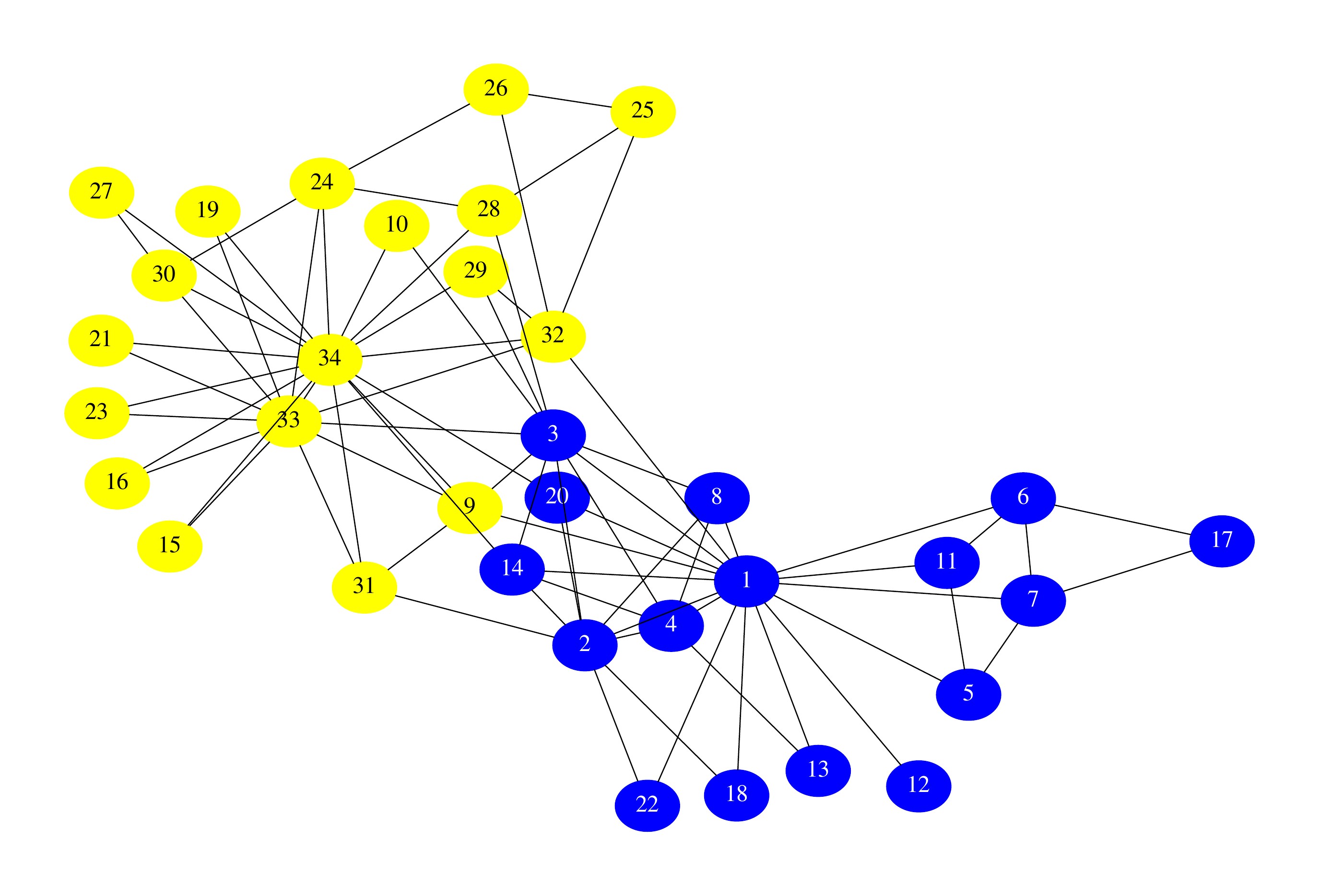}
\caption{The first two eigenvectors  $u_1$ (blue, solid) and $u_2$
  (red, dashed) of the normalized adjacency matrix of the Karate
  club graph, original (top left) and reordered (bottom left). The magenta
  vertical line shows  the position, where the elements of reordered
  vector $u_2$   changes sign. Partitioning the reordered vertices in
  the neighborhood of the sign change gives a low cost in terms of the
  number of broken edges, see \cite[Chapter 10]{eldenmatrix19}. The
  right panel shows the partitioned graph (colorcoded). 
}  
\label{fig:karate}
\end{figure}

Spectral partitioning is a standard method for partitioning undirected
graphs that are close to being disconnected, see
e.g. \cite{psl90,chung97,lux07,rine14}. Often it is based on the computation
of the two smallest eigenvalues of the normalized graph Laplacian. We here
give a brief account of the equivalent approach (see
e.g. \cite[Chapter 10]{eldenmatrix19}), when the two largest
eigenvalues of the normalized adjacency matrix are used. Let the
adjacency matrix be $A$ and the degree matrix,
$D=\diag(d)=\diag(Ae),$
where $e$ is a vector of all ones. The normalized adjacency matrix
\begin{equation}
  \label{eq:normalize}
A_N = D^{-1/2} A D^{-1/2}  
\end{equation}
has largest eigenvalue equal to 1. If the graph is connected, then the
second largest eigenvalue is strictly smaller than 1, and the
corresponding eigenvector, the Fiedler vector, can be used to
partition the graph into two subgraphs so that the cost (essentially
in terms of the number of broken edges) is small. The eigenvectors
$u_1$ and $u_2$ of the Karate club graph \cite{zach77}, see also
\cite[Chapter 10]{eldenmatrix19}, are illustrated in Figure
\ref{fig:karate}.

In \cite{eldehg20a} spectral graph partitioning is generalized to
(1,2)-symmetric tensors, where each 3-slice can be considered as the
normalized adjacency matrix of an undirected graph. A concept of
disconnectedness of tensors is defined, and it is shown that
(1,2)-symmetric tensors that are close to being disconnected can be
partitioned so that the partitioning cost is small. The partitioning
is derived from the best rank-(2,2,2) approximation $(U,U,W)$ of the
tensor. The partitioning procedure is analogous to the matrix case,
and we describe it in connection with the example in Section
\ref{sec:topic}. For further details, see \cite{eldehg20a}.

\section{Topic search in News Text}
\label{sec:topic}

  The data for this example are described in \cite{bamr03}. We
  cite from the description:\footnote{\url{http://vlado.fmf.uni-lj.si/pub/networks/data/CRA/terror.htm}.}
  \begin{quote}
``The \emph{Reuters terror news network}  is based on all stories released
during 66 consecutive days by the news agency Reuters concerning the
September 11 attack on the U.S., beginning at 9:00 AM EST 9/11/01. The
vertices of a network are words (terms); there is an edge between two
words iff they appear in the same text unit (sentence). The weight of
an edge is its frequency. The network has n = 13332 vertices
(different words in the news) and m = 243447 edges, 50859 with value
larger than 1. There are no loops in the network.''
\end{quote} 
We organized the data in a tensor
$\cA \in \RR^{13332 \times 13332 \times 66}$, where each 3-slice was
normalized as in \eqref{eq:normalize}. Thus $a_{ijk} >0$ if terms $i$
and $j$ occur in the same sentence in the news during day $k$;
otherwise $a_{ijk}=0$. Each 3-slice is the normalized adjacency matrix
of an undirected graph, and the tensor is (1,2)-symmetric.
The tensor is illustrated to the left in Figure
\ref{fig:reuters1-2}. Each nonzero element is shown as a colored
dot.
%
The horizontal and vertical 
  directions correspond to  terms, and the lateral
  directions to days.  The color
  coding is from blue for the largest elements, via red and magenta to
green for the smallest non-zero elements. The same coding is used in
all similar figures. Due to the fact that several non-zeros
  can be hidden behind others, the coding is more informative when the
tensor is   relatively small or    has some larger scale
structure (as in the following figures). The code for producing the
plot is a modification of \texttt{spy3} from
Tensorlab\footnote{\texttt{https://www.tensorlab.net/}}. 
%
Obviously it is impossible to discern any useful structure from the image
of the tensor.

\begin{figure}[htbp!]    
\centering
\includegraphics[width=.6\textwidth]{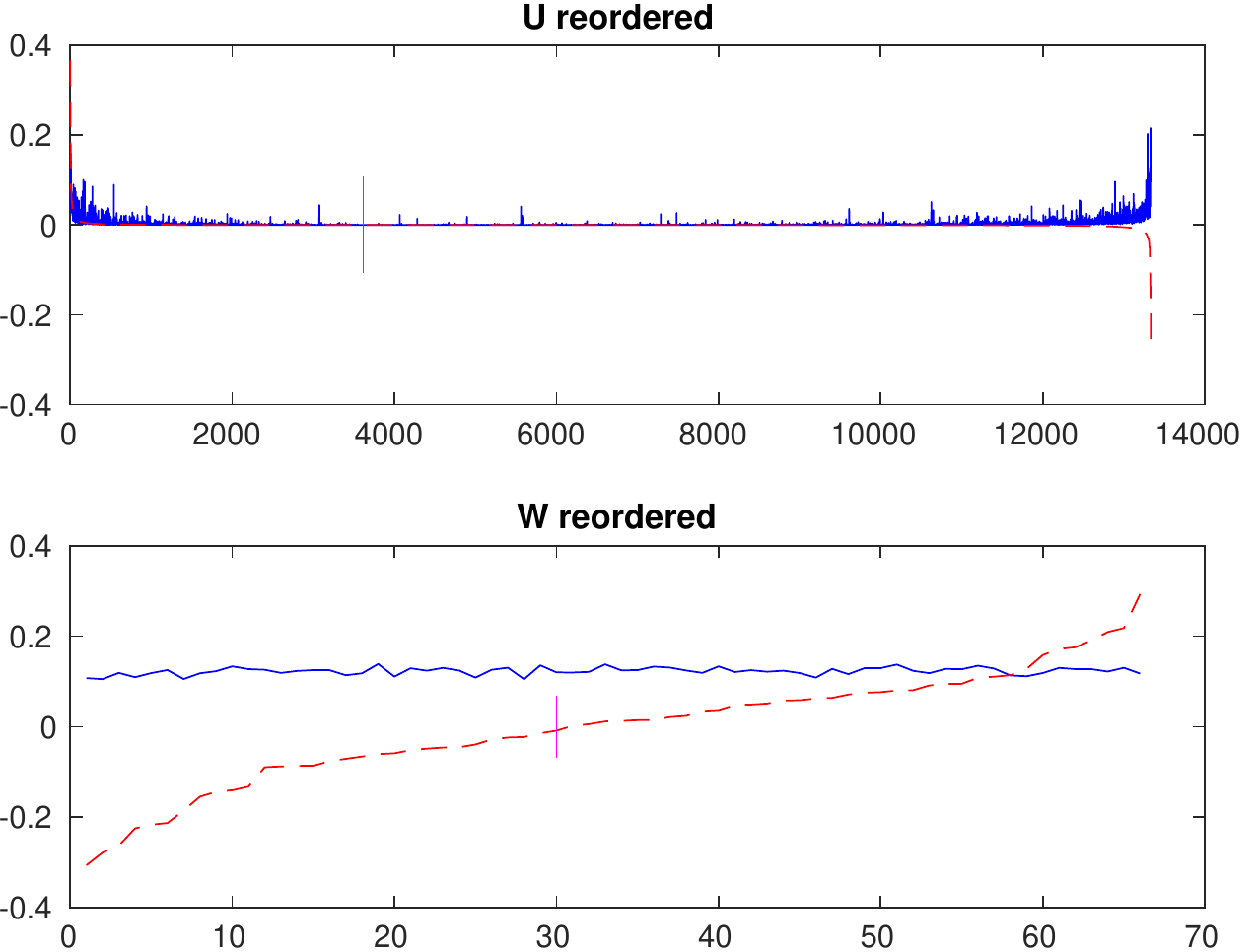}
\caption{ The column vectors of $U$ and $W$ after
  reordering; $u_1$ and $w_1$ are blue and the others red. The magenta
  vertical lines show  the position, where the elements of reordered
  vectors $u_2$  and $w_2$   change sign.  }  
\label{fig:UW66reordered}
\end{figure}

In spectral graph partitioning (Section \ref{sec:low-rank}) we
compute the two largest eigenvalues and corresponding eigenvectors of
the normalized adjacency matrix, which is equivalent to computing the
best rank-2 approximation of the matrix. The analogous concept for a
(1,2)-symmetric tensor is  the best rank-(2,2,2) approximation of the
tensor, which we computed using the BKS method\footnote{The execution
  time in Matlab on a standard desk top computer was about 44
  seconds.}.  Since the  
problem is (1,2)-symmetric, we have $U=V$. The elements of the column
vectors of $U$, which we  call term vectors, are reordered so that
those of $U(:,2)$ become monotonic; we refer to the reordered term vectors as
$(u_1,u_2)$. Similarly the elements of the  column vectors of $W$
are reordered, giving the   day vectors $(w_1,w_2)$.  The vectors for
the Reuters tensor are
illustrated in Figure \ref{fig:UW66reordered}.
A closer look (using a logarithmic plot) at the middle elements, 9000
say, of $u_1$ reveals that they are at least a couple of magnitudes
smaller than the elements at the beginning and the end.  This
indicates that the middle elements are insignificant in the mode-(1,2)
partitioning of the tensor; analogously, the corresponding terms are
insignificant in the search for topics.   In Table \ref{tab:terms66} we
list the 25 first and 
last reordered terms, as well as 25 close to the sign change of $u_2$. 
\begin{table}[htbp!]
  \centering
  \caption{Beginning terms, middle terms and end terms after
    reordering. The table is organized so that the most significant
    terms appear at the top of the first column and the bottom of the
    third (to conform with the image of the reordered tensor in Figure
    \ref{fig:reuters66reordered}). The middle column shows a few terms
    at the middle, which are relatively insignificant with respect to
    the topics in the other two columns.\label{tab:terms66} }
  \begin{tabular}{ll|ll|llll}
    &$T^1$& &Insignificant & &$T^2$&\\
    \hline
     &world\_trade\_ctr&    &planaria&         &hand&  \\      
    &pentagon&           &mars&             &ruler&       \\
    &new\_york&           &i.d&              &saudi&     \\  
    &attack&             &eyebrow&          &guest&    \\   
    &hijack&             &calendar&         &mullah&   \\ 
    &plane&              &belly&            &movement&   \\ 
    &airliner&           &auction&          &camp&   \\
    &tower&              &attentiveness&    &harbor& \\   
    &twin&               &astronaut&        &shelter&  \\   
    &washington&         &ant&              &group&      \\ 
    &sept&               &alligator&        &organization&\\
    &suicide&            &adjustment&       &leader&      \\
    &pennsylvania&       &cnd&              &rule&   \\
    &people&             &nazarbayev&       &exile&    \\   
    &passenger&          &margrit&          &guerrilla&  \\ 
    &jet&                &bailes&           &afghanistan& \\
    &mayor&              &newscast&         &fugitive&    \\
    &110-story&          &turgan-tiube&     &dissident&   \\
    &mayor\_giuliani&     &N52000&           &taliban&     \\
    &hijacker&           &ishaq&            &islamic&     \\
    &commercial&         &sauce&            &network&     \\
    &aircraft&           &sequential&       &al\_quaeda&   \\
    &assault&            &boehlert&         &militant&    \\
    &airplane&           &tomato&           &saudi-born&  \\
    &miss&               &resound&          &bin\_laden&    
  \end{tabular}
\end{table}
It is apparent from the table that we have two different topics, which
we will refer to as $T^1$ and $T^2$.

\begin{figure}[htbp!]    
\centering
\includegraphics[width=.8\textwidth]{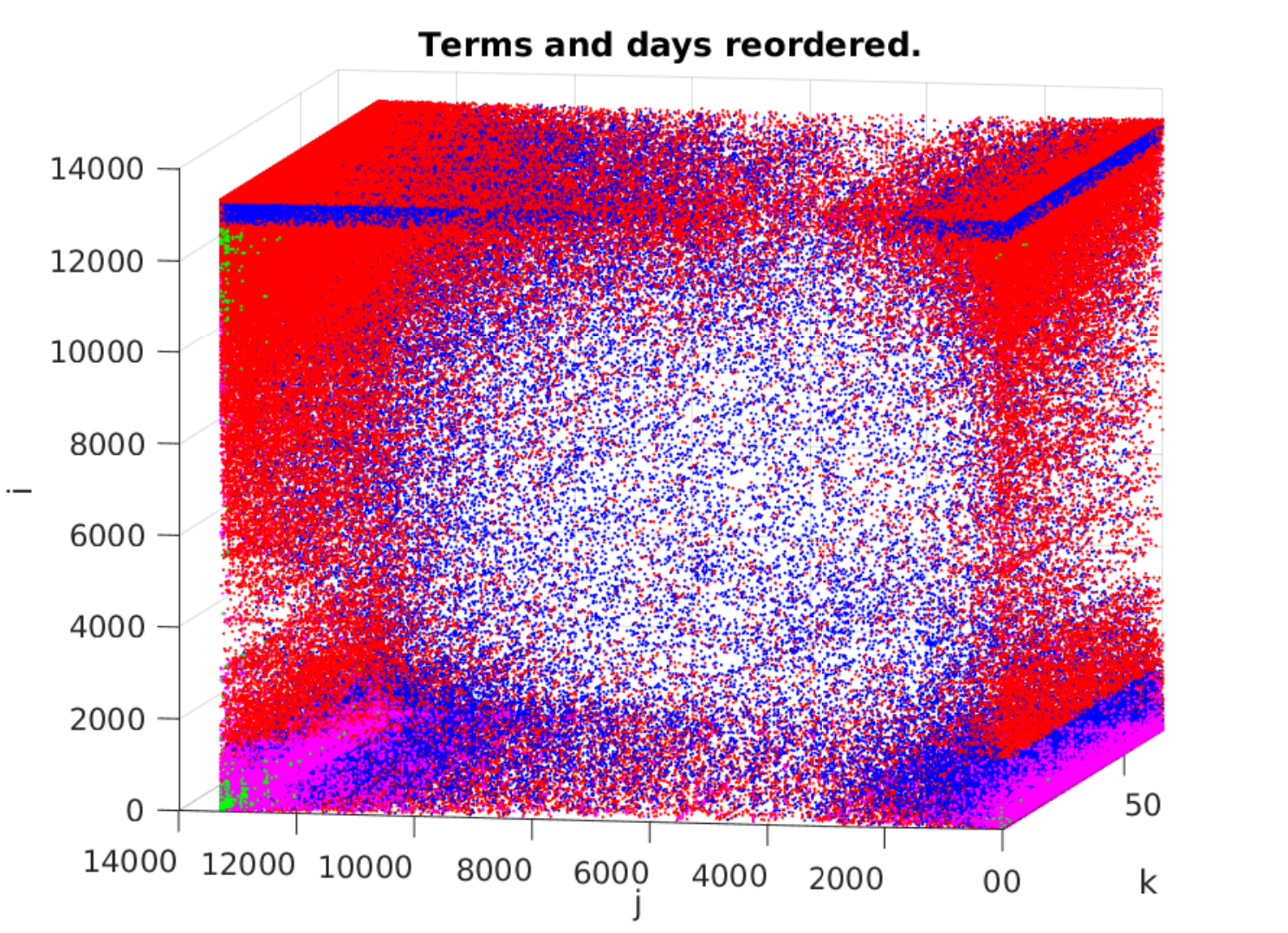}
\caption{The tensor reordered in the term and day modes. The top left terms
  correspond to those in  the first column in Table \ref{tab:terms66}, and the
  bottom right terms to  those in the third column of Table
  \ref{tab:terms66}.
It is interesting to note that many of the insignificant terms at the
middle of the tensor in the figure are marked blue. This indicates
that even if those terms are insignificant they are rather frequent.
}
\label{fig:reuters66reordered}
\end{figure}

Applying the same reordering to the original tensor, we get the tensor
illustrated in Figure \ref{fig:reuters66reordered}. 
%
It is seen that the ``mass'' is concentrated to the four corners in
term modes, which correspond to Topics 1 and 2 in Table
\ref{tab:terms66}. We computed the norms of the subtensors of
dimensions $1000 \x 1000 \x 66$ at the corners, and they were $50.6$,
$25.5$, $25.5$, 
and $35.1$, respectively. As the norm of $\cA$ was $128.6$, 15\% of
the terms accounted for about 55\%
of the total mass.

The relatively dense cluster at the top left in Figure
\ref{fig:reuters66reordered} shows that the terms from topic $T^2$
cooccur to a great extent with other terms from $T^2$. Looking at the
bottom right cluster, the analogous
statement can be made about topic $T^1$. The cluster at the top right
in the tensor indicates that many terms from topic $T^1$ cooccur with
those in $T^2$.  Apart from this, it is not possible to discern in
Figure \ref{fig:reuters66reordered} any interesting structure.

To check that the middle terms are insignificant, we removed
approximately 9000 middle terms (after the reordering), and computed a
best rank-(2,2,2) approximation of  a tensor of dimension $4333 \times
4333 \times 66$. The analysis gave exactly the same most significant
terms as in the first and third columns of Table \ref{tab:terms66}. 

We divided the 66 days into two groups, where the first, $D_A$ (30
days), corresponds to the (original) indices of $w_2$ to the left of
the sign change, see Figure \ref{fig:UW66reordered}. The second group,
$D_B$ (36 days), corresponds to the  (original) indices to the
right.  It turns out that the reordered days are not consecutive:
\medskip\medskip 
{\small
\begin{verbatim}
   DA: XXXXXXXXXX X  XX X  XX XX    X XX           XXX  X X X X   X X    
   DB:           X XX  X XX  X  XXXX X  XXXXXXXXXXX   XX X X X XXX X XXXX 
\end{verbatim}
}
\medskip\medskip
 
\begin{figure}[htbp!]    
\centering
\includegraphics[width=.45\textwidth]{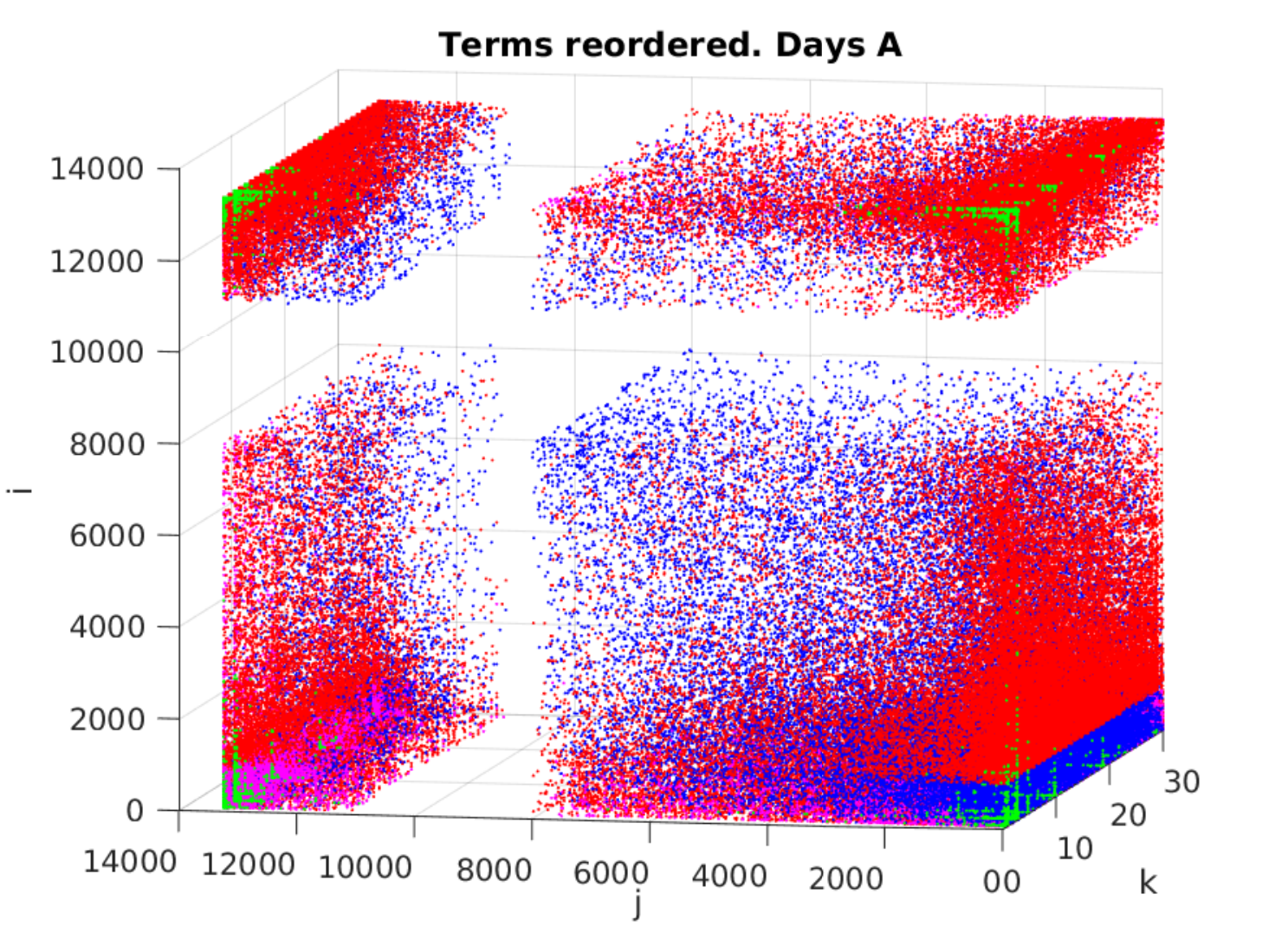}
\includegraphics[width=.45\textwidth]{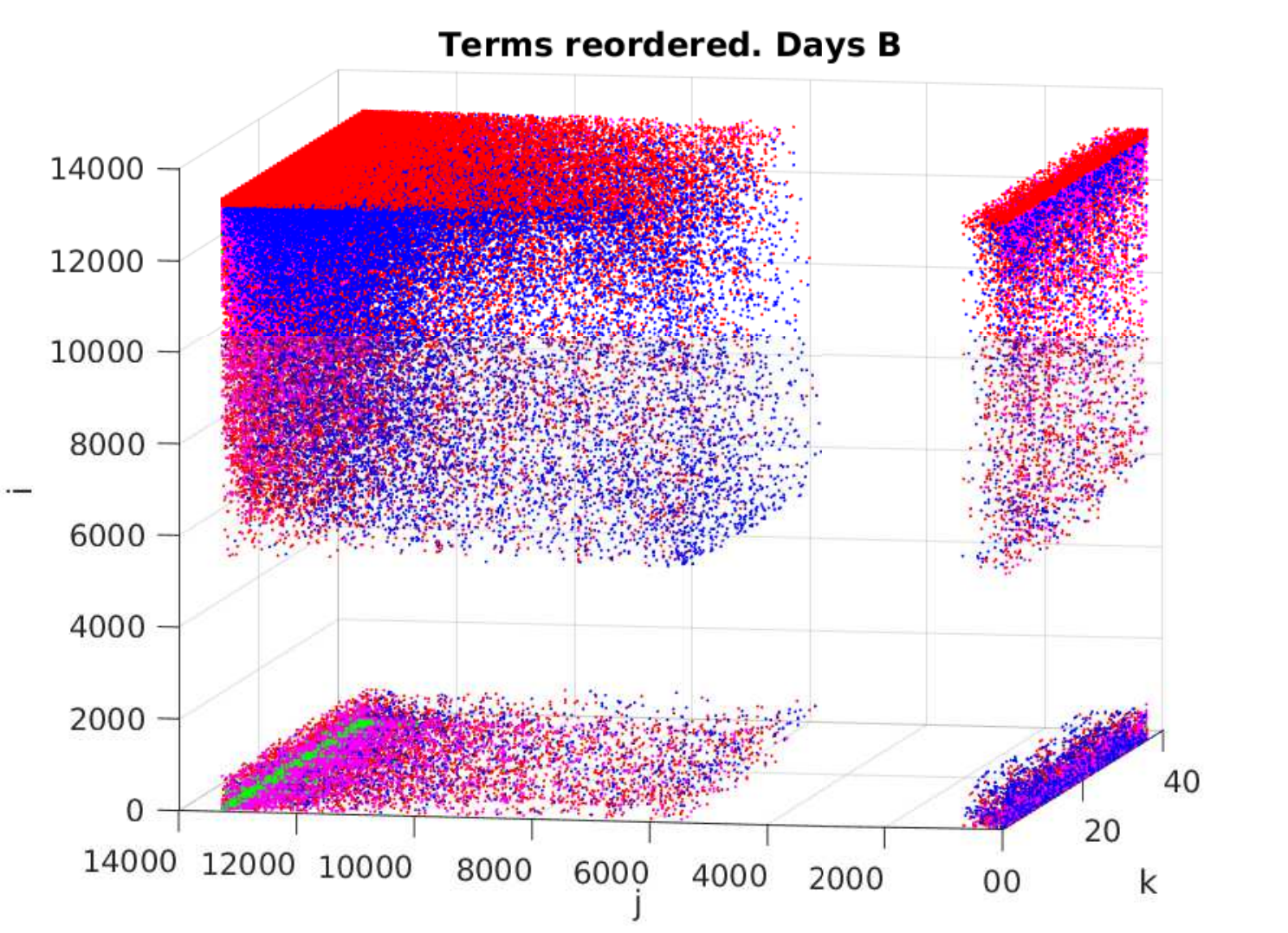}
\caption{Top row: The tensors for $D_A$ and $D_B$  reordered
  separately in the term mode. }
\label{fig:reutersDADBreordered}
\end{figure}

To see if there is any difference in vocabulary during the groups of
days, $D_A$ and $D_B$,  we analyzed the
groups  separately, by computing the best rank-(2,2,2)
approximations of the tensors $\cA(:,:,D_A) \in \RR^{13332 \times
  13332 \times 30}$  and $\cA(:,:,D_B) \in \RR^{13332 \times
  13332 \times 36}$. The
reordered tensors are shown in Figure \ref{fig:reutersDADBreordered},
and the most significant terms  in Table \ref{tab:termsDADB}.
\begin{table}[htbp!]
  \centering
  \caption{The topics of the two groups of days. The most significant
    terms appear at the top of the left two columns and at the
    bottom of the right two.\label{tab:termsDADB} }
  \begin{tabular}{llll|lllll}
    &$T^{1}_{A}$& &$T^{1}_{B}$ & &$T^2_A$& &$T^2_B$\\
    \hline
 &world\_trade\_ctr&    &attack&             &country&          &hand&        \\  
    &tower&              &world\_trade\_ctr&    &muslim&           &millionaire&   \\
    &pentagon&           &new\_york&           &force&            &pakistan-based&\\
    &twin&               &pentagon&           &united\_states&    &leader&        \\
    &new\_york&           &hijack&             &ruler&            &movement&      \\
    &plane&              &united\_states&      &reuter&           &follower&      \\
    &hijack&             &anthrax&            &law&              &saudi&         \\
    &airliner&           &tell&               &mullah&           &guest&         \\
    &attack&             &plane&              &pakistan&         &group&         \\
    &110-story&          &people&             &network&          &harbor&        \\
    &washington&         &airliner&           &dissident&        &train&         \\
    &pennsylvania&       &washington&         &terrorism&        &camp&          \\
    &commercial&         &reuter&             &group&            &organization&  \\
    &passenger&          &official&           &military&         &taliban&       \\
    &jet&                &sept&               &tell&             &exile&         \\
    &mayor&              &reporter&           &war&              &shelter&       \\
    &suicide&            &city&               &al\_quaeda&        &guerrilla&     \\
    &hijacker&           &news&               &leader&           &dissident&     \\
    &mayor\_giuliani&     &suicide&            &islamic&          &islamic&       \\
    &ruin&               &worker&             &militant&         &fugitive&      \\
    &rural&              &security&           &rule&             &network&       \\
    &crash&              &office&             &afghanistan&      &al\_quaeda&     \\
    &rubble&             &day&                &saudi-born&       &militant&      \\
    &people&             &tower&              &taliban&          &saudi-born&    \\
    &flight&             &letter&             &bin\_laden&        &bin\_laden&       \end{tabular}
\end{table}
From Table \ref{tab:termsDADB} we see that the topics are similar 
 for days $D_A$ and $D_B$, but the terms used differ considerably. In
 addition, Figure 
\ref{fig:reutersDADBreordered} shows that the vocabulary used for
topic $T^1$ is much more concentrated during $D_A$ than during
$D_B$. The opposite is  true of topic $T^2$.
It is also seen that the number of
insignificant terms that are not used at all, or used very
infrequently, during $D_A$ and $D_B$ are about 3200 and 4800,
respectively.

It is possible to apply  the procedure recursively: As an example, we
took the subtensor of dimension $7333 \times 7333 \times 36$  at 
the top left of the right image in Figure
\ref{fig:reutersDADBreordered}, 
 that contains mainly
terms from topic $T^1$ during days $D_B$, and computed its best
rank-(2,2,2) approximation. This produced the two subtopics given in
Table \ref{tab:termsT1DB}.

\begin{table}[htbp!] 
  \centering
  \caption{Two subtopics of topic $T^1_B$.  The most significant terms appear at
    the top of the first column and at the bottom of the fourth. It
    is seen that the two left columns are somewhat  more ``politically''
    oriented, while the two right columns are concerned with the actual attack.  \label{tab:termsT1DB}} 
  \begin{tabular}{llll|llllllll}
        &$T^{11}_B$& && &$T^{12}_B$\\
    \hline
 &tell&                 &special&          &jet&             &tower&  \\        
    &reuter&               &foreign&          &people&          &september&\\     
    &reporter&             &war&              &airplane&        &washington&     \\
    &force&                &news&             &month&           &suicide&        \\
    &northern\_alliance&    &military&         &wake&            &plane&          \\
    &opposition&           &security&         &assault&         &sept&           \\
    &rumsfeld&             &defense&          &passenger&       &airliner&       \\
    &official&             &pakistan&         &11&              &new\_york&       \\
    &sec&                  &united\_states&    &devastate&       &pentagon&       \\
    &pres\_bush&            &conference&       &hijacker&
                        &hijack&         
    \\
    &kabul&                &terrorism&        &pennsylvania&    &attack&         \\
    &minister&             &government&       &twin&            &world\_trade\_ctr&\\
  \end{tabular}
\end{table}
It is necessary to ask the question if the tensor method  gives more
information than an analogous spectral matrix/graph analysis. To
investigate this, we
computed a matrix by adding all the 3-slices of the tensor
and then normalizing as in \eqref{eq:normalize}. This is the
normalized adjacency matrix  of the graph that has an edge between
two terms if they occur in the same sentence any time during the 66
days, weighted by the number of occurrences. Spectral partitioning
using the reordered eigenvector corresponding to the second largest
eigenvalue gives mixed terms from topics $T^1$ and $T^2$  at one end, and
seemingly random, insignificant terms at the other. Thus the matrix
approach does not 
separate the terms from topics $T^1$  and $T^2$.
We repeated the experiment, where we removed about 8000  insignificant
terms. The results were very similar to those with
the full data set.

\section{NeurIPS Conference Papers}
\label{sec:NeurIPS}

\begin{figure}[htbp!]    
\centering
\includegraphics[width=.45\textwidth]{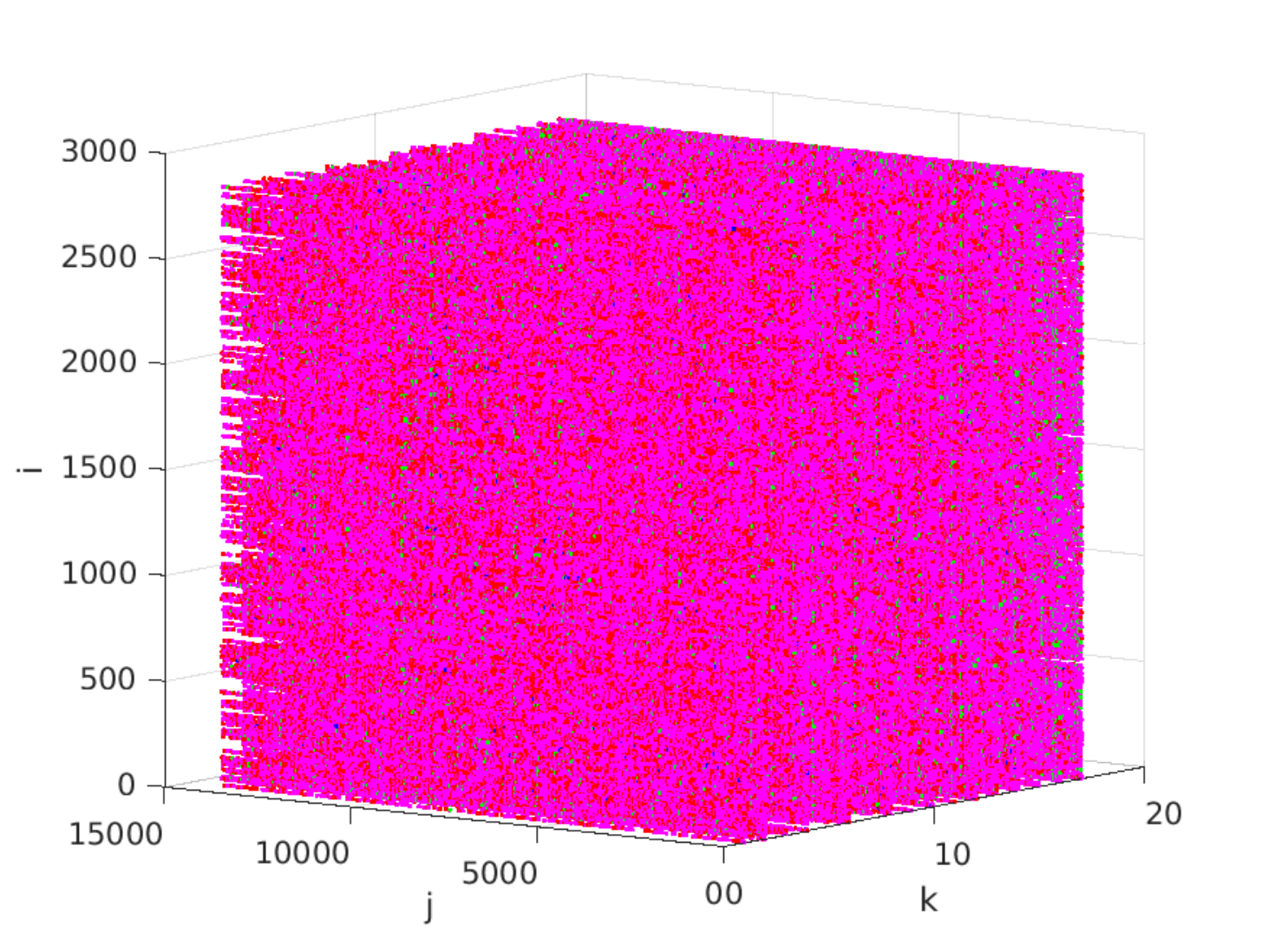}
\includegraphics[width=.45\textwidth]{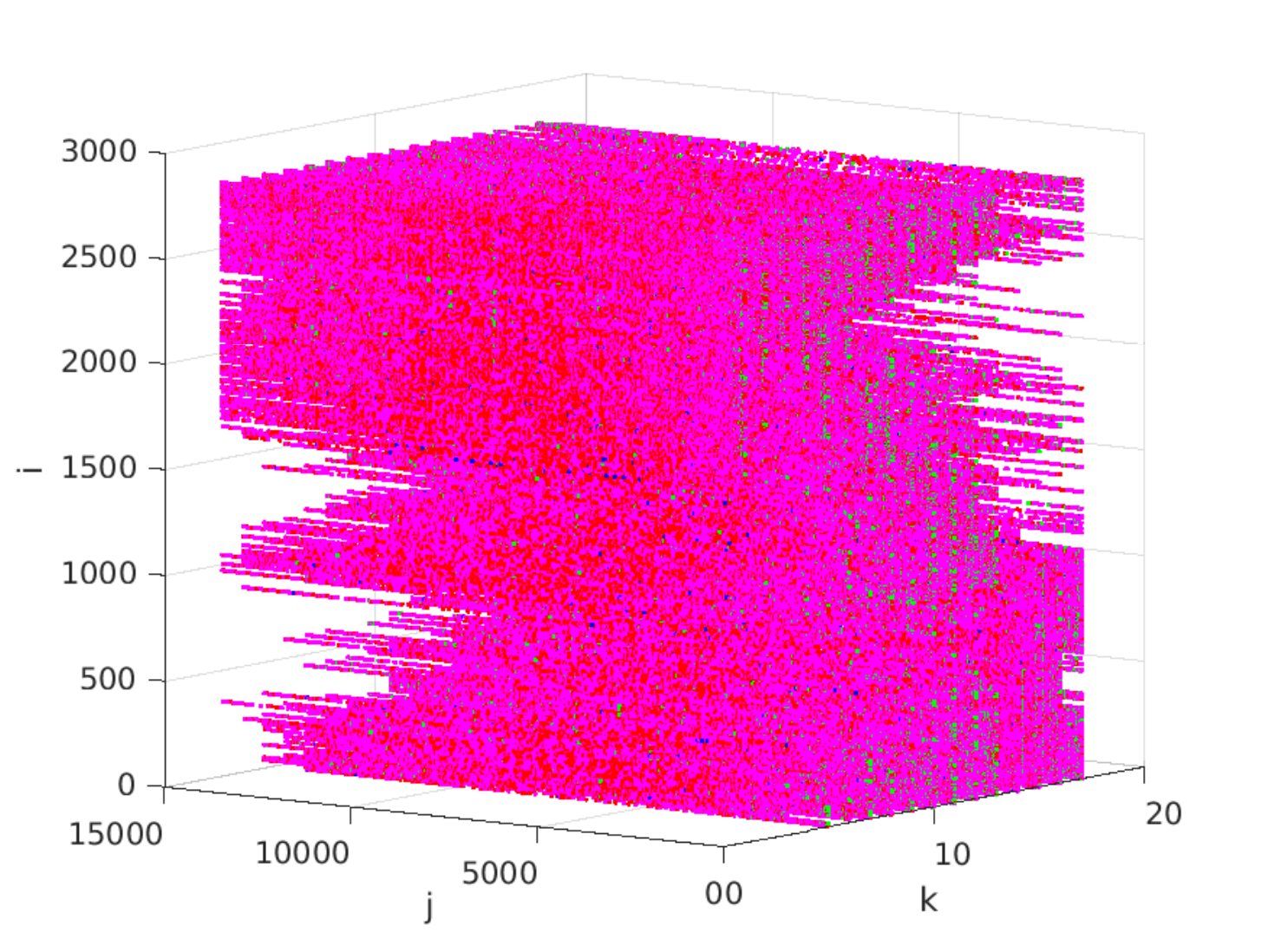}
\caption{Left: Original NeurIPS tensor. The author mode is vertical, the
 term mode is horizontal, and the year mode is lateral. The same
  color coding as in the previous section is used. Right: Reordered
  NeurIPS tensor. The tensor is reordered in the author and term modes. }
\label{fig:OrigNIPS} 
\end{figure}

Experiments with data from all the papers at  the Neural Information
Processing Systems Conferences\footnote{The acronym for the
  conferences used to be NIPS, but in 2018 it was  changed to
  NeurIPS. We will use the latter.} 1987-2003 are described in
\cite{gcpt07}.  We downloaded the data from \texttt{http://frostt.io/}
\cite{frosttdataset}, and formed a sparse tensor of dimension
$2862 \times 14036 \times 17$, where the modes represent
(author,terms,year), and the values are term counts. We performed a
non-symmetric normalization of  the 3-slices of the tensor (which is
equivalent to the symmetric normalization \eqref{eq:normalize} applied
to the tensor \eqref{eq:0AA0}). The tensor is
illustrated in Figure \ref{fig:OrigNIPS}.

\begin{figure}[htbp!]    
\centering
\includegraphics[width=.6\textwidth]{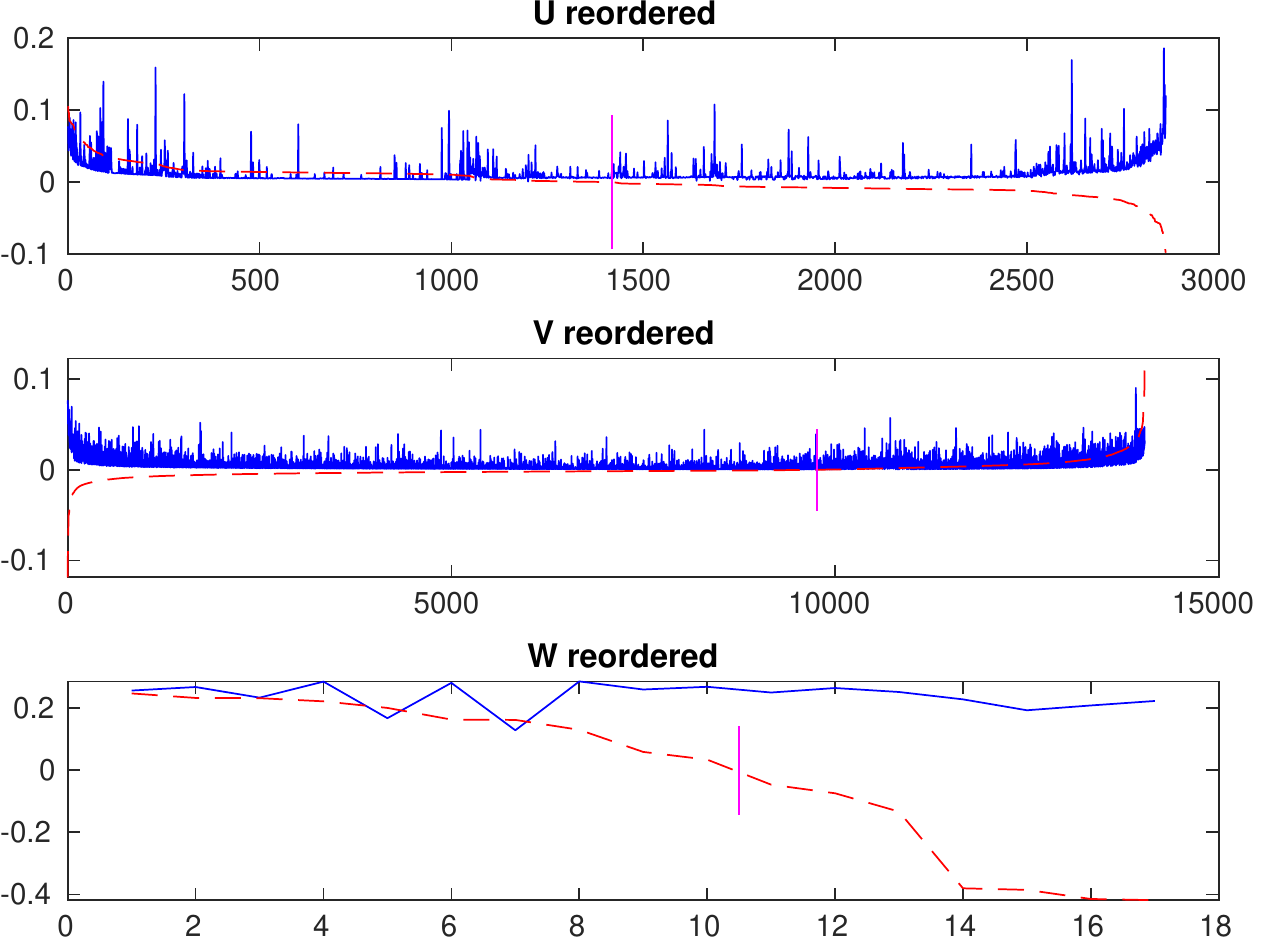}
\caption{Reordered vectors in the best rank-(2,2,2) approximation of
  the NeurIPS tensor. }
\label{fig:UVWreord1-17}
\end{figure}

We computed the best rank-(2,2,2) approximation of the tensor and
reordered the vectors as in the example in Section \ref{sec:topic},
see Figure \ref{fig:UVWreord1-17}.  The  tensor reordered in the
author and term modes 
 is
illustrated in Figure \ref{fig:OrigNIPS}. Some structure  can be
seen, but not very clearly. 
 The  reordered term vector suggests  that the two main topics in 
 the conferences during years 1-17 are those given in Table
 \ref{tab:NIPStopics1-17}.
\begin{table}[htbp!] 
  \centering
  \caption{Two main topics at NeurIPS conferences during years
    1-17. The terms are ordered so that the most important are at the
    top of column 1 and at the bottom of column 4.   \label{tab:NIPStopics1-17}}   
  \begin{tabular}{llll|llllllll}
        &$T^1$& && &$T^2$\\
    \hline
    &network&        &error&        &flies&           &neurobiology&\\
    &learning&       &figure&       &individuals&     &jerusalem&   \\
    &training&       &layer&        &independence&    &van&         \\
    &input&          &weights&      &princeton&       &steveninck&  \\
    &units&          &time&         &engineering&     &bialek&      \\
    &networks&       &set&          &huji&            &william&     \\
    &output&         &algorithm&    &israel&          &code&        \\
    &model&          &function&     &hebrew&          &rob&         \\
    &state&          &hidden&       &center&          &tishby&      \\
    &weight&         &signal&       &ruyter&          &brenner&     \\
    &performance&    &image&        &school&          &naftali&     \\
    &recognition&    &unit&         &nec&             &universality&\\
  \end{tabular}
\end{table}

From the spy-plot of the tensor in Figure \ref{fig:OrigNIPS} (right),
we see structure in the  time-author modes.  The reordering of the
year vector is 
\begin{verbatim}
           4 6 3 5 2 7 1 8 9 10 11 12 13 14 17 15 16
\end{verbatim}
This indicates (as expected) that   the structure differs between the
first and the last years of the period. To investigate
that we analyzed separately the tensors for years 1-9, and 10-17.  

\begin{figure}[htbp!]    
\centering
\includegraphics[width=.6\textwidth]{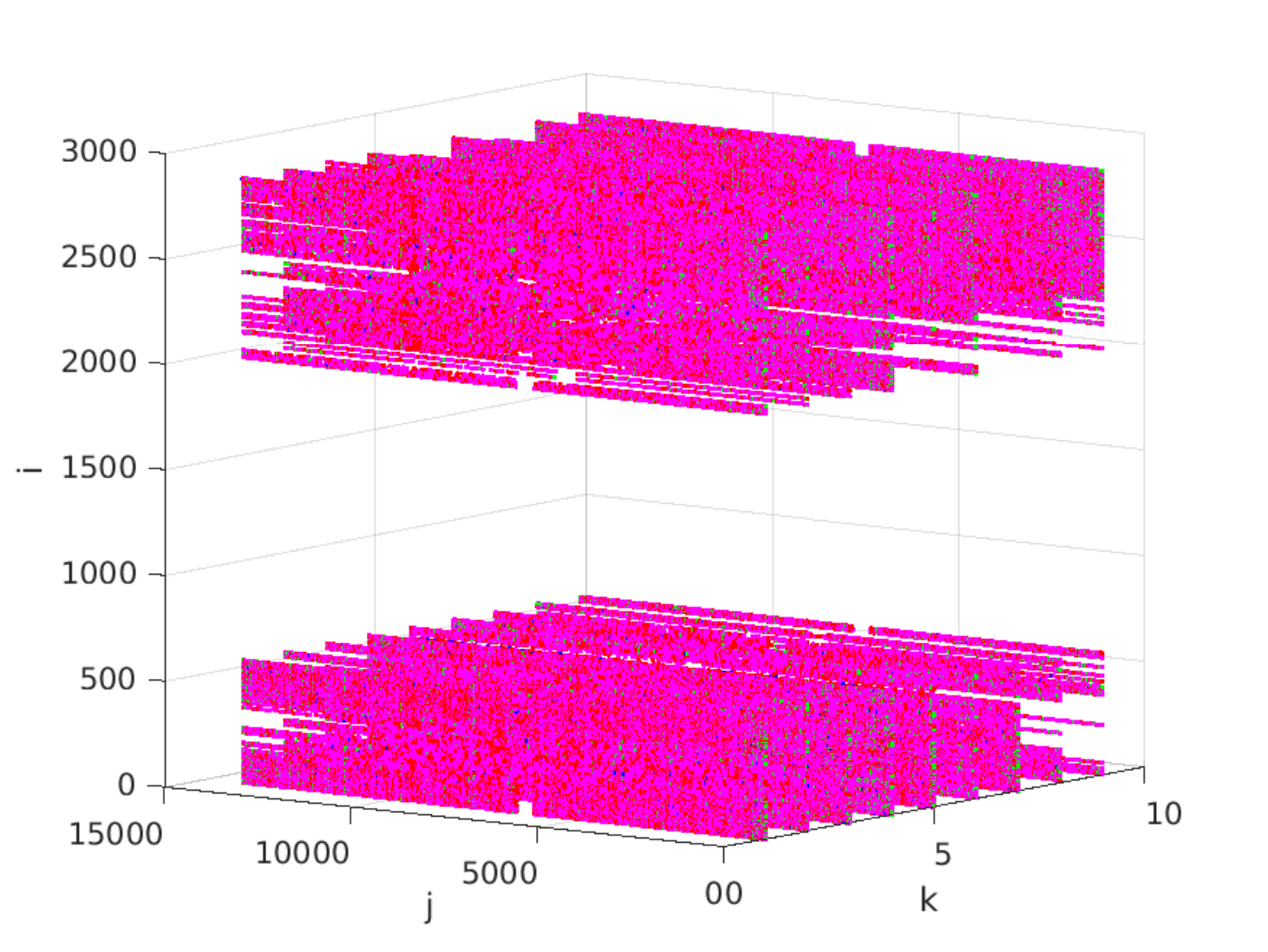}
\caption{The  tensor for the years 1-9, reordered in the author and
  term modes. }
\label{fig:A-NIPS-reord1-9}
\end{figure}

In Figure \ref{fig:A-NIPS-reord1-9} we illustrate the reordered tensor
for years 1-9.  
The gap in the middle of Figure \ref{fig:A-NIPS-reord1-9} corresponds
to authors who were not present during years 1-9. The main topics are
given in Table \ref{tab:NIPStopics1-9}. 
\begin{table}[htbp!] 
  \centering
  \caption{Two main topics at NeurIPS conferences during years
    1-9 (denoted A).    \label{tab:NIPStopics1-9}}   
  \begin{tabular}{llll|llllllll}
        &$T^1_A$& && &$T^2_A$\\
    \hline
     &propagation&    &neuron&         &schraudolph&    &mixtures&\\
    &units&          &hopfield&       &eye&            &gamma&   \\
    &back&           &patterns&       &jabri&          &gesture& \\
    &connections&    &pattern&        &reward&         &policy&  \\
    &network&        &layer&          &leen&           &singh&   \\
    &analog&         &chips&          &pca&            &mixture& \\
    &chip&           &fig&            &density&        &jordan&  \\
    &connection&     &hidden&         &experts&        &model&   \\
    &synapse&        &classifiers&    &validation&     &tresp&   \\
    &circuit&        &parallel&       &prediction&     &missing& \\
    &bower&          &vowel&          &data&           &em&     \\ 
    &associative&    &input&          &tangent&        &dayan&
  \end{tabular}
\end{table}

\begin{figure}[htbp!]    
\centering 
\includegraphics[width=.45\textwidth]{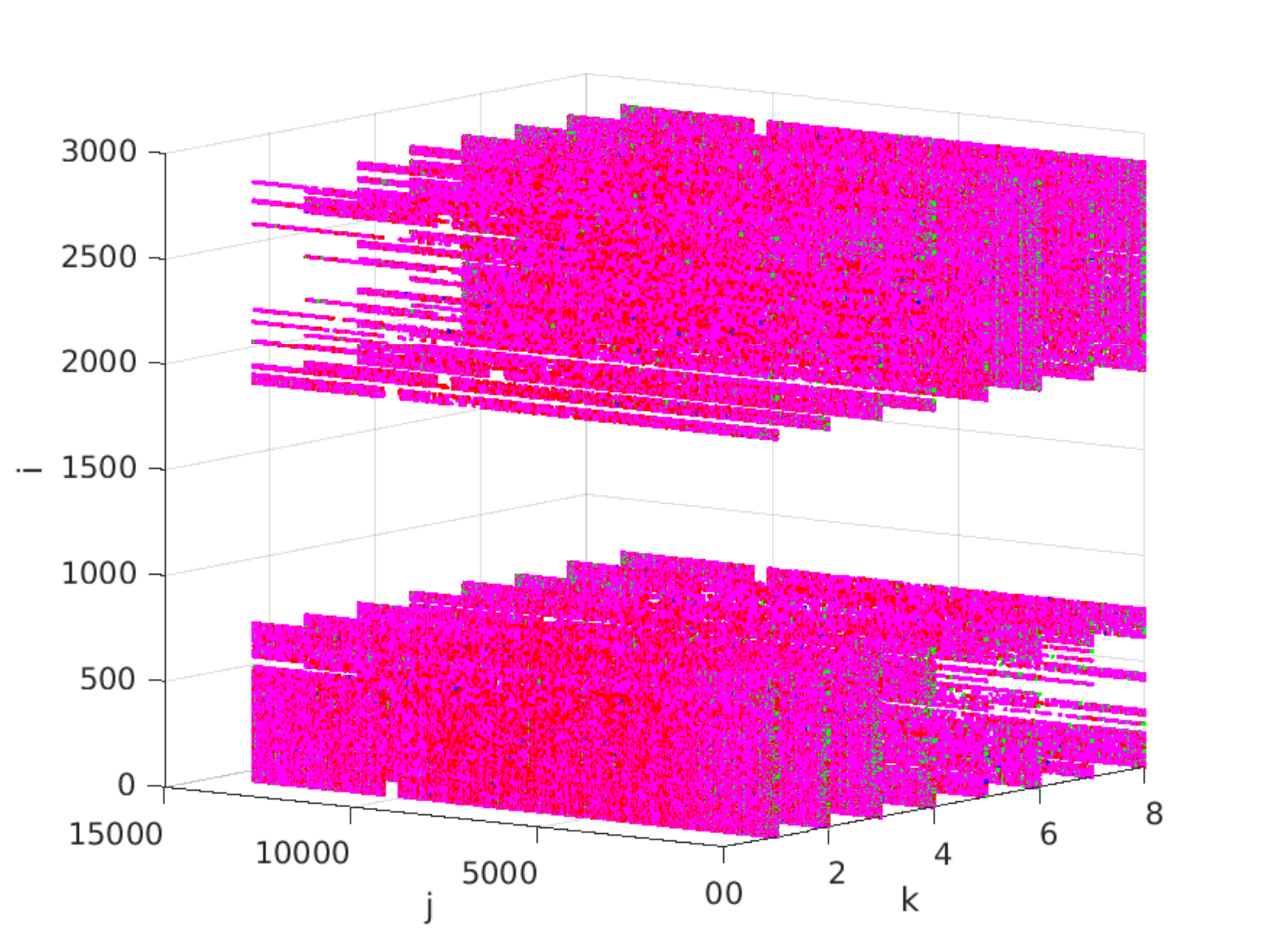}
\includegraphics[width=.45\textwidth]{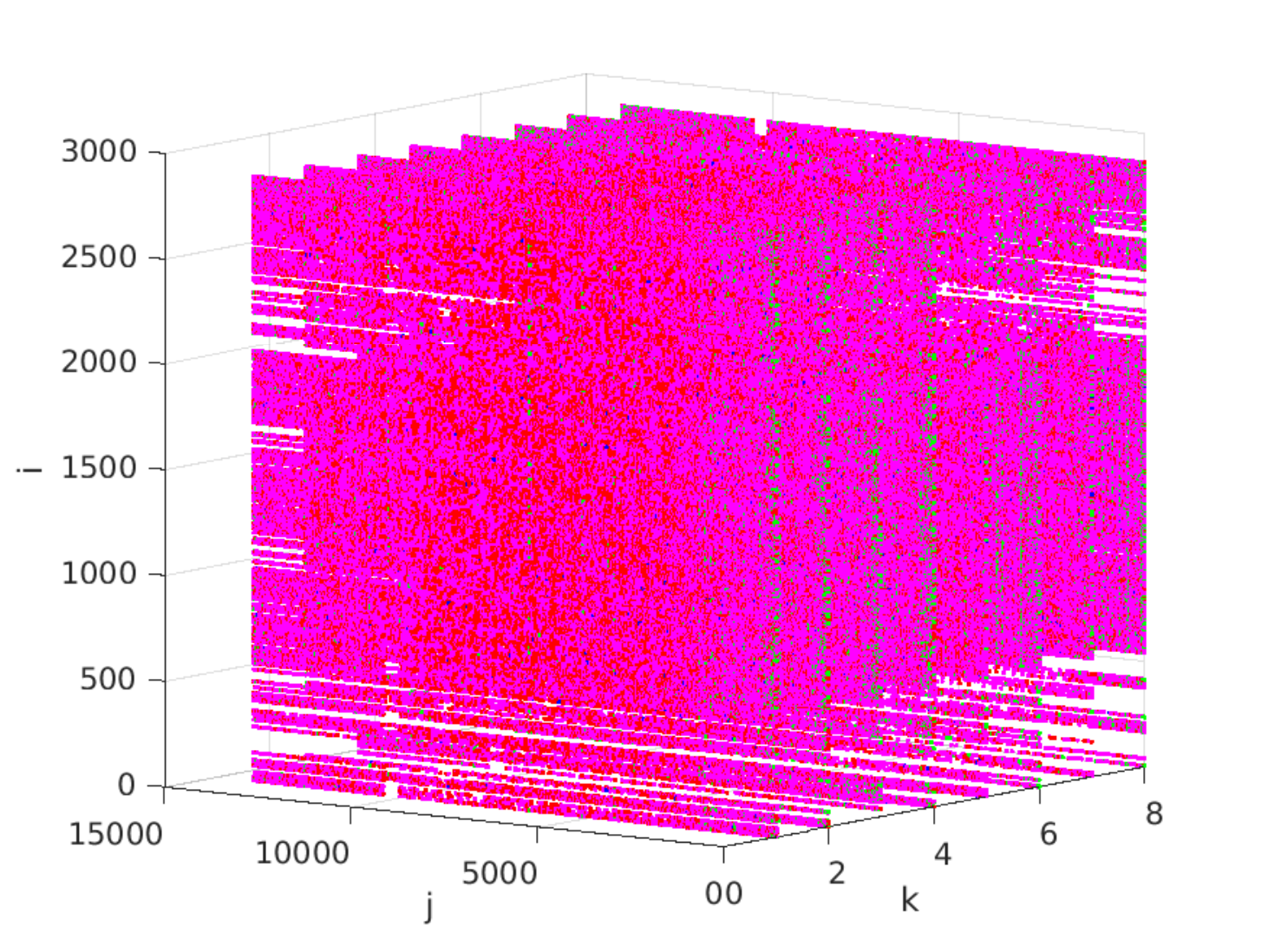}
\caption{Left: The tensor for the years 10-17,  reordered in the author and
  term modes. Right: The
  same tensor, but  the authors (vertical
  mode) have been ordered  as in Figure \ref{fig:A-NIPS-reord1-9}. }
\label{fig:A-NIPS-reord10-17}
\end{figure}

We made the same analysis for years 10-17. The results are given in
Figure \ref{fig:A-NIPS-reord10-17} and Table
\ref{tab:NIPStopics10-17}.
\begin{table}[htbp!] 
  \centering
  \caption{Two main topics at NeurIPS conferences during years
    10-17 (denoted B).    \label{tab:NIPStopics10-17}}   
  \begin{tabular}{llll|llllllll}
        &$T^1_B$& && &$T^2_B$\\
    \hline
 &model&        &models&         &usa&             &center&      \\
    &learning&     &figure&         &flies&           &school&      \\
    &algorithm&    &vector&         &independence&    &nec&         \\
    &training&     &error&          &israel&          &bialek&      \\
    &data&         &input&          &engineering&     &steveninck&  \\
    &network&      &case&           &huji&            &code&        \\
    &function&     &number&         &hebrew&          &william&     \\
    &state&        &em&             &princeton&       &rob&         \\
    &noise&        &likelihood&     &ruyter&          &tishby&      \\
    &set&          &functions&      &van&             &brenner&     \\
    &networks&     &probability&    &jerusalem&       &naftali&     \\
    &gaussian&     &class&          &neurobiology&    &universality&
  \end{tabular}
\end{table}
In the left panel of Figure \ref{fig:A-NIPS-reord10-17} we can see
that there are more authors in the years 10-17. In the right panel we
have used the same author ordering as in Figure
\ref{fig:A-NIPS-reord1-9}. The new authors, who did not contribute
during the years 1-9, are seen at the middle of the spy-plot. It is
also clear that many of the authors from years 1-9 are not present,
see the slightly  thinner pattern at the top and bottom.

Comparing Tables \ref{tab:NIPStopics1-17} and
\ref{tab:NIPStopics10-17} we see that it is the topics of the last 8
years that dominate for the whole period. The neural network topic is,
of course, the most significant one during the whole period, but note
that one can see slight differences between  Tables \ref{tab:NIPStopics1-9} and
\ref{tab:NIPStopics10-17}.

\section{Approximating  (1,2)-Symmetric Tensors
  by an   Expansion of  Rank-(2,2,1) Terms} 
\label{sec:expansion}

In some applications the purpose of the analysis is to find dominating
patterns in tensor data. In the network traffic example in Section
\ref{sec:network-logs} below there are ip-addresses, from which
many malicious attacks are sent to a large number of target 
addresses. The data are organized as a temporal sequence of graphs
that represent communication between ip-addresses.
In order to identify the attackers one needs to find dominating
subgraphs. Also, it may be interesting to study temporal communication
patterns.

We assume that the nonnegative tensor $\cA \in \RR^{m \times m \times n}$ is
$(1,2)-$symmetric and that the third mode is the temporal mode: each
3-slice $\cA(:,:,k)$ is the adjacency matrix of an undirected
communication graph. Our analysis is based on Proposition
\ref{prop:0BB0}, 
proved in \cite{eldehg20a}. Let $\cC \in \RR^{m_1 \times m_2 \times
  n}$, and  define $\cC^{'} \in \RR^{m_2 \times m_1 \times n}$ to be the
tensor such that each 3-slice is the transpose of the corresponding
slice in $\cC$.

\medskip

\begin{proposition}\label{prop:0BB0}
  Let  the  (1,2)$-$symmetric, nonnegative tensor $\cA \in \RR^{m
    \times m \times    n}$ have the structure
  \[
    \cA =
    \begin{pmatrix}
      0 & \cC \\
      \cC^{'} & 0
    \end{pmatrix}, \qquad \cC \in \RR^{m_1 \times m_2 \times n}, \quad
    m_1 + m_2 =m.
  \]
  Assume that the best rank-$(1,1,1)$ approximation of
  $\cC$ is unique, given by nonnegative
  $(u,v,w)$, and let $\tmr[1,2]{\cC}{u,v} = \tau
  w$, where $\tau>0$. Define the matrix $U$,
  \begin{equation}
    \label{eq:U-struct}
    U =
    \begin{pmatrix}
      0 & u\\
       v & 0
    \end{pmatrix}. 
  \end{equation}
  Then the best rank-$(2,2,1)$ approximation of $\cA$ is given by
  $(U,U,w)$, and
  \[
    \RR^{2 \times 2 \times 1} \ni \cF =
    \begin{pmatrix}
      0 & \tau \\
      \tau & 0
    \end{pmatrix}.
  \]  
\end{proposition}

 The rank-$(2,2,1)$ approximation $\tml{U,U,w}{\cF}$  can be written
\[
  \tml{U,U,w}{\cF} = \tml[3]{w}{\cB}, \qquad \cB = \tml[1,2]{U,U}{\cF}
  \in \RR^{m \times m \times 1},
\]
where $\cB$ can be identified with a matrix $B \in \RR^{m \times
  m}$. We can also write
\begin{equation}
  \label{eq:F-diag}
  B = \widehat{U} \widehat{F} \widehat{U}^T,
  \qquad \widehat{U} =
  \frac{1}{\sqrt{2}}\begin{pmatrix}
    u & -u\\
    v & v
  \end{pmatrix}, \qquad
  \widehat{F} =
  \begin{pmatrix}
    \tau & 0\\
    0    & -\tau
  \end{pmatrix}.
\end{equation}

It is easy to see that if we embed $\cC$ symmetrically in a larger
\emph{zero tensor} $\cZ \in \RR^{M \x M \x n}$, with $M > m$, then the
corresponding results hold. Let the solution of the best
rank-$(2,2,1)$ approximation for $\cZ$ be $(U_0,U_0,w_0)$, which
consists of the elements of $(U,U,w)$ embedded in larger zero matrices
and vectors.  Furthermore, if we embed $\cC$ symmetrically in a
larger (1,2)-symmetric but nonnegative  tensor
$\cA_0 \in \RR^{M \x M \x n}$, which has much smaller norm,
$\| \cA_0 \| \ll \| \cC \|,$ then since the best approximation depends
continuously on the elements of the tensor \cite{elsa11}, the solution
of the best rank-$(2,2,1)$ approximation problem for $\cA=\cA_0 + \cC$
will be
\[
  (U,U,w) = (E,E,e) + (U_0,U_0,w_0), \qquad
  \|(E,e)\| \approx K \| \cC \| / \| \cA \|,
\]
where $E$ and $e$ are perturbations, and $K$ depends on the
conditioning of the approximation problem for $\cA$ \cite{elsa11}. The
estimate is essentially the same as when we compute a best rank-$1$
approximation of a matrix (principal components, singular value
decomposition). This implies that we can use the best rank-$(2,2,1)$
approximation to determine the dominating information in a data
tensor.

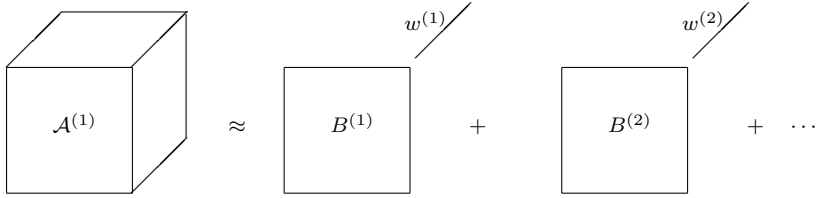
\begin{figure}[htbp!]    
{\setlength{\unitlength}{0.7pt}  
  \centering\footnotesize\bigskip
 \begin{picture}(200,100)(10,20)
\put (40,27){\line(1,0){68}} 
\put (40,27){\line(0,1){68}} 
\put (40,95){\line(1,0){68}} 
\put (108,27){\line(0,1){68}} 
\put (40,95){\line (1,1){30}}
\put (107.5,95){\line (1,1){30}}
\put (107.5,27){\line (1,1){30}}
\put (137,57){\line (0,1){67}}
\put (70,125){\line (1,0){68}}
%
\put (190,27){\line (0,1){68}}
\put (258,27){\line (0,1){68}}
\put (190,95){\line (1,0){68}}
\put (190,27){\line (1,0){68}}
\put (261,100){\line (1,1){30}}
%
\put (340,27){\line (0,1){68}}
\put (408,27){\line (0,1){68}}
\put (340,95){\line (1,0){68}}
\put (340,27){\line (1,0){68}}
\put (411,100){\line (1,1){30}}
\put (290,60){\mbox{$+$}}
\put (440,60){\mbox{$+ \quad \cdots$}}
\put (65,60){\mbox{$\cA^{(1)}$}}
\put (215,60){\mbox{$B^{(1)}$}}
\put (365,60){\mbox{$B^{(2)}$}}
\put (255,115){\mbox{$w^{(1)}$}}
\put (405,115){\mbox{$w^{(2)}$}}
\put (160,60){\mbox{$\approx$}}
\end{picture}
}
\caption{Symbolic illustration of the rank-(2,2,1) expansion of a tensor.}
\label{fig:expansion-symb}
\end{figure}

In the matrix SVD case, due to the fact that the singular value
decomposition diagonalizes the matrix, we  get the best
rank-$1$, rank-$2$, rank-$3$, etc.  approximations directly from the
SVD, see, e.g., 
\cite[Chapter 6]{eldenmatrix19}, \cite[Chapter 2]{govl13}. The
corresponding is not the case for the best rank-$(r_1,r_2,r_3)$
approximation of a tensor. Therefore, we suggest a procedure, where we
explicitly deflate the tensor\footnote{In the matrix case a similar
  procedure  would 
  give the SVD.}. Assuming that the best rank-(2,2,1)
approximation of the tensor $\cR^{(1)}:=\cA$ is given by
$(U^{(1)},U^{(1)},w^{(1)})$, with core tensor $\cF^{(1)}$ and
$B^{(1)}=\tml[1,2]{U^{(1)},U^{(1)}}{\cF^{(1)}}$, we compute
\[
  \cR^{(2)} := \cR^{(1)} - \tml[3]{w^{(1)}}{B^{(1)}}.
\]
Repeating this process, i.e. computing a rank-(2,2,1) approximation of
$\cR^{(\nu)},\, \nu=2,3,\ldots$, and putting $\cR^{(\nu+1)} = \cR^{(\nu)} -
\tml[3]{w^{(\nu)}}{B^{(\nu)}}$, then after $q-1$ steps   we have an expansion 
  \[
    \cA = \sum_{\nu=1}^q \tml[3]{w^{(\nu)}}{B^{(\nu)}} + \cR^{(q)}.
  \]
Clearly the matrices
$B^{(\nu)}$ have rank 2, so we can call the 
expansion a \emph{rank-(2,2,1) expansion.} The expansion is illustrated
symbolically in Figure \ref{fig:expansion-symb}. 
We present further details of the algorithm in connection with an
example in  Section \ref{sec:network-logs}. There the expansion is applied to a
nonnegative and sparse tensor, which has dominating patterns, so  for
this tensor the motivation above is applicable. However, in order to
avoid destroying sparsity  and nonnegativity modifications of the
procedure are introduced. 

The example described in Section \ref{sec:Reuters221}, does
not have that property of strongly dominating patterns,  i.e., the
eigenvalues of the matrices $F^{(\nu)}$ do not appear in pairs with
different signs. Still the expansion gives meaningful results. 

The algorithm is summarized in Section \ref{sec:sum-comp}. There we
also point out the features of the algorithm that ensure that for
sparse, nonnegative  tensors, the expansion is also sparse and
nonnegative.

\subsection{Rank-(2,2,1)  Expansion of Network Traffic Logs}
\label{sec:network-logs}
\begin{figure}[htbp!]    
\centering
\includegraphics[width=.45\textwidth]{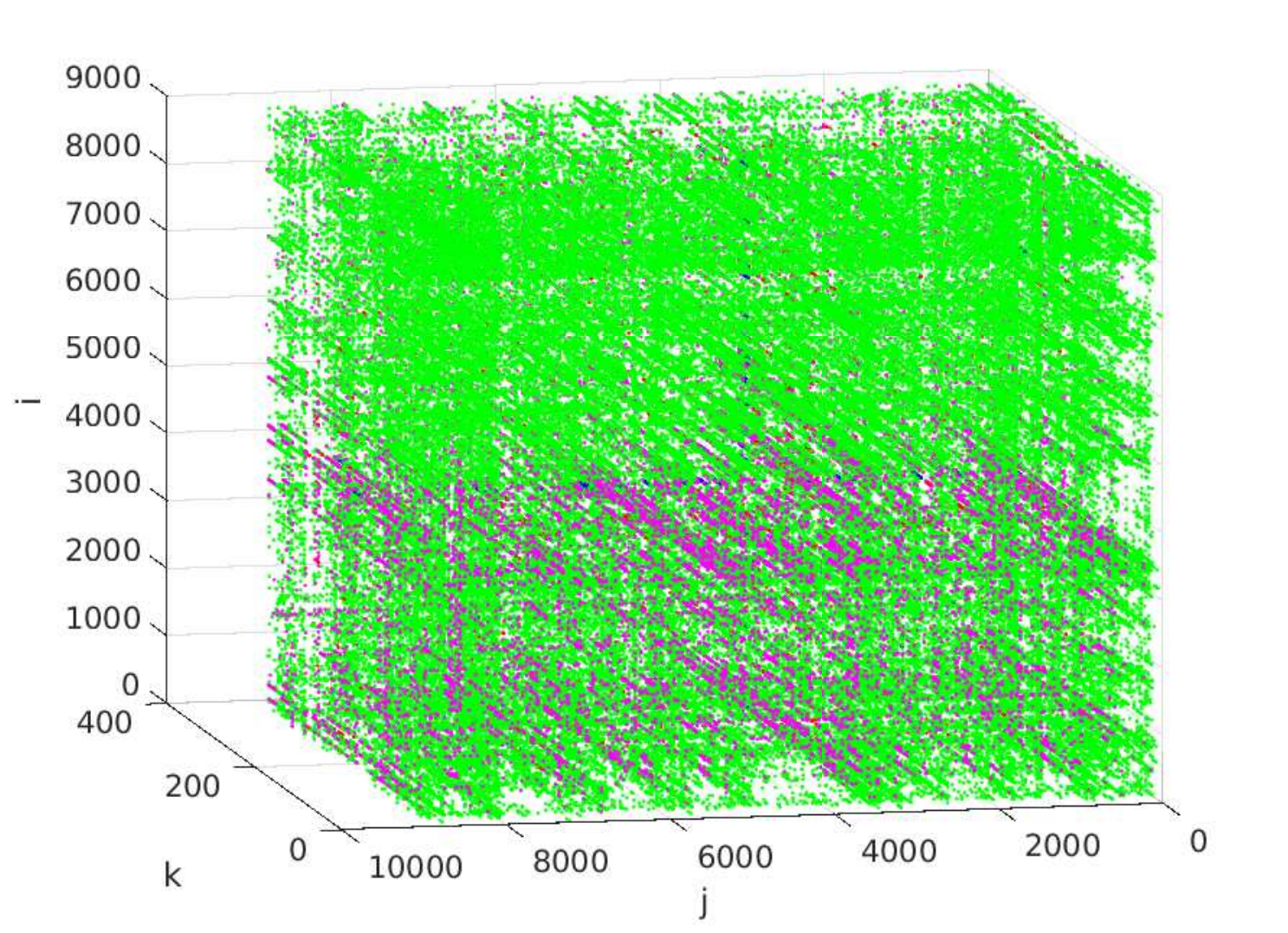}
\caption{1998DARPA tensor, original ordering.}
\label{fig:1998-tensor}
\end{figure}

\begin{figure}[htbp!]    
\centering
\includegraphics[width=.45\textwidth]{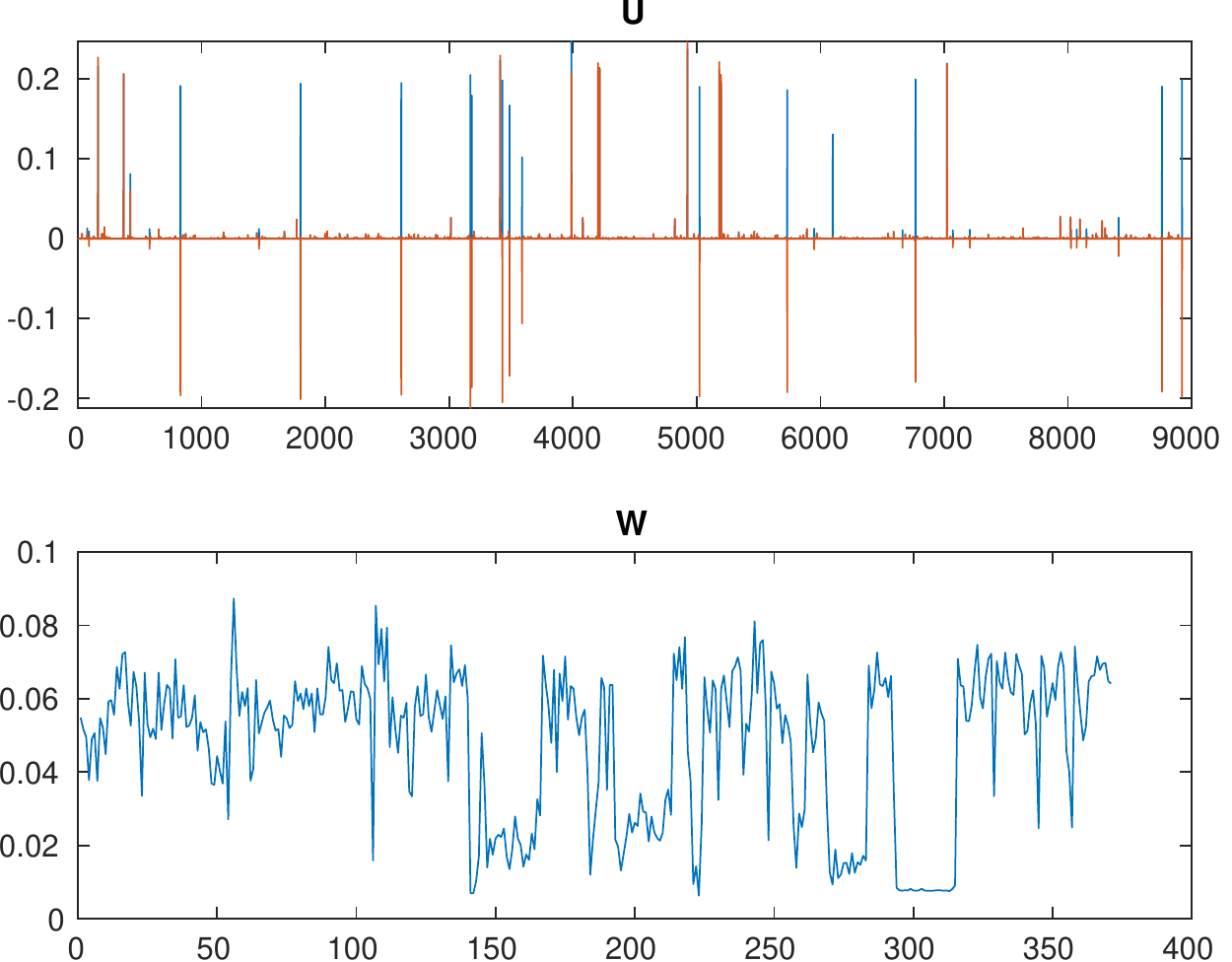}
\includegraphics[width=.45\textwidth]{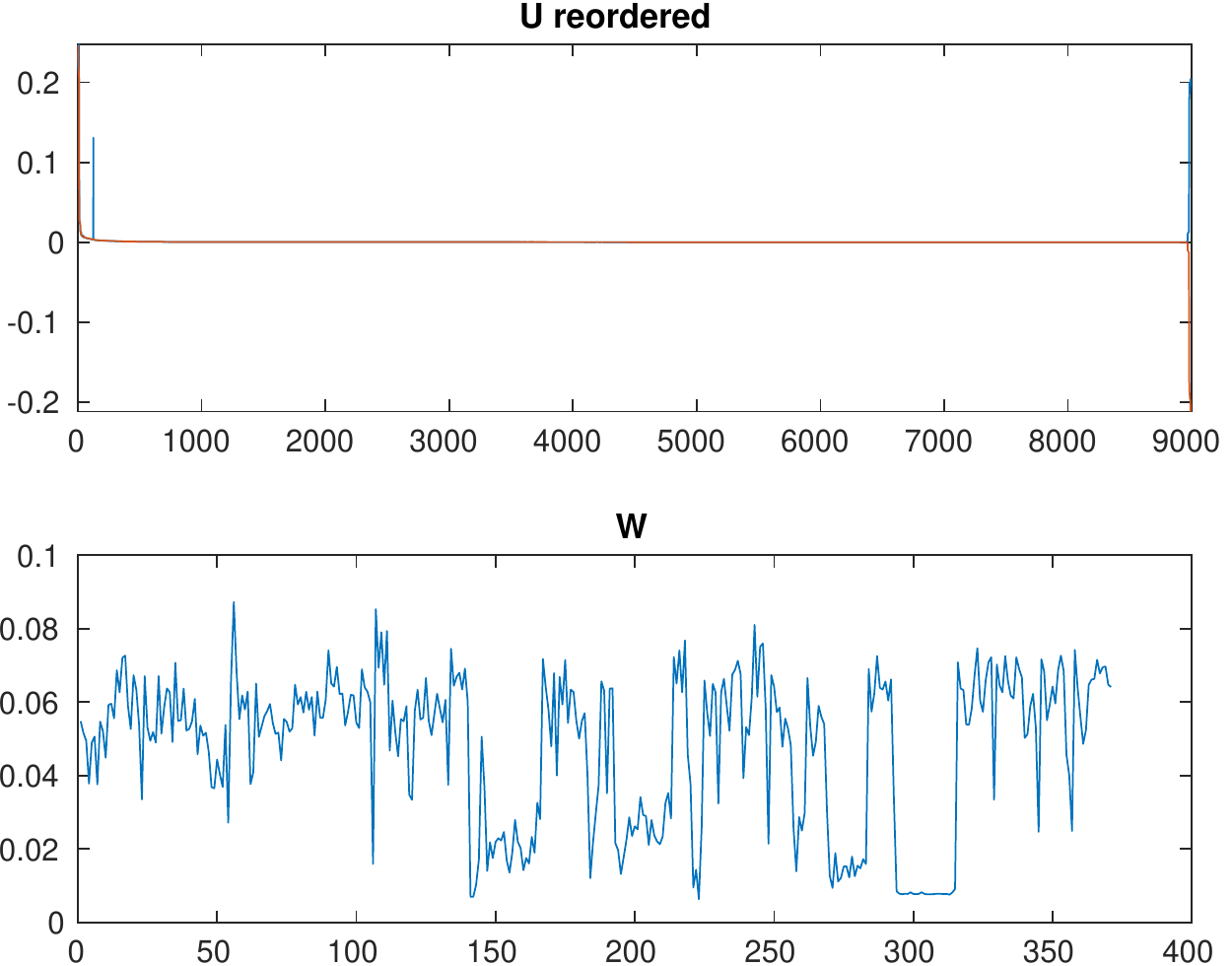}
\caption{Top panels: the vectors $u_1^{(1)}$ (blue),   $u_2^{(1)}$ (red). Bottom
  panels: $w^{(1)}$.  Left: Original
  ordering. Right: $(u_1^{(1)}\; u_2^{(1)})$ are  reordered  so that the elements
  of $u_2^{(1)}$ are monotonically decreasing.  } 
\label{fig:1998-UW}
\end{figure}

In \cite{jpfsk16} network traffic logs are analyzed in order to
identify malicious attackers. The data are called the 1998 DARPA
Intrusion Detection Evaluation Dataset and were first published by the
Lincoln Laboratory at
MIT\footnote{\url{http://www.ll.mit.edu/r-d/datasets/}\\
  \url{1998-darpa-intrusion-detection-evaluation-dataset.}}. 
We downloaded the data set from
\url{https://datalab.snu.ac.kr/haten2/} in October 2018. The records
consist of (source IP, destination IP, port number, timestamp). In the
data file there are about 22000 different IP addresses. We chose the subset
of 8991 addresses that both sent and received messages. The time span
for the data is from June 1 1998 to July 18, and the number of
observations is about 23 million. We merged the data in time by 
collecting every 63999 consecutive observations into one bin. Finally
we symmetrized the tensor $\cA  \in \RR^{m \times m \times n}$,
where $m=8891$  and $n=371$,  
so that
\[
  a_{ijk}=
  \begin{cases}
    1 & \mbox{ if } i  \mbox{ communicated with } j \mbox{ in time
      slot } k \\
    0 & \mbox{ otherwise.}
  \end{cases}
  \]
In this example we did not normalize the slices of the
tensor: The 3-slices are extremely  sparse, and are likely to have
non-zero patterns corresponding to several disconnected subgraphs. The
normalization \eqref{eq:normalize} would make a small, insignificant
subgraph have the same largest eigenvalue as a large subgraph, which
would give the small subgraph too much weight. Instead we scaled the
slices to have 
  Frobenius norm   equal to 1. This problem is quite well-conditioned
  and the Krylov-Schur algorithm converged in 5-17 iterations. We  refer
  to this as  the 1998DARPA tensor, and it is illustrated in 
Figure \ref{fig:1998-tensor}. 
For this type of data one relevant task is to identify the dominating
IP addresses, that communicate with many others. These  may be spammers. 
We will now demonstrate how we can extract such information by computing
the rank-$(2,2,1)$ approximation of the tensor.
The
$U^{(1)}=(u_1^{(1)}\;u_2^{(1)})$   and $w^{(1)}$ vectors are shown in Figure
\ref{fig:1998-UW}. 
The eigenvalues of the core tensor $\cF^{(1)}$  were 6.35 and -6.12, and 
 the rank-$(2,2,1)$ approximation
accounted for quite a 
large part of $\cA$: $\| \cA \| \approx 19.26$  and $\|
\cF^{(1)}\| \approx 8.82$. This indicates that there is a dominating
pattern in $\cA$  of the type  of Proposition \ref{prop:0BB0}. 

\begin{figure}[htbp!]    
  \centering
  \begin{minipage}{0.5\textwidth}
    \centering
     \includegraphics[width=0.8\textwidth]{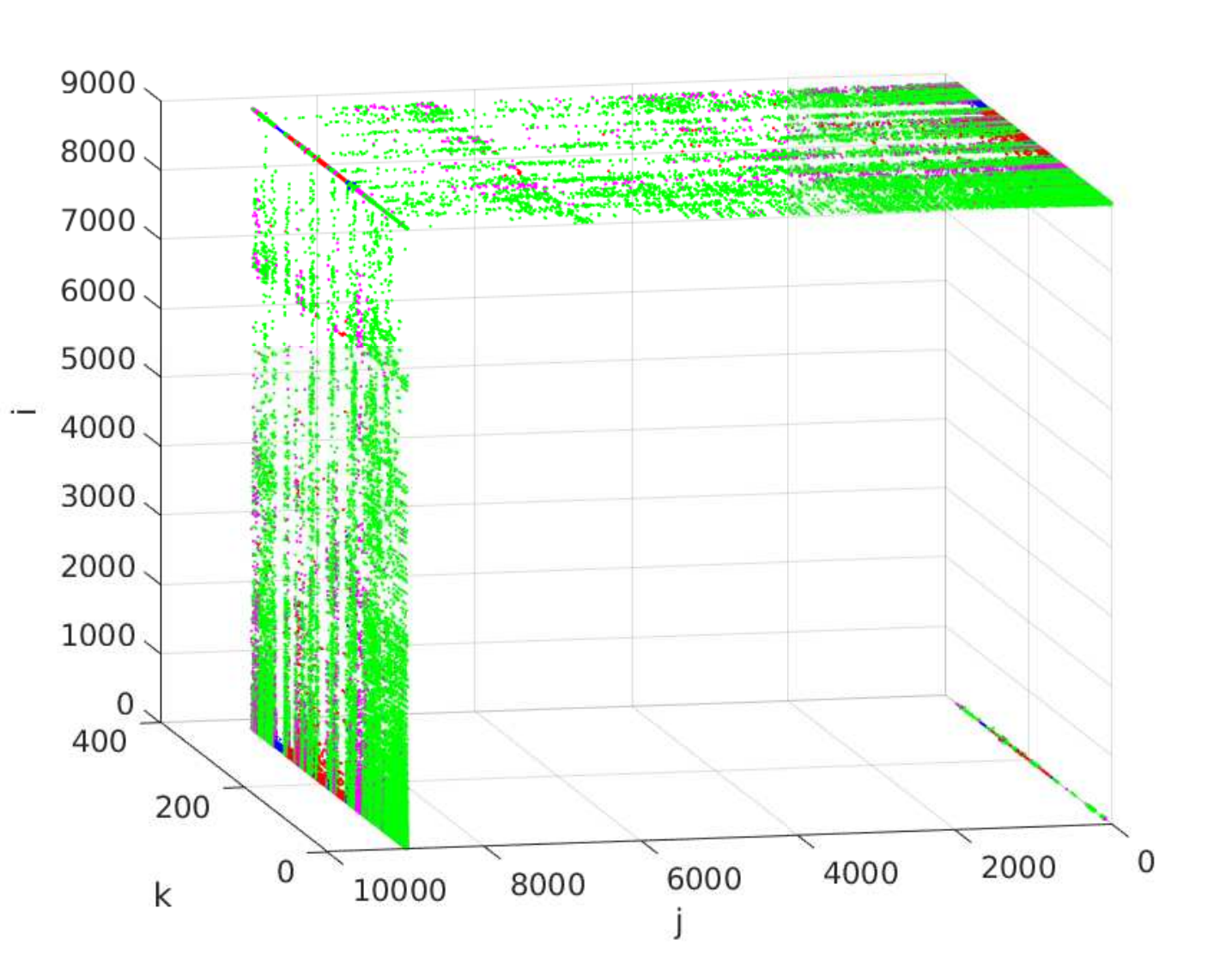}
  \end{minipage}%
  \begin{minipage}{0.5\textwidth}
    \centering
  $$\begin{bmatrix}
     1.5  &    7.0  &   11.6\\
    7.0   &    0.11  &       0.03\\
   11.6   &       0.03  &  2.8
 \end{bmatrix}$$
  \end{minipage}
\caption{1998DARPA tensor,  reordered in the
  (1,2) modes.  Right: Norms of  tensor blocks $\cA(I,J,:)$,
  where $I$ and $J$ are sequences of indices, such that   the outermost blocks
  are   200 elements wide.
}
\label{fig:1998-tensor-reordered}
\end{figure}

The $U^{(1)}$ vectors are very sparse and spiky, which shows that
very few addresses are dominating the communication.  When we reordered
so that $u_2^{(1)}$ 
became monotonic, the most active addresses (the spikes)
are placed at 
the beginning and the end. Applying the same reordering to the tensor
in the (1,2)-modes, we got the left panel in  Figure
\ref{fig:1998-tensor-reordered}. 
In the sequel, for notational convenience, we let $\cA$  denote the reordered
tensor and  $U^{(1)}$ the correspondingly reordered matrix. 
Clearly there is a relatively small group of addresses (upper left
corner in Figure \ref{fig:1998-tensor-reordered}) that
communicate with the majority of  the others, and another very small group
(lower right corner) that communicate only within the  small group
and with the first group (upper
right corner). The majority of addresses (middle) do not
communicate with each other. 

\begin{figure}[htbp!]    
  \centering
  \begin{minipage}{0.6\textwidth}
    \centering
     \includegraphics[width=\textwidth]{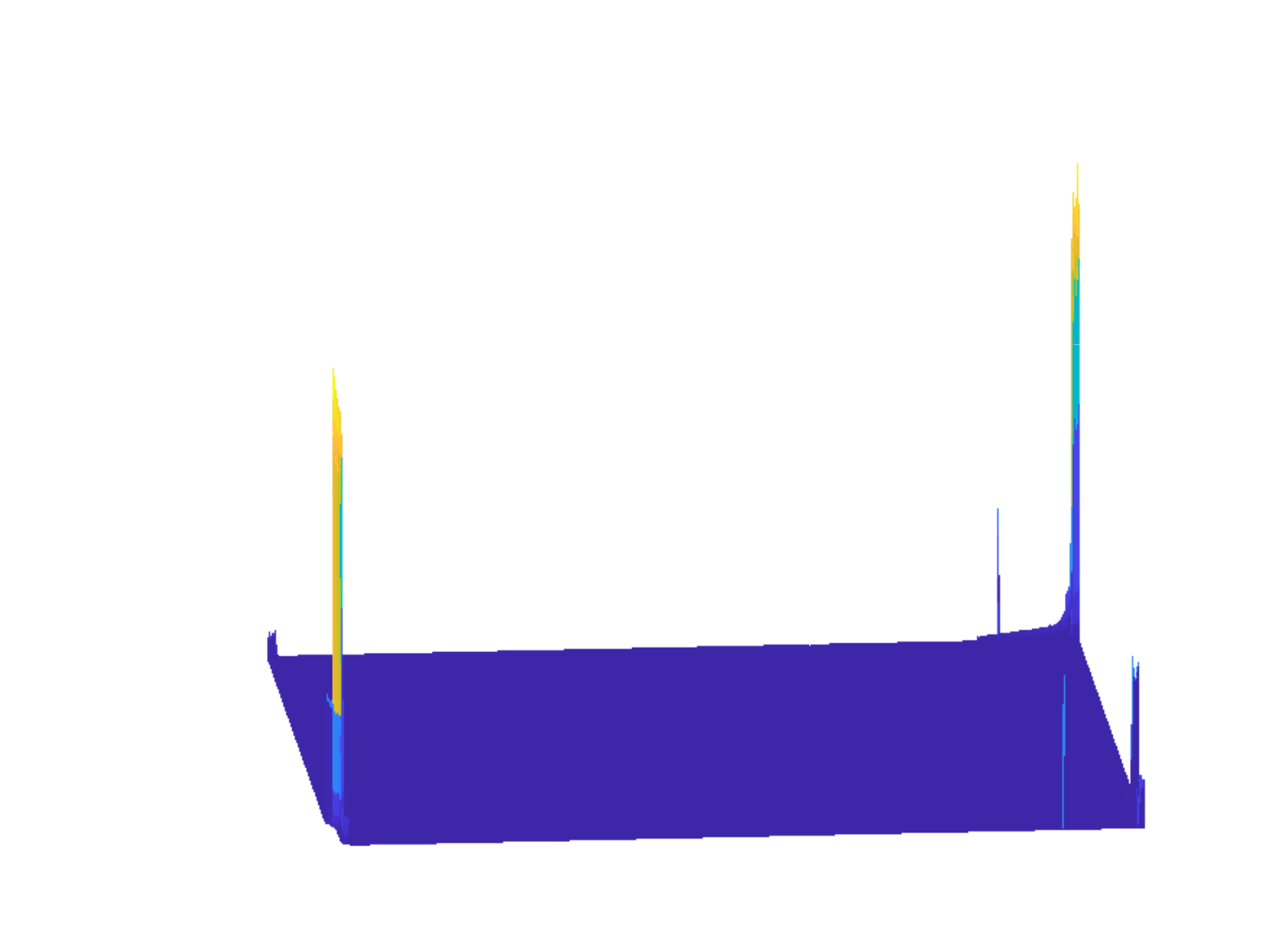}
  \end{minipage}%
  \begin{minipage}{0.4\textwidth}
    \centering
  $$\begin{bmatrix}
     0.14  &    0  &   6.2\\
    0   &    0  &    0\\
   6.2   &       0  &  0.85 
 \end{bmatrix}$$
  \end{minipage}
  \caption{  Left: Mesh plot of the matrix $\widehat{B}^{(1)}$. For
    visibility the scale is much higher at the corners. 
    Right: Norms of $200 \times 200$ blocks at the four corners; the
    0's are the norms of the much larger other blocks. }
\label{fig:1998-B}
\end{figure}

\begin{figure}[htbp!]    
  \centering
     \includegraphics[width=0.5\textwidth]{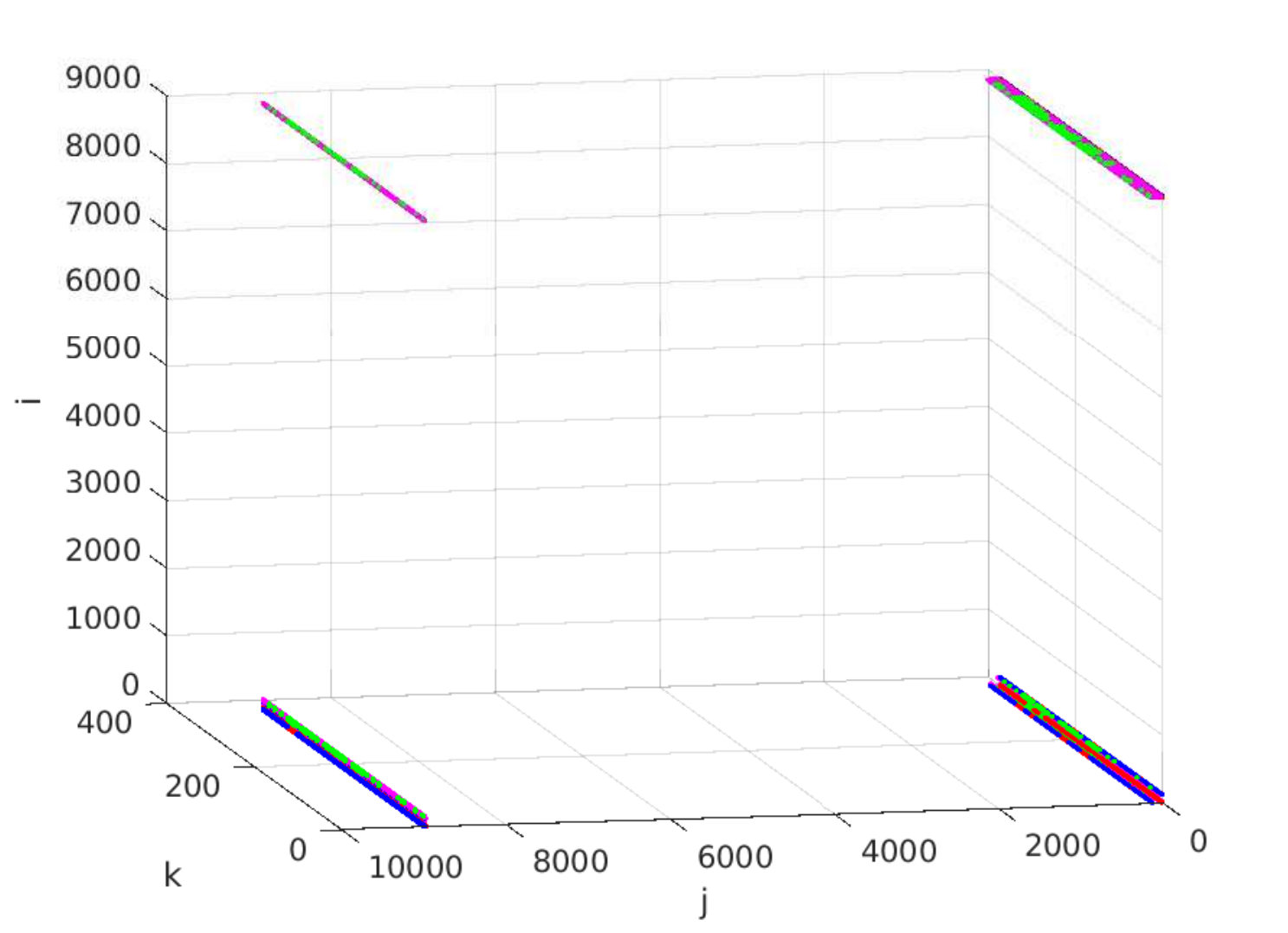}
\caption{The approximation tensor $\tml[3]{w^{(1)}}{\widehat{B}^{(1)}}$.  }
\label{fig:1998-Aapprox}
\end{figure}

As the columns of $U^{(1)}$ are orthogonal, we can not expect the
matrix $B^{(1)} = \tml[1,2]{U^{(1)},U^{(1)}}{\cF^{(1)}}$ to be sparse;
in this example it even  has small negative
elements (the smallest element of $B^{(1)}$ is $-0.018$
and the largest is $0.63$). In order to have a sparse and nonnegative
approximation we define a new matrix $\widehat{B}^{(1)}$:  we let
$b_{\max}$ be the 
largest element in $B^{(1)}$, and for $\theta=0.01$ we  put
\begin{equation}
  \label{eq:b-theta}
  \widehat{b}_{ij}^{(1)}=
  \begin{cases}
    {b}_{ij}^{(1)} & \mbox{if } {b}_{ij}^{(1)}>\theta
    \,b_{\max},\\
    0, & \mbox{otherwise.}
  \end{cases}
  \end{equation}
  Thus we approximated
  \[
    \cR^{(1)}:= \cA = \tml[3]{w^{(1)}}{\widehat{B}^{(1)}} + \cR^{(2)},
  \]
  where   $\cR^{(2)}$ is a residual tensor. The approximation is illustrated
  in Figures \ref{fig:1998-B} and \ref{fig:1998-Aapprox}. Note that
  the low rank approximation is 
  concentrated in this blocks where the reordered tensor  has highest
  ``density'', cf. Figure \ref{fig:1998-tensor-reordered}.

The addresses that dominate the communication can be identified in the
matrix $\widehat B^{(1)}$. In fact, $\widehat B^{(1)}$ can be interpreted as the
adjacency matrix of a communication graph for the dominating
communication, and the vector $w_1^{(1)}$ (Figure \ref{fig:1998-UW}) shows
 the strength of the communication at the different time slots.

Next we deflated the tensor,
\[
  \cR^{(2)} =  \cR^{(1)} - \tml[3]{w^{(1)}}{\widehat{B}^{(1)}}. 
\]
The norms of the  blocks of  the deflated tensor $\cR^{(2)}$
(partitioned as in Figure  \ref{fig:1998-tensor-reordered}) were
\[
  \begin{bmatrix}
    1.5 & 7.0 & 9.6\\
    7.0 & 0.11& 0.036\\
    9.6 & 0.036& 2.8
  \end{bmatrix}
\]
  
\begin{figure}[htbp!]    
\centering
\includegraphics[width=.6\textwidth]{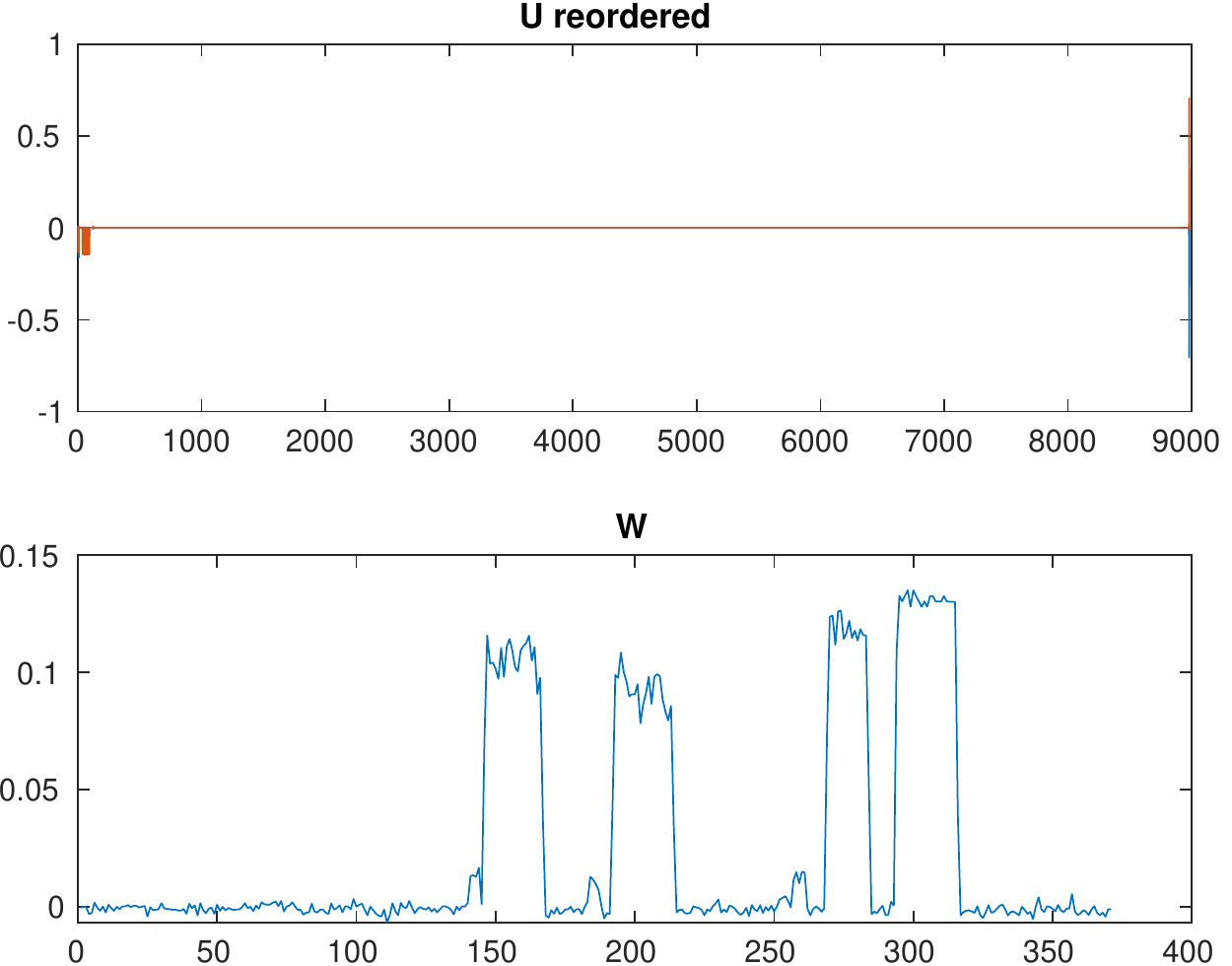}
\caption{Analysis of the deflated tensor $\cR^{(2)}$ . Top: the vectors
  $u_1^{(2)}$ (blue),   $u_2^{(2)}$ (red). Bottom: $w^{(2)}$. } 
\label{fig:1998-UW-1}
\end{figure}

We performed the same analysis on $\cR^{(2)}$ with the analogous 
thresholding for $\widehat{B}^{(2)}$. 
%
%
 The eigenvalues of the  core tensor were $5.1$ and $-5.1$ 
The rank-$(2,2,1)$ approximation accounted for quite a
large part of $\cR^{(2)}$: $\| \cR^{(2)} \| \approx 17.13$  and $\| \cF^{(2)}\|
\approx 7.22$.  The $U^{(2)}=(u_1^{(2)}\;u_2^{(2)})$   and $w^{(2)}$ vectors are shown in Figure
\ref{fig:1998-UW-1}.
Comparing $w^{(2)}$ with the plot of $w^{(1)}$ in Figure
\ref{fig:1998-UW}, we see that the second most dominating
communication shown here took place essentially when the level of
communication in first  group was low (the angle between 
 $w^{(1)}$ and 
$w^{(2)}$ is approximately $86^\circ$).
\begin{figure}[htbp!]    
  \centering
  \begin{minipage}{0.6\textwidth}
    \centering
     \includegraphics[width=\textwidth]{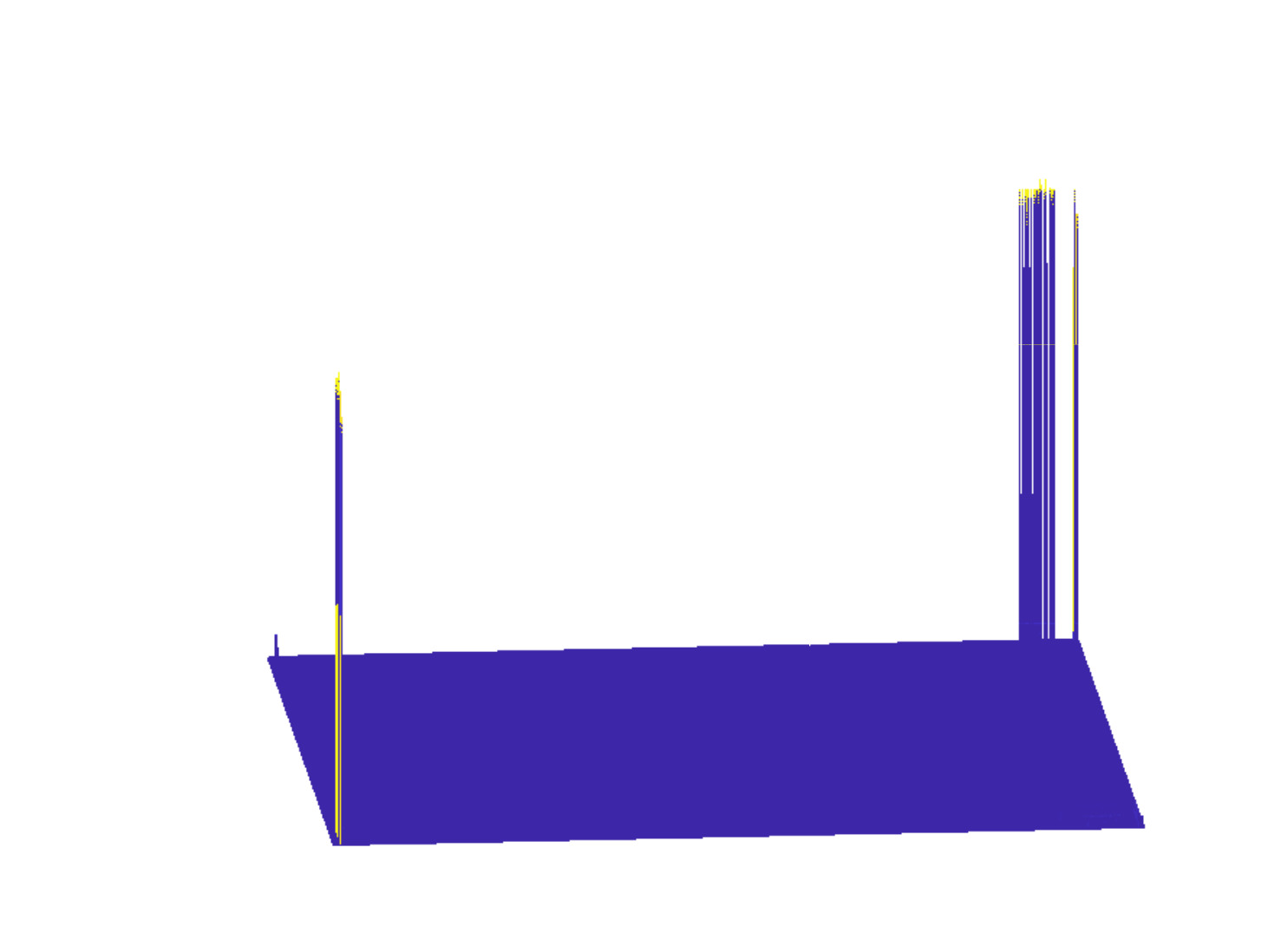}
  \end{minipage}%
  \begin{minipage}{0.4\textwidth}
    \centering
  $$\begin{bmatrix}
     0.077  &    0  &   5.1\\
    0   &    0  &     0\\
   5.1   &       0  &  0.19
 \end{bmatrix}$$
  \end{minipage}
\caption{ Left: Mesh plot of the matrix $\widehat{B}^{(2)}$. For
    visibility the scale is much higher at the corners. 
    Right: Norms of $200 \times 200$ blocks at the four corners; the
    0's are the norms of the much larger other blocks. }
\label{fig:1998-B2}
\end{figure}  

The matrix $\widehat{B}^{(2)}$ is shown in Figure
\ref{fig:1998-B2}. The approximation  tensor
$\tml[3]{w^{(2)}}{\widehat{B}^{(2)}}$  looked very similar to 
   $\tml[3]{w^{(1)}}{\widehat{B}^{(1)}}$ in Figure \ref{fig:1998-Aapprox}.

\begin{figure}[htbp!]    
\centering
\includegraphics[width=.6\textwidth]{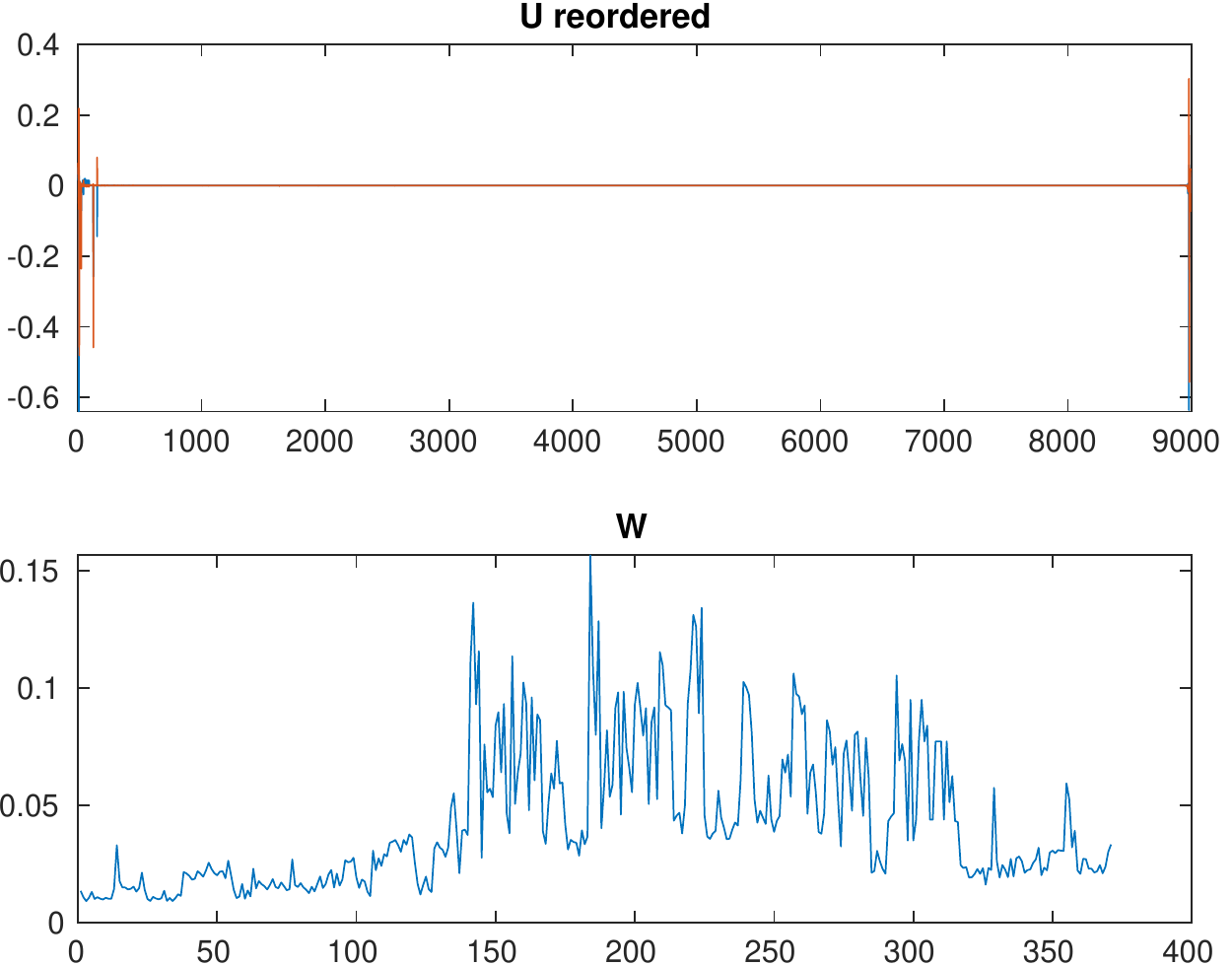}
\caption{Analysis of the deflated tensor $\cR^{(3)}$ . Top: the vectors
  $u_1^{(3)}$ (blue),   $u_2^{(3)}$ (red). Bottom: $w^{(3)}$. } 
\label{fig:1998-UW-2}
\end{figure}

  Again we deflated the tensor,
and performed another step of analysis on $\cB^{(3)}$.  
The eigenvalues of the  core tensor were $2.06$ and  $1.82$; 
this indicates that there was no longer any significant structure as
that in Proposition \ref{prop:0BB0}. 
$\| \cR^{(3)} \| \approx 15.54$ and $\| \cF^{(3)}\| \approx 2.76$. The
$U$ and $W$ vectors are shown in Figure \ref{fig:1998-UW-2}.  
\begin{figure}[htbp!]    
  \centering
  \begin{minipage}{0.6\textwidth}
    \centering
     \includegraphics[width=\textwidth]{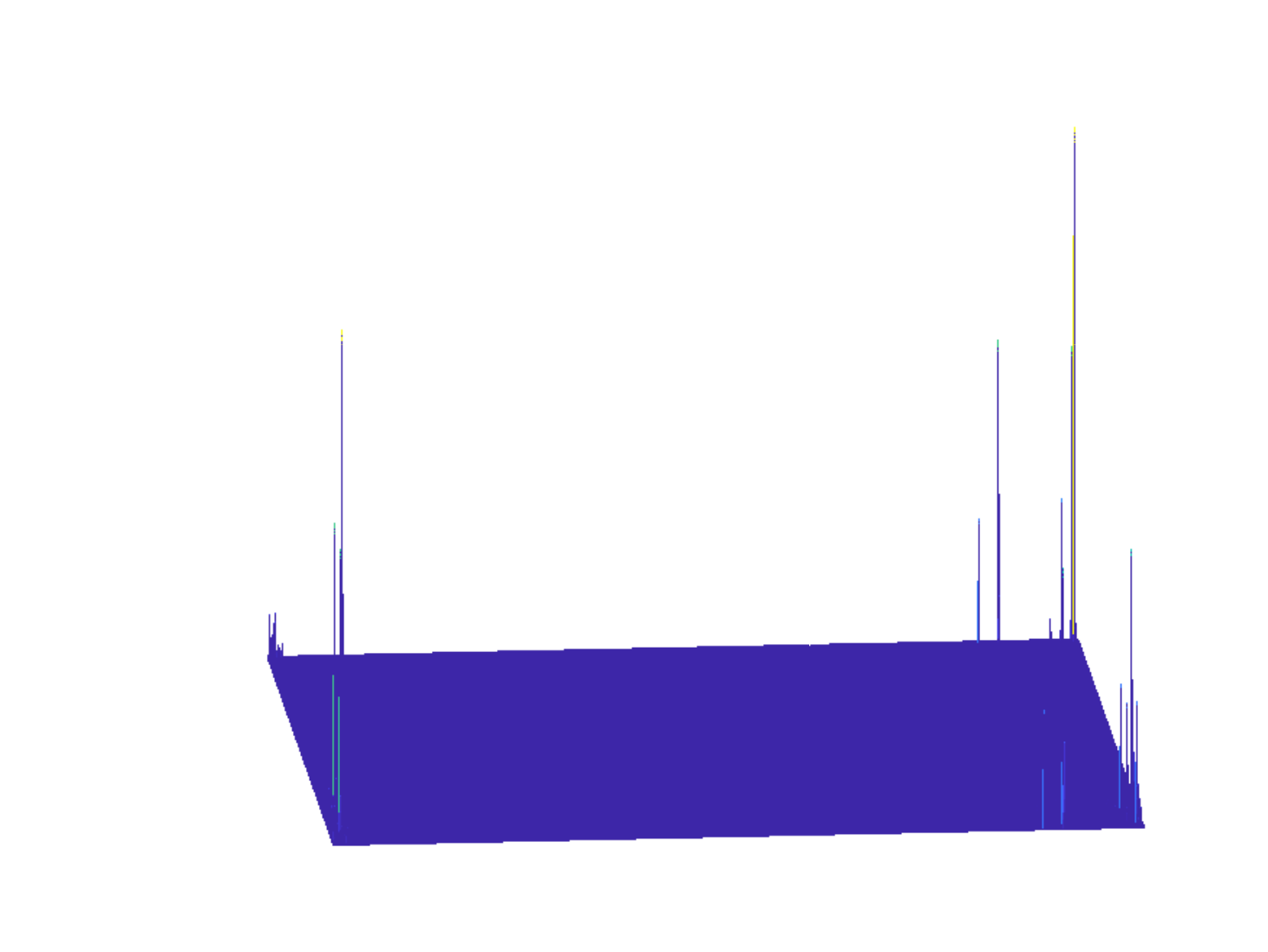}
  \end{minipage}%
  \begin{minipage}{0.4\textwidth}
    \centering
  $$\begin{bmatrix}
     0.18  &    0  &   1.3\\
    0   &    0  &       0\\
   1.3   &       0  &  0.95
 \end{bmatrix}$$
  \end{minipage}
\caption{ Left: Mesh plot of the matrix $\widehat{B}^{(3)}$. For
    visibility the scale is much higher at the corners. 
    Right: Norms of $200 \times 200$ blocks at the four corners; the
    0's are the norms of the much larger other blocks.}
\label{fig:1998-B3}
\end{figure}  
The matrix $\widehat{B}^{(3)}$ is illustrated in Figure \ref{fig:1998-B3}. The
approximation tensor $\tml[3]{w^{(3)}}{B^{(3)}}$ looked very similar
to those  in the previous steps.
 
  The procedure so far can be written in the form
  \[
    \cA^{(1)} \approx \sum_{\nu=1}^3 \tml[3]{w^{(\nu)}}{\widehat{B}^{(\nu)}},
  \]
  where the terms are nonnegative tensors of rank-$(2,2,1)$. Note that
neither  the vectors $w^{(\nu)}$  nor the
matrices $\widehat{B}^{(\nu)}$  are  constructed to be orthogonal
(with  Euclidean  inner 
products). However, for this example where the tensor $\cA$ is sparse and
nonnegative, we have computed a few terms of a
\emph{sparse, nonnegative low rank expansion}.

In the case when the 3-slices of the tensor $\cA$ are adjacency matrices of
graphs, we can interpret the matrices $\widehat{B}^{(\nu)}$ in  the expansion as
adjacency matrices of dominating or \emph{salient subgraphs}. The
element  $w_k^{(\nu)}>0$ of the vector $w^{(\nu)}$ can be seen as a
measure of how much of the salient subgraph  corresponding to
$\widehat{B}^{(\nu)}$ is present in  slice $k$ of the tensor.

To check how much the different terms in the expansion
overlapped, we computed the cosine of the angles between the matrices,
\[
  \cos(\widehat B^{(\nu)},\widehat B^{(\lambda)}) = \frac{\<\widehat  B^{(\nu)},\widehat  B^{(\lambda)}
    \>}%
  {\| \widehat  B^{(\nu)} \| \, \|\widehat  B^{(\lambda)} \|  },
\]
and give them below Table \ref{tab:BFnorms-netw} (the cosine is a measure of overlap between edges
in the communication graphs). 
\begin{table}[htbp!]
  \centering
  \caption{Network Traffic example. Top: Norms of
    $\widehat{B}^{(\nu)}$ and 
    $\cF^{(\nu)}$,  the largest and 
    smallest values of ${B}^{(\nu)}$ (i.e. before thresholding), and the
    eigenvalues of $\cF^{(\nu)}$ . Bottom:
    Cosines of the angles between the matrices $\| \widehat  B^{(\nu)} \|$,
    $C_B({\mu,\lambda})= \cos(\widehat B^{(\nu)},\widehat B^{(\lambda)})$.   
    \label{tab:BFnorms-netw}}
\begin{tabular}{l|ccccccc}
  \hline
  \noalign{\vskip 2pt}
  $\| \widehat{B}^{(\nu)} \|$&8.8&7.2&2.1\\
  $\| \cF^{(\nu)} \| $ &8.8&7.2&2.7\\
    max& 0.63&1.0&0.97\\
  min&$-0.02$&$-0.07$&$-0.08$\\
  $\lambda_1$&6.3&5.1& 2.1\\
  $\lambda_2$&$-6.1$&$-5.1$&1.8   \\
 \hline
\end{tabular}
\medskip
\[
     C_B=
  \begin{bmatrix}
    1&0.07&0.06\\
    &1&0.006\\
    & &1
  \end{bmatrix}
\qquad .
\]
\end{table}
There we also give the norms on $\widehat{\cB}^{(\nu)}$ and
$\cF^{(\nu)}$. As the identity $\|{B}^{(\nu)}\| = \| \cF^{(\nu)} \|$
holds, it is seen that  the elements in ${\cB}^{(1)}$ and
${\cB}^{(2)}$ that are removed in   the thresholding
\eqref{eq:b-theta} are insignificant.  

Summarizing the experiment, we have computed an approximate  expansion
of the tensor that corresponds to three relatively small and
almost disjoint groups of IP addresses that dominate the
communication. The communication corresponding to  the first two
terms in the expansion is taking place at different times.  If we are
 looking only for structure of the type of Proposition \ref{prop:0BB0},
then we should not use the third term.

\subsection{Low Rank Expansion of News Text}
\label{sec:Reuters221}

In this section we will  compute a low rank expansion of seven terms,
  \[
    \cA \approx \sum_{\nu=1}^7 \tml[3]{w^{(\nu)}}{\widehat{B}^{(\nu)}},
  \]
for  the Reuters
news text tensor that we analyzed in Section \ref{sec:topic}. 
 The matrices $\widehat{B}^{(\nu)}$ are considered as adjacency matrices of
subgraphs of the large graph corresponding to cooccurrence of all 13332
terms in the sentences of the texts, cf. Section \ref{sec:topic}. 

\begin{figure}[htbp!]    
 \centering
  \begin{minipage}{0.7\textwidth}
    \centering
  \includegraphics[width=1\textwidth]{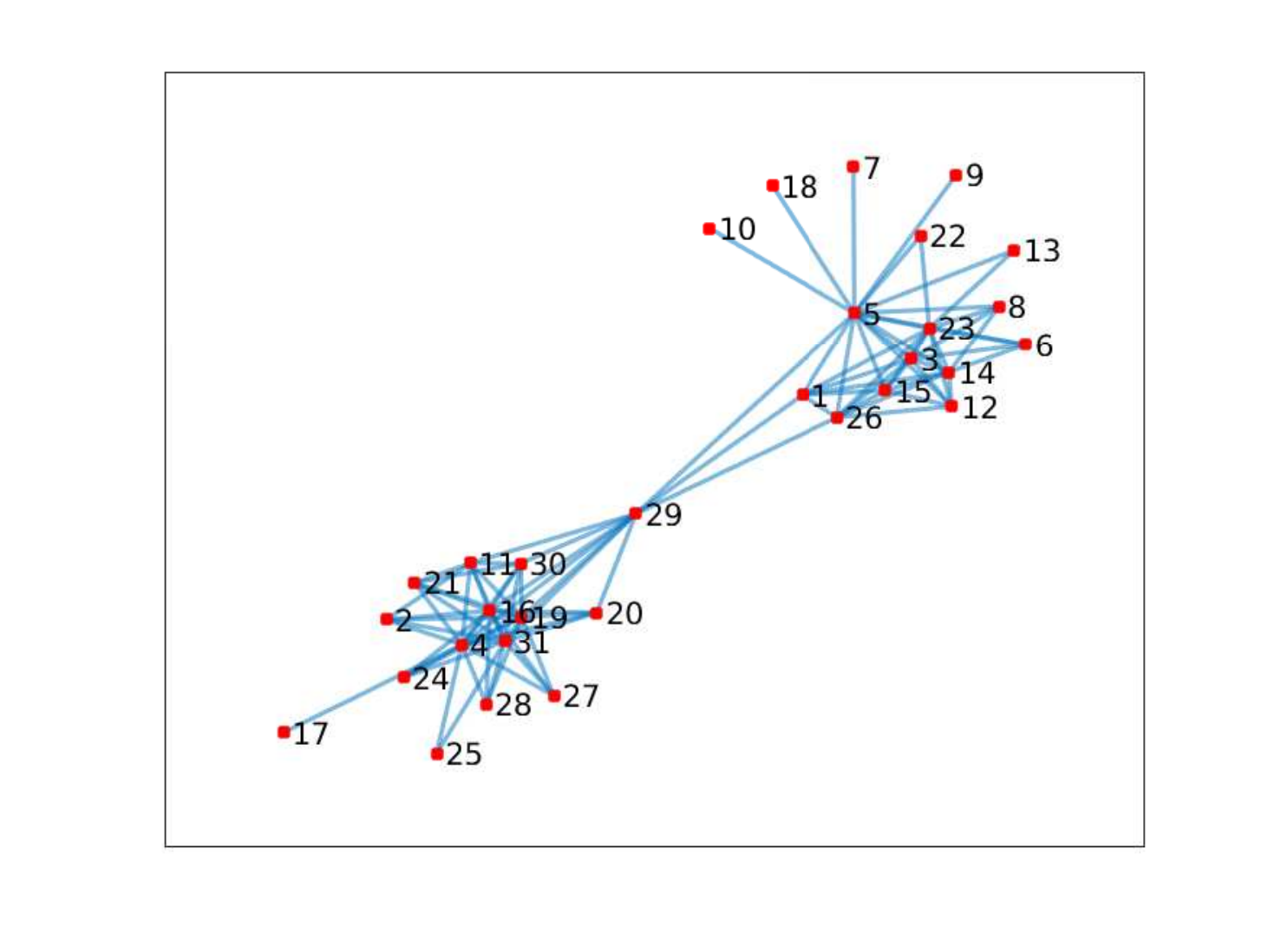}
  \end{minipage}%
  \begin{minipage}{0.3\textwidth}
    \centering
  \includegraphics[width=0.7\textwidth]{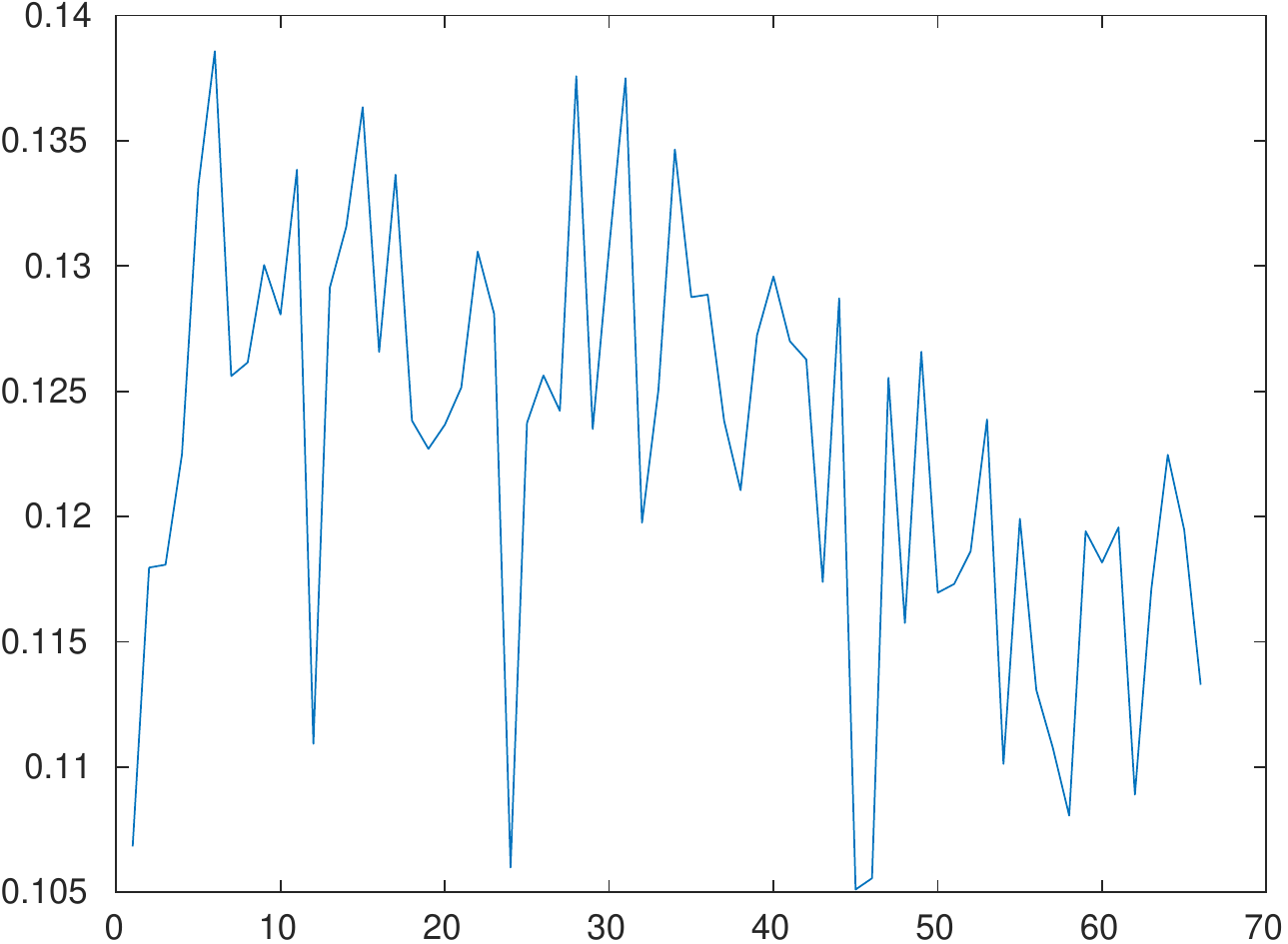}
\end{minipage}
       \begin{minipage}{1\textwidth}
    \footnotesize
    \begin{tabular}{llllll}
1. afghanistan  &  8. fugitive & 15. network & 22. rule & 29. united\_states\\
2. airliner     &  9. group    & 16. new\_york& 23. saudi-born& 30. washington\\
3. al\_quaeda  & 10. guerrilla& 17. official& 24. sept & 31. world\_trade\_ctr\\
4. attack       & 11. hijack   & 18. organization& 25. suicide \\
5. bin\_laden   & 12. islamic  & 19. pentagon& 26. taliban\\
6. dissident    & 13. leader   & 20. people  & 27. tower\\
7. exile        & 14. militant & 21. plane   & 28. twin
    \end{tabular}
       \end{minipage}
\caption{First term in the expansion of the news text
  tensor. We illustrate the graph, the vector $w^{(1)}$, and the keywords
  corresponding to the vertices in the graph. Two subgraphs are
  visible, connected via the term ``united\_states''; the  two subtopics
  are close to $T^1$ and $T^2$  in Table \ref{tab:terms66}. }
\label{fig:Reuters-1}
\end{figure}

\begin{figure}[htbp!]    
 \centering
  \begin{minipage}{0.7\textwidth}
    \centering
  \includegraphics[width=1\textwidth]{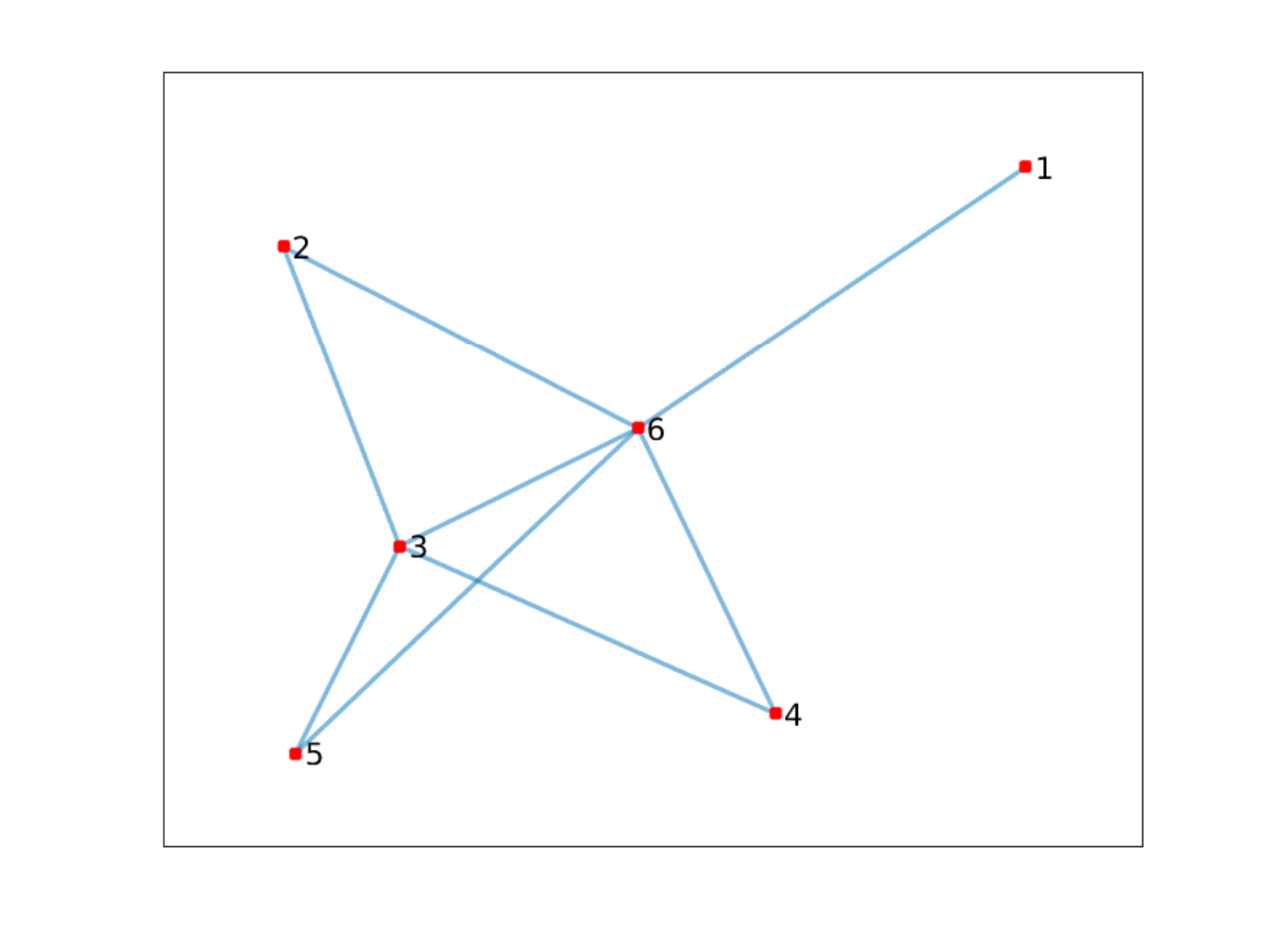}
  \end{minipage}%
  \begin{minipage}{0.3\textwidth}
    \centering
  \includegraphics[width=0.7\textwidth]{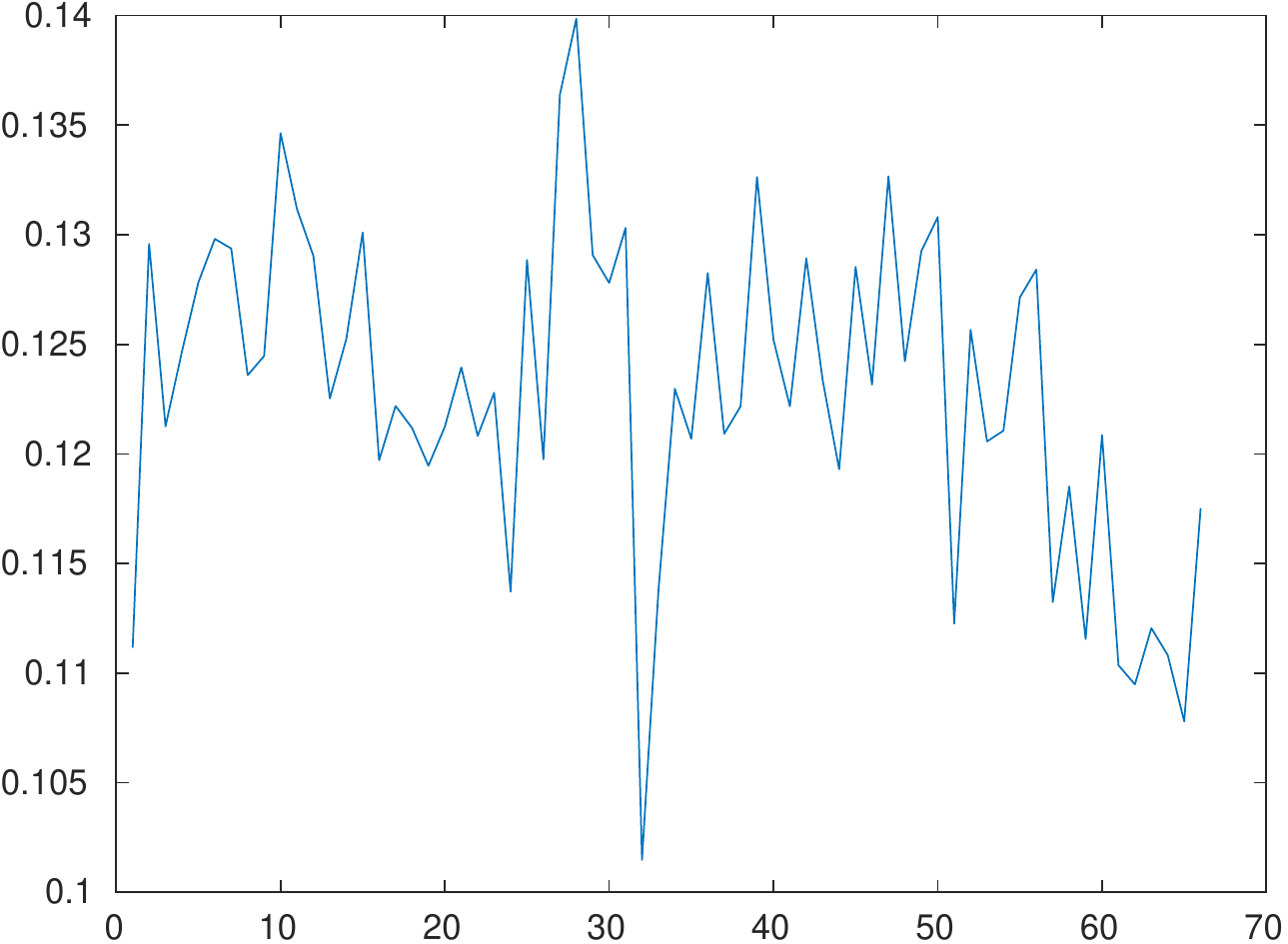}
\end{minipage}
       \begin{minipage}{1\textwidth}
    \footnotesize
    \begin{tabular}{llllll}
    1.    news & 2. reporter & 3. reuter & 4. rumsfeld &
    5. sec & 6.    tell      
    \end{tabular}
       \end{minipage}
\caption{Second term.}
\label{fig:Reuters-2}
\end{figure}

\begin{figure}[htbp!]    
 \centering
  \begin{minipage}{0.7\textwidth}
    \centering
  \includegraphics[width=1\textwidth]{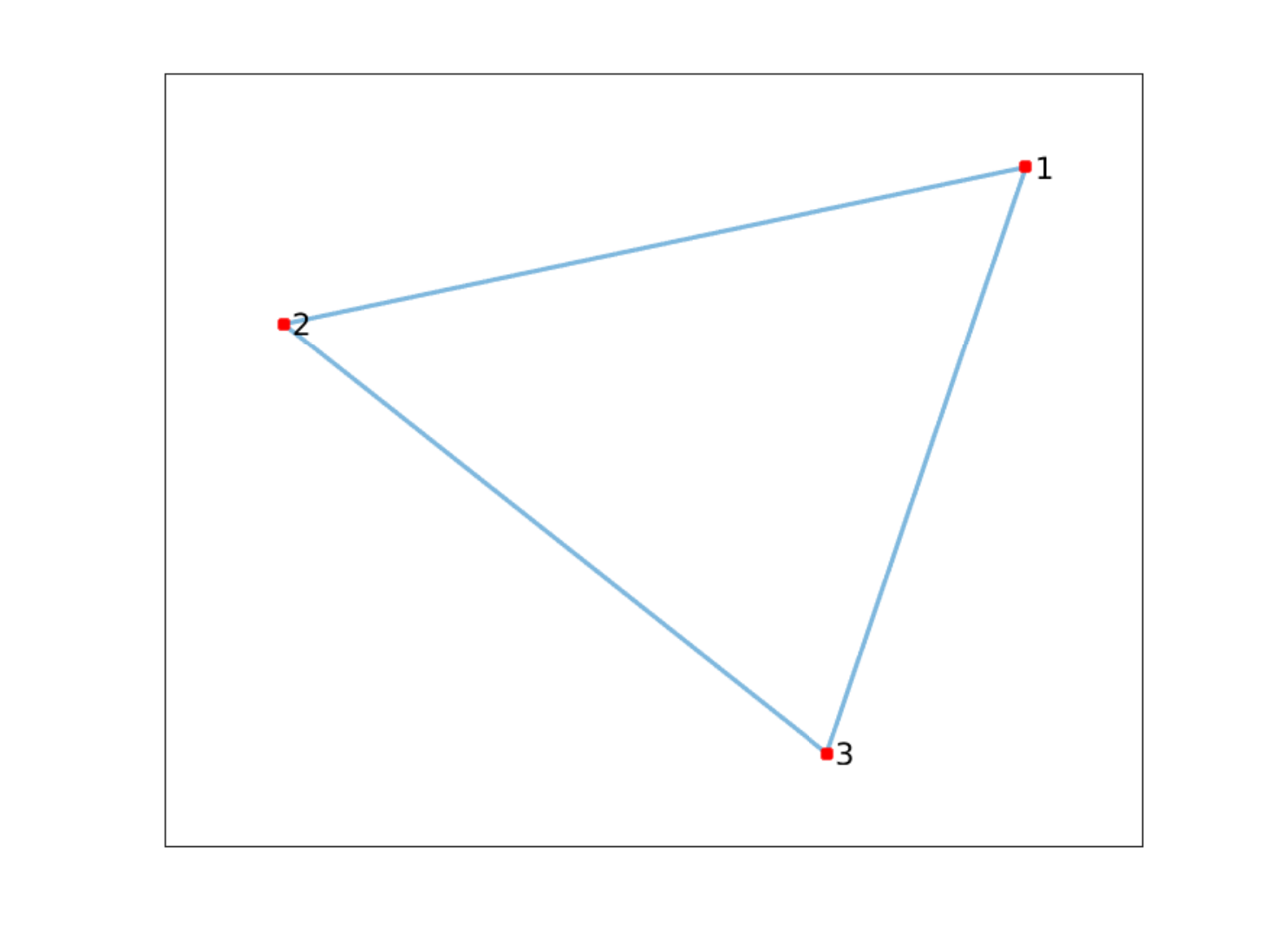}
  \end{minipage}%
  \begin{minipage}{0.3\textwidth}
    \centering
  \includegraphics[width=0.7\textwidth]{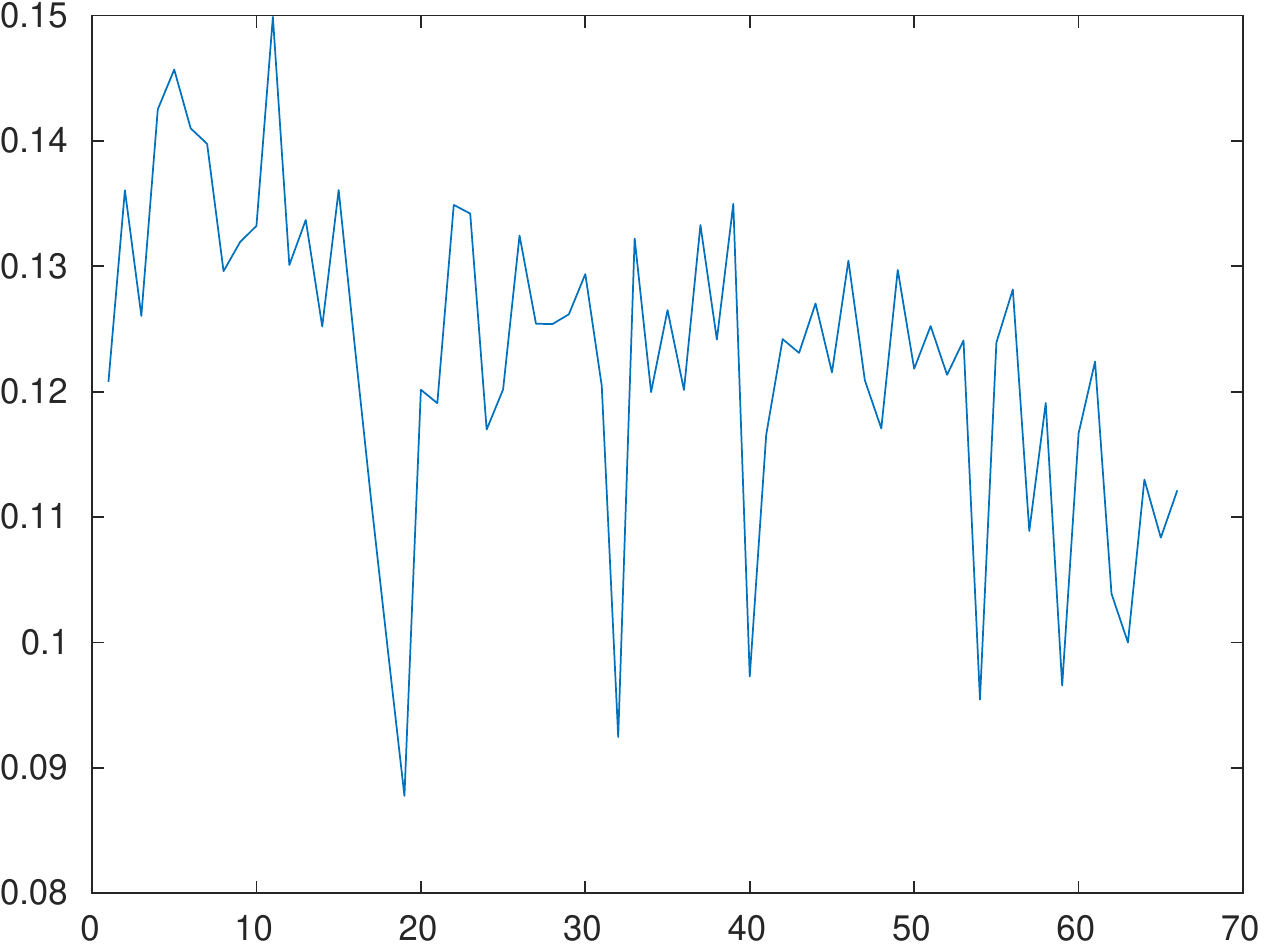}
\end{minipage}
       \begin{minipage}{1\textwidth}
    \footnotesize
    \begin{tabular}{llllll}
      1.    powell  &
    2.    rumsfeld& 
    3.    sec     
    \end{tabular}
       \end{minipage}
\caption{Third term.}
\label{fig:Reuters-3}
\end{figure}

\begin{figure}[htbp!]    
 \centering
  \begin{minipage}{0.7\textwidth}
    \centering
  \includegraphics[width=1\textwidth]{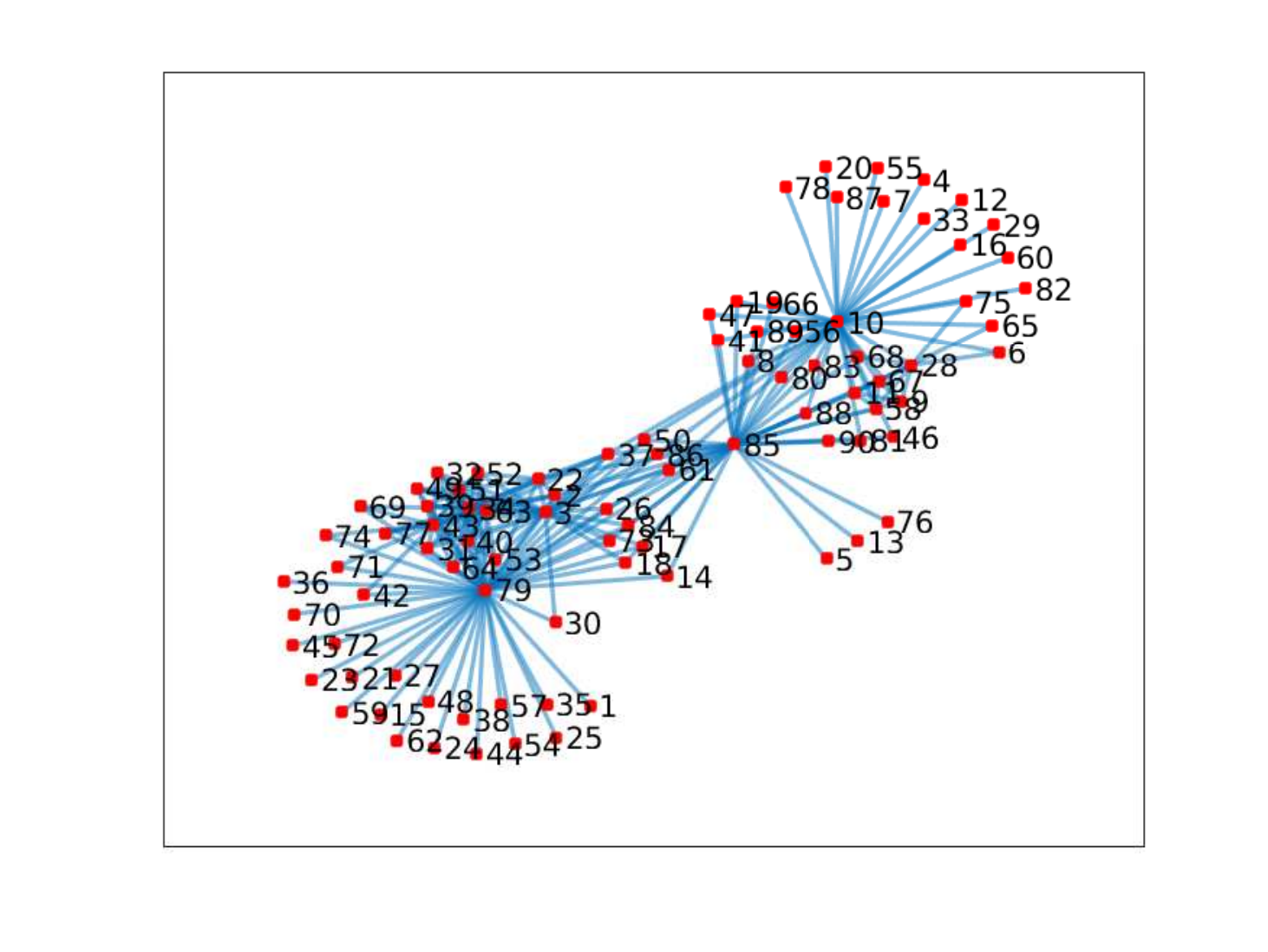}
  \end{minipage}%
  \begin{minipage}{0.3\textwidth}
    \centering
  \includegraphics[width=0.7\textwidth]{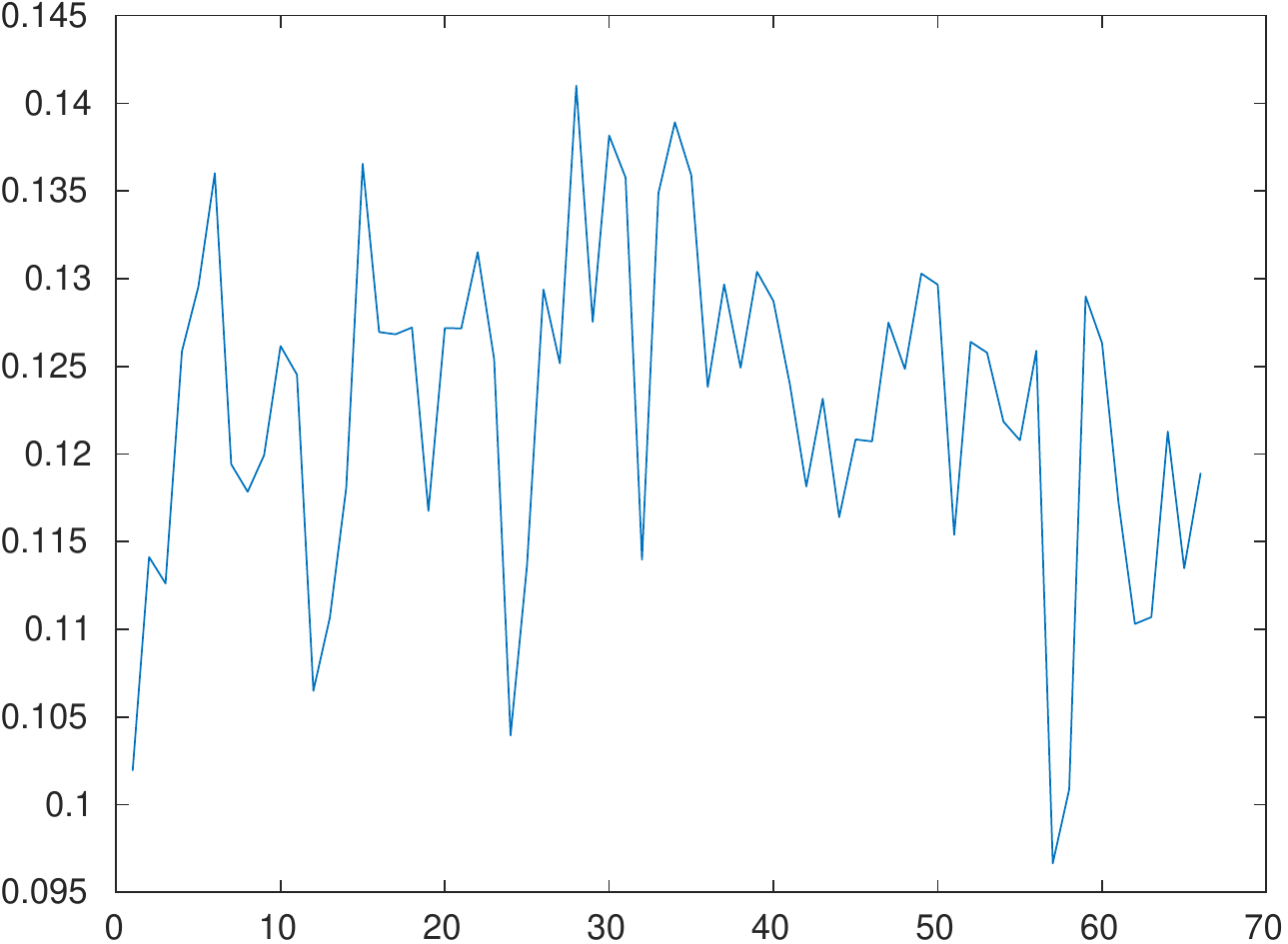}
\end{minipage}
       \begin{minipage}{1\textwidth}
    \footnotesize
    \begin{tabular}{llllll}
1.  abdul\_salam\_zaeef &25. general &49. northern\_alliance&73. strike\\
2.  afghan  &26. government&50. official&74.    stronghold   \\
3.  afghanistan&27. group&51. omar&75. suicide      \\
4.  agent      &28. hijack&52. opposition&76. support      \\
5.  air       &29.  hijacker&53. pakistan&77. supreme      \\
6.  airliner  &30.  islamic &54. pakistani&78. suspect      \\
7.  america   &31.  kabul   &55. passenger&79. taliban      \\
8.     american&32.    kandahar&56.   people&80.    terrorism    \\
9.     anthrax&33.    law&57.    pervez&81.    terrorist\\    
10.    attack&34.    leader&58.    plane&82.    threat       \\
11.    bin\_laden&35.    line&59.    position&83.    time         \\
12.    biological&36.    mazar-i-sharif&60.    possible&84.    troop  \\
13.    bomb&37.    military&61.    pres\_bush&85.    united\_states\\
14.    campaign&38.    minister&62.    president&86.    war\\          
15.    capital&39.    mohammad&63.    rule&87.    warfare      \\
16.    case&40.    mohammad\_omar&64.    ruler&88.    washington   \\
17.    city&41.    month&65.    saudi-born&89.    week         \\
18.    country&42.    movement&66.    security&90.    world        \\
19.    day&43.    mullah&67.    sept\\
20.    deadly&44.    musharraf&68.    september\\
21.    fighter&45.    muslim&69.    southern\\
22.    force&46.    new&70.    special\\
23.    foreign&47.    new\_york&71.    spiritual\\
24.    front&48.    northern&72.    spokesman
    \end{tabular}
       \end{minipage}
\caption{Fourth term.}
\label{fig:Reuters-4}
\end{figure}

\begin{figure}[htbp!]    
 \centering
  \begin{minipage}{0.7\textwidth}
    \centering
  \includegraphics[width=1\textwidth]{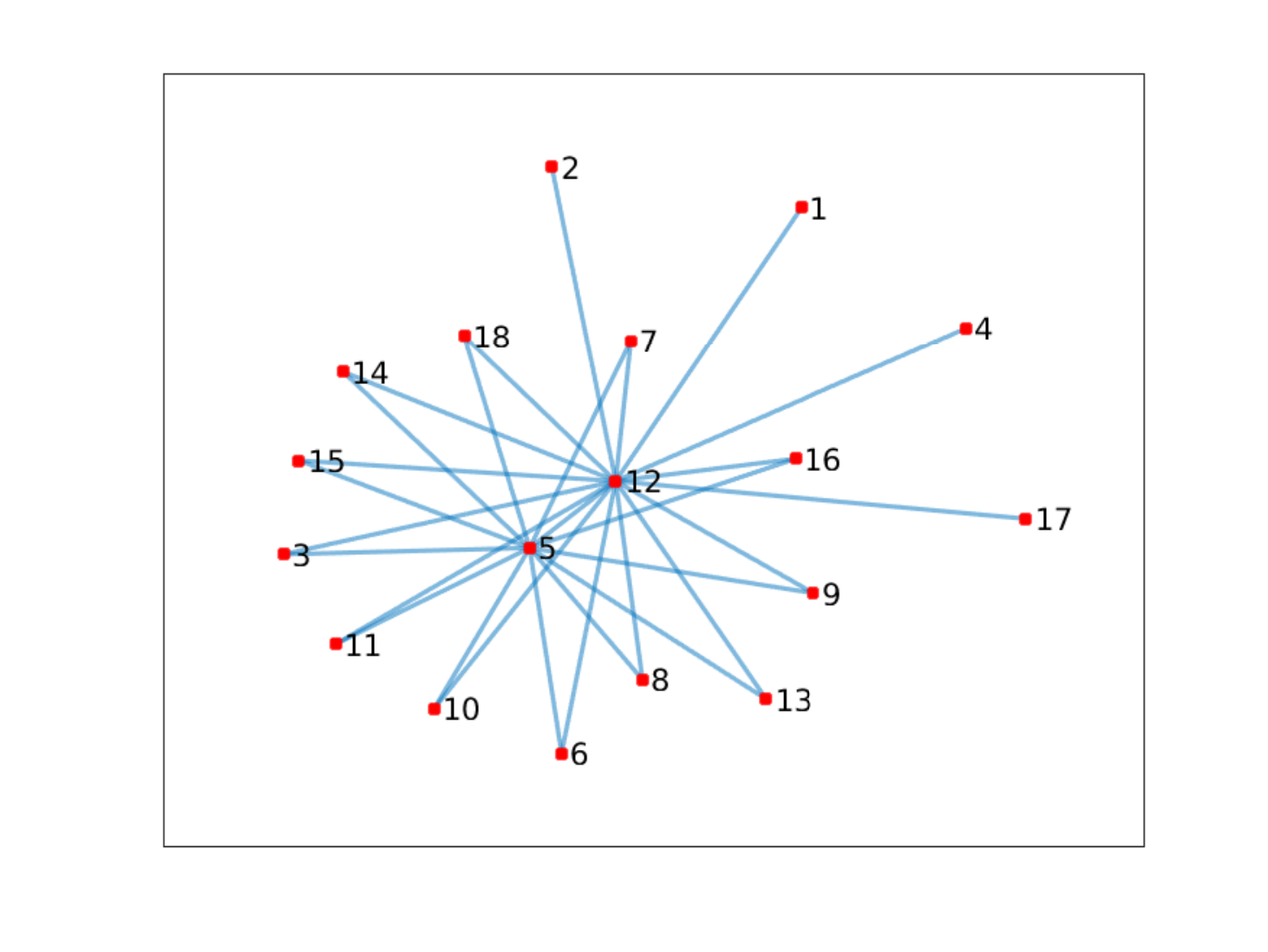}
  \end{minipage}%
  \begin{minipage}{0.3\textwidth}
    \centering
  \includegraphics[width=0.7\textwidth]{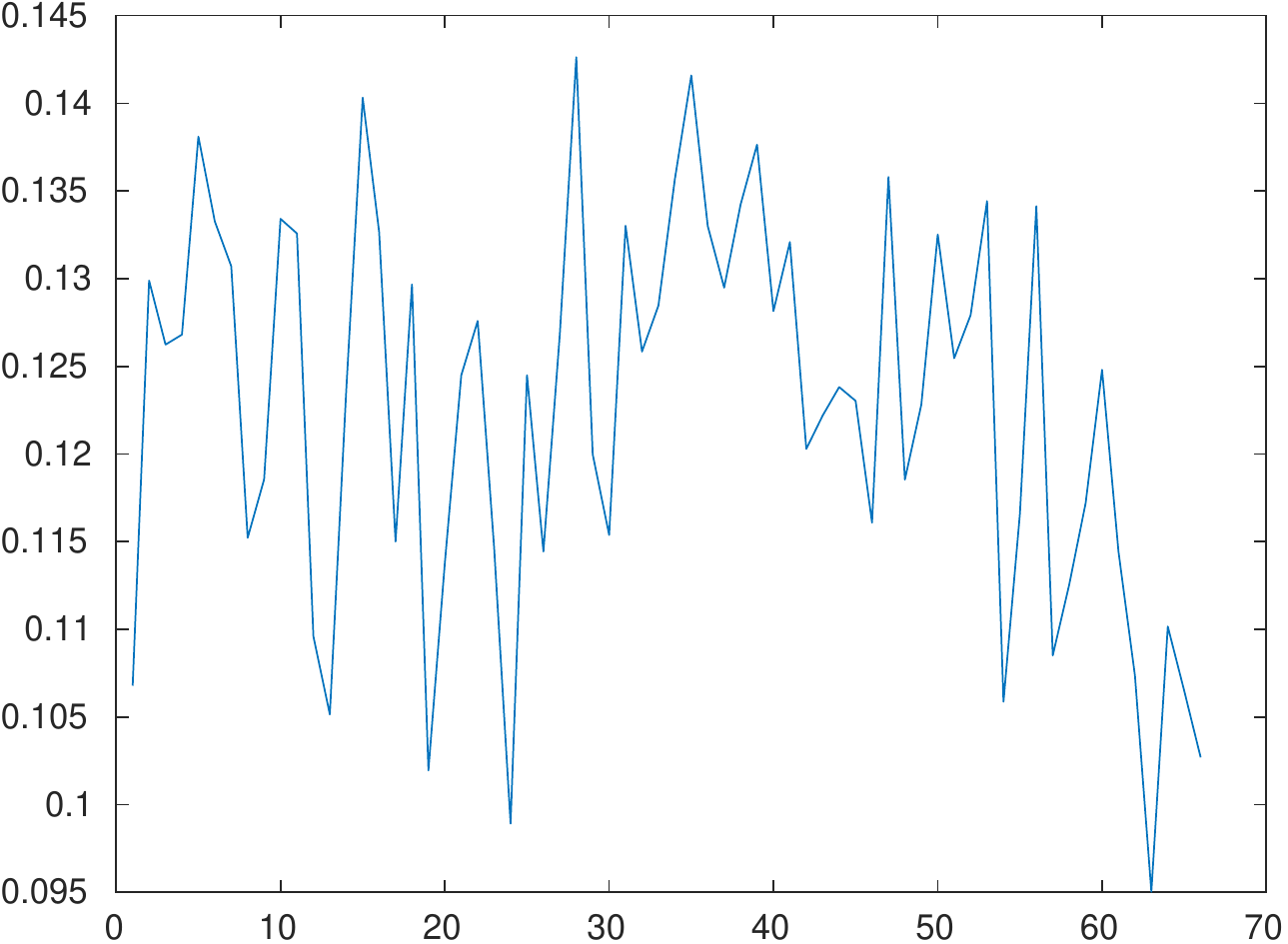}
\end{minipage}
       \begin{minipage}{1\textwidth}
    \footnotesize
    \begin{tabular}{llllll}
1.    agency&6.     enforcement&11.    new\_york&16.    washington     \\
      2.    anchor&7.     law&12.    news&17.    white\_house    \\
3.    brief&8.     mayor&13.    official&18.    world\_trade\_ctr\\
4.    brokaw&9.     mayor\_giuliani&14.    pentagon\\
5.    conference&10.    nbc&15.    tom
    \end{tabular}
       \end{minipage}
\caption{Fifth term.}
\label{fig:Reuters-5}
\end{figure}

\begin{figure}[htbp!]    
 \centering
  \begin{minipage}{0.7\textwidth}
    \centering
  \includegraphics[width=1\textwidth]{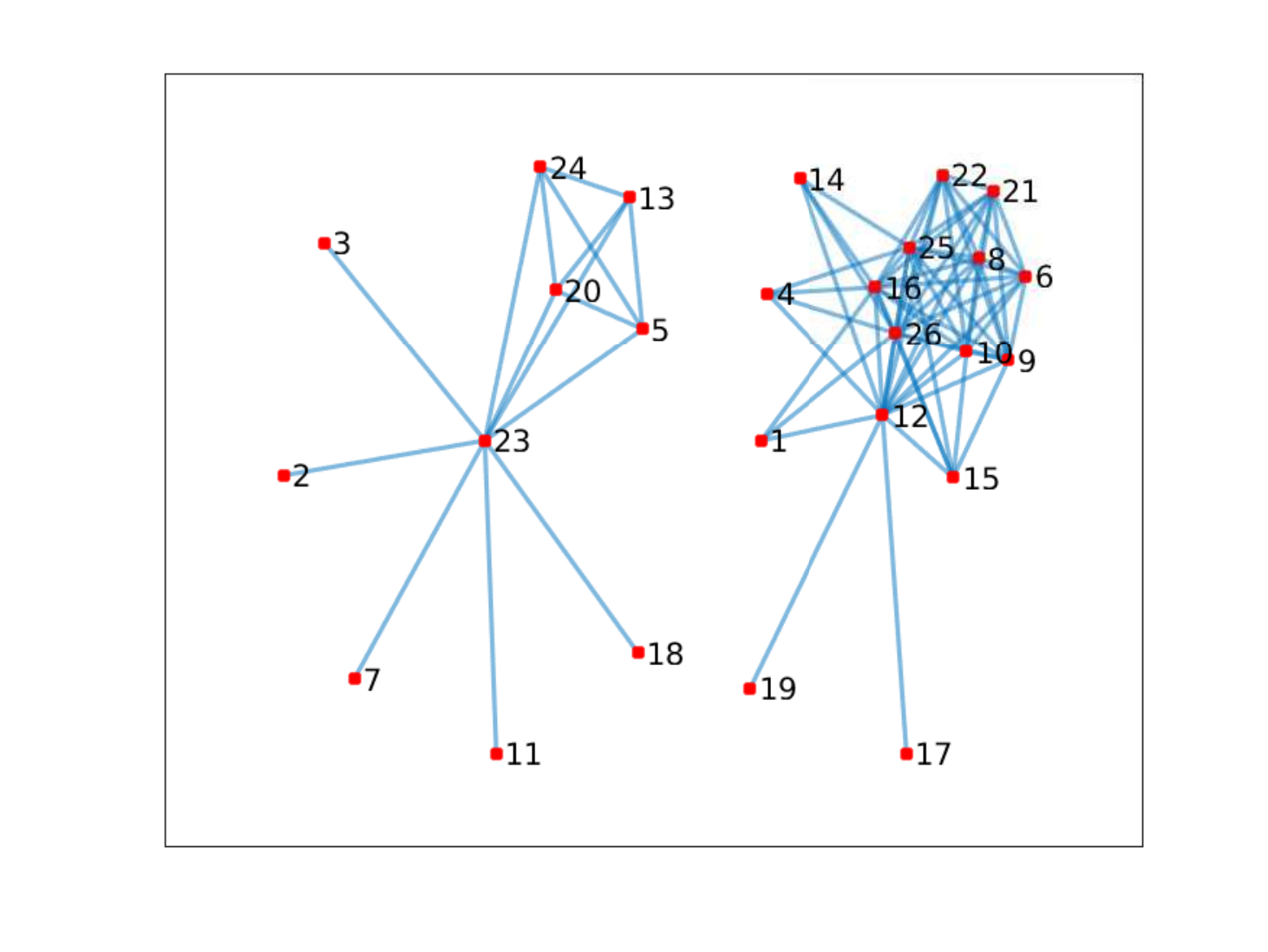}
  \end{minipage}%
  \begin{minipage}{0.3\textwidth}
    \centering
  \includegraphics[width=0.7\textwidth]{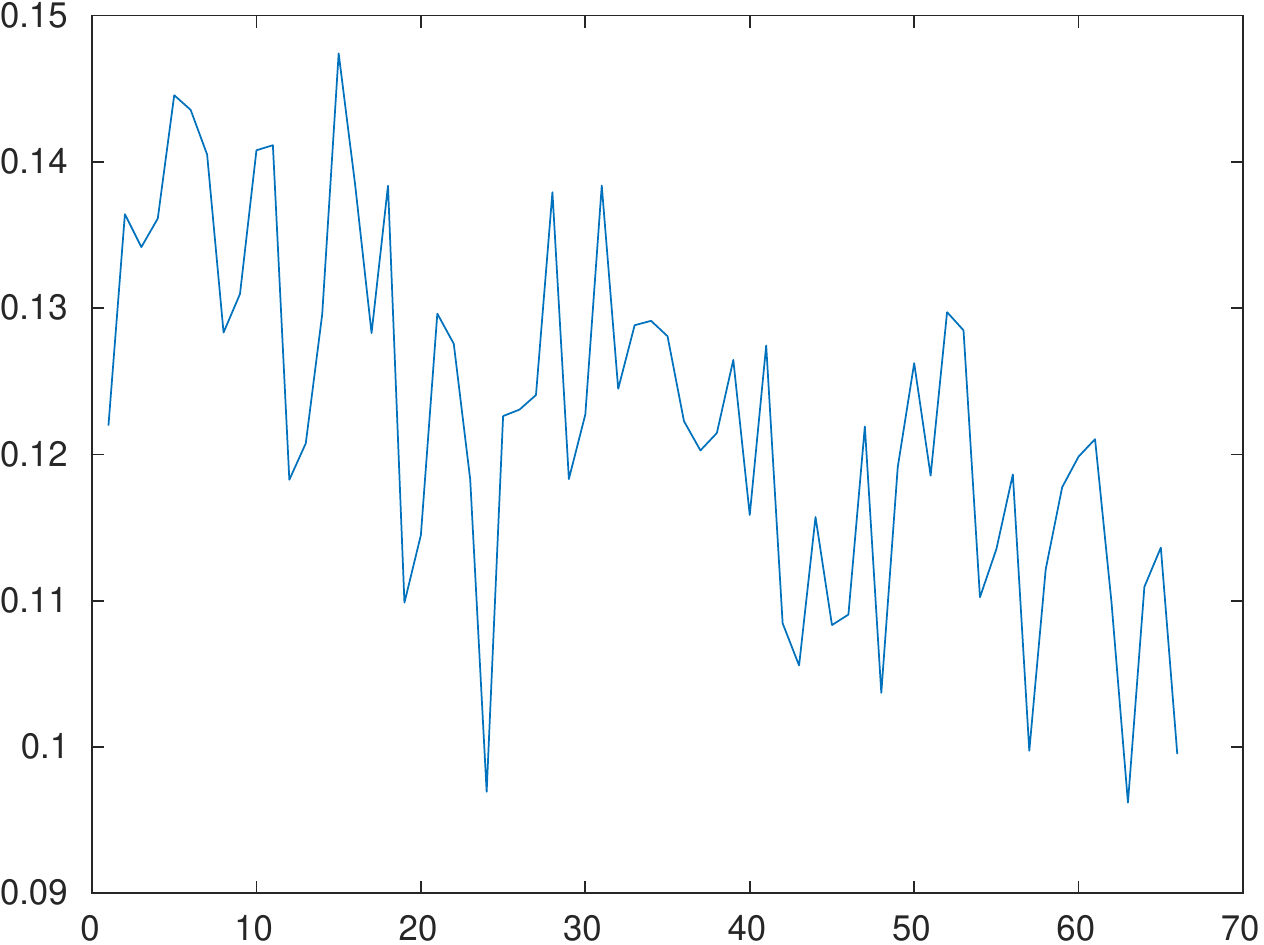}
\end{minipage}
       \begin{minipage}{1\textwidth}
    \footnotesize
    \begin{tabular}{llllll}
1.    110-story&8.     law&15.    pennsylvania&22.    twin           \\
2.    afghanistan&9.     mayor&16.    pentagon&23.    united\_states  \\
3.    attack&10.    mayor\_giuliani&17.    people&24.    war          \\
4.    city&11.    military&18.    pres\_bush&25.    washington     \\
5.    conference&12.    new\_york&19.    tell&26.    world\_trade\_ctr\\
6.    enforcement&13.    news&20.    terrorism\\
7.    force&14.    official&21.   tower
    \end{tabular}
       \end{minipage}
\caption{Sixth term.}
\label{fig:Reuters-6}
\end{figure}

\begin{figure}[htbp!]    
 \centering
  \begin{minipage}{0.7\textwidth}
    \centering
  \includegraphics[width=1\textwidth]{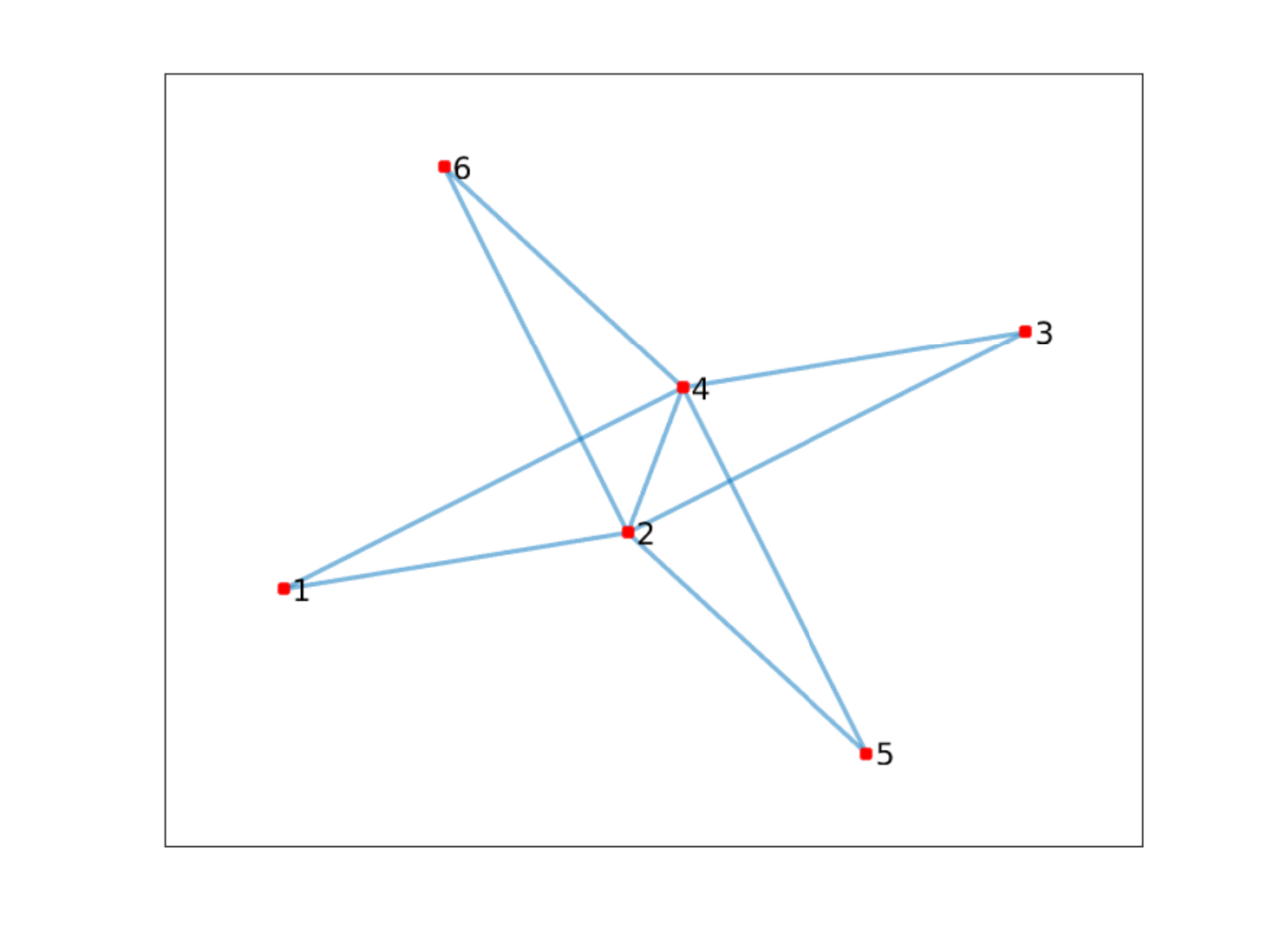}
  \end{minipage}%
  \begin{minipage}{0.3\textwidth}
    \centering
  \includegraphics[width=0.7\textwidth]{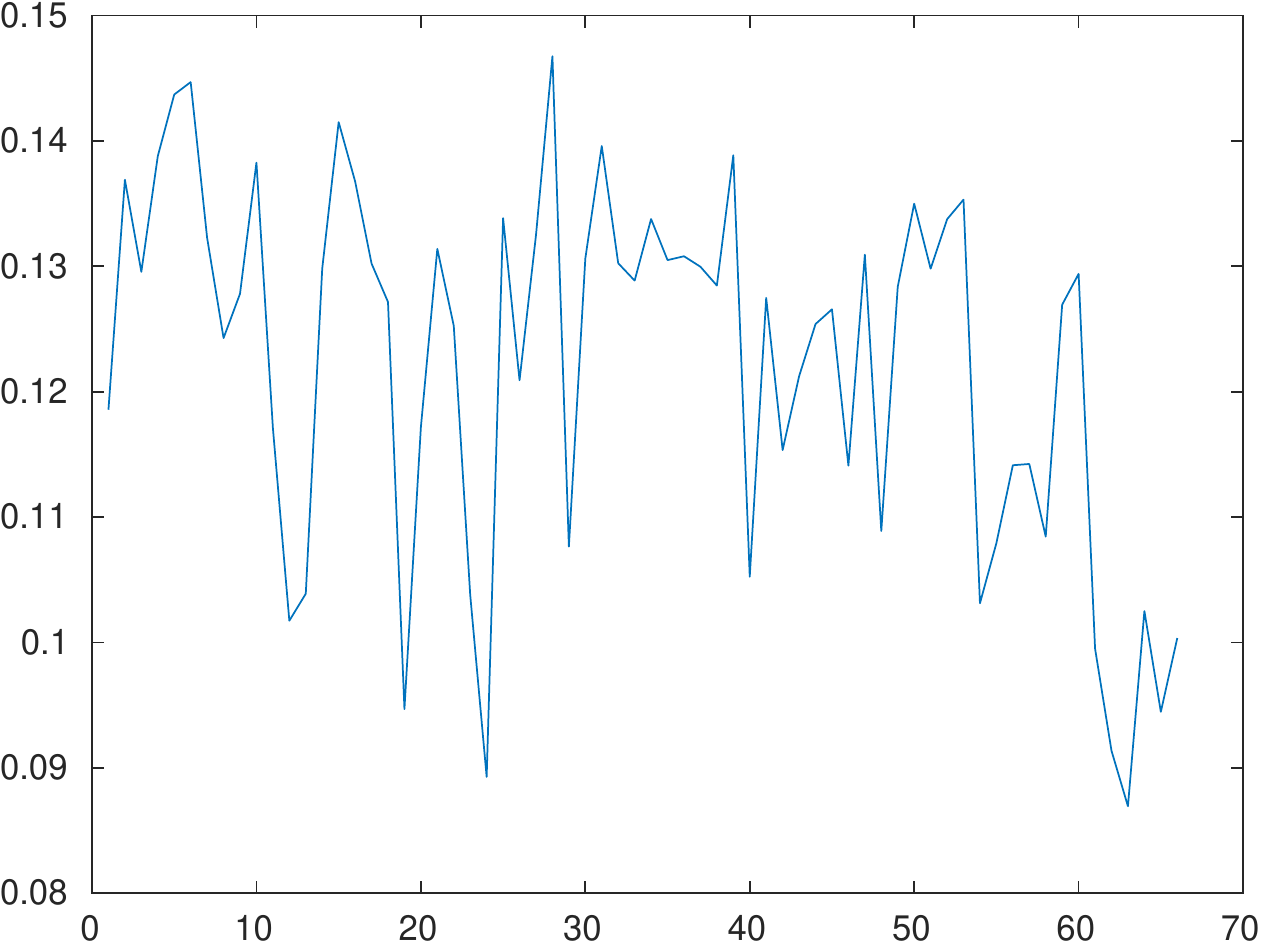}
\end{minipage}
       \begin{minipage}{1\textwidth}
    \footnotesize
    \begin{tabular}{llllll}
 1.    agency&3.    federal&5.    official\\
    2.    enforcement  &4.    law &6.    united\_states
    \end{tabular}
       \end{minipage}
\caption{Seventh term.}
\label{fig:Reuters-7}
\end{figure}

In the construction of the $\widehat{B}^{(\nu)}$ matrices we used the value
$\theta=0.25$  for the cut-off parameter \eqref{eq:b-theta}. The seven
terms are illustrated in Figures
\ref{fig:Reuters-1}--\ref{fig:Reuters-7}. 

In Table \ref{tab:BFnorms-text} we give the norms of the terms in the
expansion, i.e., the norms of the matrices $\widehat{B}^{(\nu)}$. For comparison
we also give the norms of the tensors $\cF^{(\nu)}$, which are equal
to the norms of the terms before thresholding. It is seen that the
thresholding removes a considerable number of elements from
$\cB^{(\nu)}$. Note that the  
terms are  given approximately in descending order of magnitude.
\begin{table}[htbp!]
  \centering
  \caption{Top: Norms of $\widehat{B}^{(\nu)}$, $\cF^{(\nu)}$,  the largest and
    smallest values of ${B}^{(\nu)}$ (i.e. before thresholding), and the
    eigenvalues of the $\cF^{(\nu)}$. Bottom:
    Cosines of the angles between the matrices $\widehat B^{(\nu)} $. 
    \label{tab:BFnorms-text}}
\begin{tabular}{l|ccccccc}
  \hline
  \noalign{\vskip 2pt}
$\| \widehat{B}^{(\nu)} \|$&2.7&1.9&1.8&1.9&1.2&1.6&1.0\\
    $\| \cF^{(\nu)} \|$&6.1&5.6&5.4&5.4&5.1&5.1&5.0\\
    max& 0.46&0.84&1.0&0.24&0.43&0.27&0.54\\
  min&$-0.24$&$-0.11$&$-0.37$&$-0.12$&$-0.13$&$-0.25$&$-0.17$\\
    $\lambda_1$&5.1&4.7&4.7&4.7&4.4&4.4&4.4\\
  $\lambda_2$&3.4&2.9&2.8&2.6&2.5&2.5&2.4\\
 \hline
\end{tabular}
\medskip
\[
     C=
  \begin{bmatrix}
    1&0&0&0.14&0&0.21&0\\
     &1&0&0&0&0&0\\
     & &1&0&0&0&0\\
     & & &1&0&0.04&0\\
     & & & &1&0.11&0\\
     & & & & &1&0.06\\    
     & & & & & &1\\         
  \end{bmatrix}
\qquad C_{\mu\lambda}= \cos(\widehat B^{(\nu)},\widehat B^{(\lambda)}).
\]
\end{table}
To check how much the different terms in the expansion
overlapped, we computed the cosine of the angles between the
matrices $\widehat B^{(\nu)}$ (the cosine is a measure of overlap between edges
in the graphs). 
There was some overlap, as seen in Table 
\ref{tab:BFnorms-text}, but mostly not significant. The vectors $w^{(i)}$, on
the other hand, are  nonnegative and some of them are almost linearly
dependent.   

For comparison we also ran a test, where we used a threshold on the
absolute values in \eqref{eq:b-theta}, 
with $\theta=0.25$ as before. As expected, 
the subgraphs were larger (words corresponding to negative elements in
the ${B}^{(\nu)}$ were retained), but otherwise
the expansion was quite similar.

\subsection{Summary of the Algorithm and Computational Issues}
\label{sec:sum-comp}

The low rank expansion algorithm is summarized below. The 3-slices of initial
tensor $\cA=\cR^{(1)}$ are normalized.

\vbox{
  \vspace*{10pt}
\hrule\medskip
\noindent{\bf   Indata}: $\cA=\cR^{(1)}$, $\theta$, $q$, \quad
\textbf{Outdata:} Expansion $\cA \approx \sum_{\nu=1}^q
\tml[3]{w^{(\nu)}}{\widehat{B}^{(\nu)}}$ 
\smallskip\hrule
\smallskip
\begin{description}\setlength\itemsep{5pt}
  \item {\textbf{for }$\nu=1:q$ }
    \begin{description}\setlength\itemsep{5pt}
      \item Compute best rank-(2,2,1) approximation of $\cR^{(\nu)}$,\\
        giving  $(U^{(\nu)},U^{(\nu)},w^{(\nu)})$, and $\cF^{(\nu)}$ 
      \item Compute ${B}^{(\nu)} =
        \tml[1,2]{U^{(\nu)},U^{(\nu)}}{\cF^{(\nu)}}$
              \item Apply threshold \eqref{eq:b-theta}, giving
                $\widehat{B}^{(\nu)}$ 
       \item $\cR^{(\nu+1)} = \cR^{(\nu)} - \tml[3]{w^{(\nu)}}{\widehat{B}^{(\nu)}}$ 
       \end{description}
   \item \textbf{end}
   \end{description}
   \medskip
\hrule\bigskip}

In the  example in Section \ref{sec:network-logs} the sparse tensor
$\cA$ requires approximately 6 
megabytes storage. The way we construct the matrices $\widehat B^{(\nu)}$ they
are very sparse: between 4 and 77 kilobytes storage.  The deflation producing 
 $\cR^{(2)},\cR^{(3)},\ldots$ incurs some fill-in so that, e.g.,
 $\cR^{(3)}$ requires 
 56 megabytes. However, it is not necessary to perform the deflation
 explicitly. In the Krylov-Schur method only the action of the tensor
 on narrow blocks of vectors is required. Therefore one can keep the
 deflated tensor in the form of the original tensor minus the low rank
 terms, e.g.,
 \[
    \cR^{(3)} = \cA - \sum_{\nu=1}^2
    \tml[3]{w^{(\nu)}}{\widehat{B}^{(\nu)}}.
  \]
  The action of $\cB^{(3)}$ on a block of vectors is then only
  marginally more expensive  than the action of $\cR^{(1)}$.  

\section{Conclusions}
\label{sec:conclusions}

In several applications there is a need to analyze data organized in
tensors.  The aim of this paper has been to show that the best
rank-(2,2,2) and rank-(2,2,1) approximations can be used to extract
useful information from large and sparse 3-tensors from a few
applications. The methods can be considered as generalizations of
spectral graph partitioning.

In the first text analysis example the
tensor was (1,2)-symmetric, and consisted of a sequence of
cooccurrence graphs for terms in news text.  The aim was to find  the main
topics. The tensor method detected structure 
that was not visible when a corresponding  spectral
matrix method was used (probably because in the latter the abundance
of insignificant terms obscured the relevant information). 

In the second  example the tensor was non-symmetric, and
represented  authors and terms at conferences over 17 years. We
demonstrated that the tensor spectral 
method could simultaneously reveal structure in all three modes.

We then presented an algorithm for computing a rank-(2,2,1) expansion
of a (1,2)-symmetric tensor. It was applied to analyze  a  tensor that
represented network traffic logs.   The dominating communication was
found, and it was shown that a small number of users accounted for
most of the communication. It was shown   that the tensor was much
sparser than could be seen initially.  In the second text analysis
example we used the same news texts as in the first. The rank-(2,2,1)
expansion found dominating subgraphs.

In this paper it has not been our intention to give a
comprehensive treatment of the use of best rank-$(2,2,1)$ and
rank-$(2,2,2)$   
approximations to tensor data from  applications. Instead, we have
demonstrated a few problems with real data, where our methods were able to find
relevant structure.  These examples  show that the proposed
methods are potentially  powerful tools  for the analysis of large and
sparse data with tensor structure. Further investigations into the use of our
methods are needed.

  \bibliography{/media/lars/ExtHard/WORK/forskning/BIBLIOGRAPHIES/general,/media/lars/ExtHard/WORK/forskning/BIBLIOGRAPHIES/LE-papers}
\bibliographystyle{plain}


\end{document}